\documentclass[11pt,leqno]{article}

\usepackage{amsmath,amsfonts,amscd,amssymb,theorem}

\long\def\comment#1\endcomment{}

\makeatletter
\begingroup
\gdef\th@dotted{\normalfont\itshape
  \def\@begintheorem##1##2{%
        \item[\hskip\labelsep \theorem@headerfont ##1\ ##2.]}%
\def\@opargbegintheorem##1##2##3{%
   \item[\hskip\labelsep \theorem@headerfont ##1\ ##2\ (##3).]}}
\endgroup
\makeatother

\theoremstyle{dotted}

\newtheorem{theorem}{Theorem}[section]
\newtheorem{lemma}[theorem]{Lemma}

\newtheorem{prop}[theorem]{Proposition}
\newtheorem{corr}[theorem]{Corollary}

\makeatletter
\begingroup
\gdef\th@upshape{\normalfont
  \def\@begintheorem##1##2{%
        \item[\hskip\labelsep \theorem@headerfont ##1\ ##2.]}%
\def\@opargbegintheorem##1##2##3{%
   \item[\hskip\labelsep \theorem@headerfont ##1\ ##2\ (##3).]}}
\endgroup
\makeatother

\theoremstyle{upshape}

\newtheorem{defn}[theorem]{Definition}
\newtheorem{remark}[theorem]{Remark}
\newtheorem{exa}[theorem]{Example}

\makeatletter
\renewcommand{\subsection}{\@startsection{subsection}{2}{0pt}{-3ex
plus -1ex minus -0.2ex}{-2mm plus -0pt minus
-2pt}{\normalfont\bfseries}} 
\renewcommand{\subsubsection}{\@startsection{subsubsection}{3}{0pt}{-3ex
plus -1ex minus -0.2ex}{-2mm plus -0pt minus
-2pt}{\normalfont\bfseries}} 
\makeatother

\makeatletter
\@addtoreset{equation}{section}
\makeatother

\newcommand{\cntrct}                
{\hspace{2pt}\raisebox{1pt}{\text{$\lrcorner$}}\hspace{2pt}}

\newcommand{\proof}[1][Proof.]{\smallskip\noindent{\em #1}}
\def\endproof{\hfill\ensuremath{\square}\par\medskip}

\renewcommand{\labelenumi}{{\normalfont(\roman{enumi})}}

\def\eqref#1{\thetag{\ref{#1}}}
\newcommand{\eqrefi}[2]{\thetag{\ref{#1}\:#2}}

\let\latexref=\ref
\def\ref#1{{\normalfont{\latexref{#1}}}}

\newcommand{\wt}{\widetilde}

\newcommand{\dg}{\dagger}

\newcommand{\idot}{{\:\raisebox{1pt}{\text{\circle*{1.5}}}}}
\newcommand{\eps}{\varepsilon}
\renewcommand{\phi}{\varphi}

\newcommand{\F}{\mathcal{F}}

\newcommand{\Hom}{\operatorname{Hom}}
\newcommand{\Hhom}{\operatorname{\mathcal{H}{\it om}}}

\newcommand{\Fun}{\operatorname{Fun}}
\newcommand{\fFun}{\operatorname{{\mathcal F}{\mathit u}{\mathit n}}}

\newcommand{\Sec}{\operatorname{Sec}}
\newcommand{\sSec}{\operatorname{{\mathcal S}{\mathit e}{\mathit c}}}

\newcommand{\Ch}{\operatorname{Ch}}

\newcommand{\id}{\operatorname{\sf id}}
\newcommand{\Tot}{\natural}

\newcommand{\C}{\mathcal{C}}
\newcommand{\A}{\mathcal{A}}

\newcommand{\hush}{\natural}
\newcommand{\hhush}{{\hush\hush}}

\newcommand{\Unf}{K}
\newcommand{\unf}{k}

\newcommand{\Sets}{\operatorname{Sets}}
\newcommand{\sSets}{\operatorname{\mathcal{S}{\it ets}}}
\newcommand{\Cat}{\operatorname{Cat}}
\newcommand{\cCat}{\operatorname{\mathcal{C}{\it at}}}
\newcommand{\Aut}{{\operatorname{Aut}}}

\newcommand{\Iso}{{\star}}

\newcommand{\ppt}{{\sf pt}}

\newcommand{\copr}{\sqcup}

\newcommand{\N}{{\mathbb N}}

\newcommand{\E}{\mathcal{E}}

\newcommand{\sk}{\operatorname{\text{\rm sk}}}

\newcommand{\Con}{\operatorname{Con}}

\newcommand{\Ho}{\operatorname{Ho}}
\newcommand{\Hho}{\operatorname{\mathcal{H}\mathit{o}}}

\newcommand{\colim}{\operatorname{\sf colim}}
\renewcommand{\lim}{\operatorname{\sf lim}}

\newcommand{\V}{{\sf V}}

\newcommand{\Y}{{\sf Y}}
\newcommand{\ZZ}{{\sf Z}}

\newcommand{\cCone}{\operatorname{\mathcal{C}{\mathit o}{\mathit
      n}{\mathit e}}}

\newcommand{\ev}{\operatorname{\sf ev}}

\newcommand{\Hh}{\mathcal{H}}

\newcommand{\Pos}{\operatorname{PoSets}}
\newcommand{\pPos}{\operatorname{\mathcal{P}{\it o}\mathcal{S}{\it ets}}}
    
\newcommand{\Posf}{\operatorname{Pos}}

\newcommand{\Rel}{\operatorname{BiPoSets}}

\newcommand{\Relf}{\operatorname{BiPos}}

\newcommand{\Cyl}{{\sf C}}

\renewcommand{\P}{\mathcal{P}}

\newcommand{\I}{\mathcal{I}}

\newcommand{\QQ}{{\sf Q}}
\newcommand{\Nn}{{\sf N}}
\newcommand{\cN}{\mathcal{N}}

\newcommand{\bb}{/\!\!/}
\newcommand{\bbd}{\setminus\!\!\!\setminus}
\newcommand{\bbdi}{\bbd_\Iso}
\newcommand{\bbi}{\operatorname{/\!\!/\!_\Iso}}

\newcommand{\mm}{\operatorname{/\!\!/\!_\flat}}
\newcommand{\mc}{\operatorname{/\!\!/\!_\sharp}}
\newcommand{\mmh}{\mm^h}
\newcommand{\mch}{\mc^h}

\newcommand{\dm}{\diamond}

\newcommand{\Ar}{\operatorname{\sf ar}}

\newcommand{\Arc}{\Ar_c}
\newcommand{\Arf}{\Ar_f}

\newcommand{\oB}{\overline{B}}
\newcommand{\oN}{\overline{\N}}

\newcommand{\ssetminus}{\smallsetminus}

\newcommand{\Comp}{\operatorname{Comp}}
\newcommand{\cComp}{\operatorname{\mathcal{C}\mathit{o}\mathit{m}\mathit{p}}}

\newcommand{\Ind}{\operatorname{Ind}}
\newcommand{\iInd}{\operatorname{\mathcal{I}\mathit{n}\mathit{d}}}

\newcommand{\trl}{\triangleleft}
\newcommand{\trr}{\triangleright}

\newcommand{\Env}{\operatorname{Env}}

\def\emptyset{\varnothing}

\title{How to enhance categories, and why?}

\author{D. Kaledin\thanks{Supported by the Russian Science
    Foundation, grant 21-11-00153.}}

\date{{\em To Sonya, for her patience and encouragement.}}

\begin{document}

\maketitle

\tableofcontents

\section*{Introduction.}

While both homology and homotopy originally appeared within
algebraic topology, it became clear quite soon that the former
properly belongs elsewhere. At least since \cite{toho}, we
understand that homological algebra is not at all about homology of
topological spaces; rather, it deals with abelian categories,
derived functors, short exact sequences and suchlike.

For homotopy, the picture even now is much less obvious, and some of
its essential parts are still missing. What seems clear, though, is
that the main subject of homotopical algebra is {\em localization}
-- the act of changing a category by formally inverting a class of
maps.

Philosophically speaking, localization aims at clarifying things by
throwing out irrelevant information. Given a set, we might declare
that some of its elements are indistinguishable by imposing an
equivalence relation and forcing them to become equal. For a
category, this is not a good idea, since objects in a category can
never be equal. Instead, we force certain maps to become
isomorphisms. The term ``localization'' comes from algebraic
geometry where localizing the ring of functions on an affine scheme
corresponds to replacing the scheme with its open subset. In that
context, localization is one of the simplest and best-behaved
operations, but this is only because the rings of functions are
commutative. Once you have a non-commutative algebra, or even simply
a monoid, inverting a set of elements typically creates something
large and hard to control. For categories -- that is, monoids with
several objects -- the situation is even worse. If the category in
question is large, localization may not even exist.

Historically, probably the first general example of a localization
problem that is controllable appeared in the definition of the
derived category given by Verdier in his thesis
\cite{ver.the}. Formally, the derived category of an abelian
category is obtained by localizing the category of chain complexes
with respect to the class of quasiisomorphisms, but Verdier does it
in two steps. First, he passes to the so-called homotopy category of
chain complexes and chain-homotopy classes of maps (this can also be
obtained by localization with respect to the class of chain-homotopy
equivalences, but a direct description is simpler). Then he shows
that the homotopy category admits an additonal structure of a {\em
  triangulated category}, and proves a general localization theorem
for those.

This constrution is not perfect, but to some extent, it
works. However, it is limited in scope: in many respects,
triangulated categories are linear objects (at the very least, they
are additive). In a non-linear situation, things are much more
complicated. The great breakthrough here is Quillen's \cite{qui.ho}
that actually gave a name to the subject, and essentially determined
how we think about it even today. No non-linear analog of a
triangulated structure is known; instead, Quillen controls
localization by introducing the notion of a {\em model category},
another additional structure one can put on a category one wants to
localize. Unlike Verdier localization, this can be only used once,
since the localization of a model category is no longer a model
category -- in fact, axiom $0$ is that a model category is finitely
complete and cocomplete, and localization by its very nature tends
to destroy this. However, if we do have a model category to begin
with, its localization is at least as controllable as the Verdier
localization, and one can use this control to construct adjoint
pairs of functors (this is known as ``Quillen adjunction'') and
mutually inverse equivalences between localizations (``Quillen
equivalence''). In a sense, the latter is the main point of the
exercise: this gives a formal meaning to saying that two categories
``model the same homotopy theory''.

\medskip

The ``homotopy theory'' in question includes the localized category
itself, but there is more. For example, if a category $h(\C)$ is
obtained by localizing a model category $\C$, one can define not
only the sets $\Hom(c,c')$ of morphisms between two objects $c,c'
\in h(\C)$, but a whole homotopy type $\Hhom(c,c')$. The ``naive''
$\Hom(c,c')$ is the set $\pi_0(\Hhom(c,c'))$ of its connected
components, but it can also have non-trivial higher homotopy
groups. Invariants such as these are preserved by Quillen
equivalences, but they cannot be recovered simply from $h(\C)$. Thus
the result of the localization procedure is not just a category; it
should come equipped with an additional structure known as
``enhancement''.

In fact, a triangulated structure in the sense of \cite{ver.the} is
already a form of enhancement, and it allows one to recover
something -- for instance, the higher homotopy groups
$\pi_\idot(\Hhom(c,c'))$ appear as group of maps from $c$ to
homological shifts of $c'$. But this enhancement is much too weak:
triangulated functors do not form a triangulated category, one
cannot do gluing constructions, and so on and so forth (for a longer
but still not exhaustive list of problems, see
e.g.\ \cite[Introduction]{ka.glue}). Moreover, it only makes sense
in the linear setting anyway. For a general notion of an enhancement
one should look elsewhere.

However, a moment's reflection shows that the problem is not easy
from the purely logical point of view. Namely, we say that we want a
``homotopy type'' of morphisms $\Hhom(c,c')$; but what is a homotopy
type? It is either a topological space considered up to a homotopy
equivalence, or a simplicial set considered up to a weak
equivalence, or maybe something else, but still defined ``up to''
something -- that is, an object in a category obtained by
localization. But shouldn't we enhance this ambient category, too?
-- and if yes, in what sense?

\medskip

Among many people bothered by this over the years, the most
prominent was perhaps Grothendieck, with his famous ``thousand pages
long letter'' to Quillen \cite{champ}. Unfortunately, there are no
theorems there, but there are lots of fantastic ideas. One such is
that of ``higher categories''. Roughly speaking, homotopy types that
only have $\pi_0$ and $\pi_1$ are very efficiently parametrized by
small groupoids: objects correspond to points, and morphisms
correspond to homotopy classes of paths. Shouldn't there be a higher
notion of an ``$\infty$-groupoid'', and maybe even
``$\infty$-category'', that not only has objects and morphisms, but
also morphisms between morphisms, morphisms between morphisms
between morphisms, and so on?  Well, maybe, but let us spell out
something that is usually not spelled out explicitly:
\begin{itemize}
\item {\em this idea does not work}.
\end{itemize}
Or at least, nobody has been able to make it work, in the
forty-something years that \cite{champ} has been around. The problem
is the ``and so on'' part: there is no ``so on''. It is not possible
to describe the resulting full structure of higher compatibilities
and constraints.

In retrospect, this is not too surprising, since doing this would in
particular give a purely combinatorial effective description of
homotopy groups of spheres, and it is safe to assume that this
cannot be done (although those with a surfeit of spare time and
ambition are welcome to try). In practice, all the notions of
``$\infty$-groupoids'' that appeared in the literature were still,
at the end of the day, defined up to a weak equivalence of some
sort, and then it was not clear why bother -- fibrant simplicial
sets of Kan have been around since 1950-ies. Effectively, in modern
usage, ``$\infty$-groupoid'' is simply a fancy name for a fibrant
simplicial set. It is doubtful that such rebranding helps one to
prove something.

\medskip

In any case, as far as enhancements are concerned, one option would
be to just declare that an enhanced category is the same thing as a
model category ``considered up to a Quillen equivalence'', but this
is not enough since not all localizations needed in practice appear
in this way (essentially because of axiom $0$ mentioned above). Thus
the current thinking goes along more-or-less the following lines.
\begin{enumerate}
\item ``Quillen-equivalent model categories have the same homotopy
  theory''; this is accepted as an article of faith and not
  discussed.
\item One constructs a ``category of models'' for enhanced small
  categories; this category of models is equipped with a model
  structure and produces all the desiderata; an ``enhanced
  category'' is then simply defined as an object in the
  corresponding localized category.
\item Models are not unique at all, and neither are ``categories of
  models'', but one checks that they are all Quillen-equivalent, so
  see \thetag{i}.
\end{enumerate}
There are two obvious issues with this kind of thinking. Firstly, it
is very set-theoretical in nature and feels like a throughback to
19-th century -- a category, something that should be a fundamental
notion, is treated as a special type of a simplicial set, or
``space'', whatever it is, or something like that. The idea of {\em
  symmetry} so dear to people like Grothendieck is thrown out of the
window. Secondly, a worse problem is the inherent circularity of the
argument. Of all the avaliable models, it is best seen in the
approach of \cite{BK} based on relative categories. By definition, a
{\em relative category} is a small category $\C$ equipped with a
class of maps $W$. Barwick and Kan propose putting a model structure
on the category of relative categories, and showing that it is
Quillen-equivalent to all the other existing models. Then in this
particular model, the result of localizing a category $\C$ with
respect to a class of maps $W$ is the relative category $\langle
\C,W \rangle$. Effectively, it looks pretty much as if in this
approach -- and ipso facto in all the others, since they are all
Quillen-equivalent -- one solves the localization problem by
declaring it solved.

\medskip

Of course, all of the above might be the nature of things: if one
cannot do better, it is better to do something rather than
nothing. The situation might just be similar to the definition of an
algorithmically computable function, where we have several
alternative definitions, none of them too natural, but they are all
provably equivalent, so Church's Thesis just declares that this is
it. The point of the present paper, however, is that one {\em can}
do better. This is based on another fantastic idea of \cite{champ},
and this idea indeed works.

\medskip

Roughly speaking, what one wants to do is the following. If it is
inevitable that enhanced categories are only defined up to an
equivalence of some sort, let us at least make this equivalence as
easy to control as possible. Then observe that there is another type
of controlled localization that is so common and widespread that it
usually goes unnoticed by its users: the category $\Cat$ of small
categories, and the class $W$ of, well, equivalences of
categories. In principle, this can be localized by using model
category techniques, but this is akin to smelling roses through a
gas mask. The answer is actually much simpler, and similar to the
homotopy category of chain complexes: objects are small categories,
morphisms are isomorphism classes of functors. Moreover, we can also
consider families of small categories indexed by some category
$I$. This is conveniently packaged by the Grothendieck construction
of \cite{SGA} into a {\em Grothendieck fibration} $\C \to I$ with
small fibers, with morphisms between fibrations given by functors
$\C \to \C'$ cartesian over $I$. Then again, localizing the category
of fibrations with respect to equivalences gives the category with
the same objects, and isomorphism classes of cartesian functors as
morphisms (for precise definitions, see below
Subsection~\ref{yo.subs}).

Now, whatever an enhanced category $\C$ is, it should come equipped
with its underlying usual category $h(\C)$, but there is more: for
any small category $I$, we should also have the enhanced category
$\C^I$ of functors $I \to \C$, and its underlying usual category
$h(\C^I)$. Thus we actually have a whole family of categories
indexed by $\Cat$. This has been described in \cite{champ} under the
name of a {\em derivator}; the question was, is it enough to recover
$\C$?  Our answer is: with some modifications, yes.

\medskip

The notion of a derivator has been studied by many people over the
years, and seems to have become a standard term with a well-defined
meaning (see e.g.\ \cite{grotz} for an overview); to avoid
gratuitous rebranding, let us just call our objects {\em enhanced
  categories}. The main modification compared to \cite{champ} is
that it is not necessary, nor in fact desirable to index our
enhanced categories over the whole $\Cat$ -- it is sufficient to
consider the category $\Pos$ of partially ordered sets. We actually
use an even smaller category $\Posf^+$ of {\em left-bounded}
partially ordered sets, see Definition~\ref{dim.def} and the
paragraph below it, but any enhanced category defined over $\Posf^+$
extends canonically to a family over $\Pos$. What seems important
here is that partially ordered sets considered as categories are
{\em rigid}, that is, the only isomorphisms are the identity
maps. Therefore $\Pos$ has no $2$-categorical structure, and
treating it as simply a category is a reasonable thing to do. One
thing that seems {\em not} possible to do is to cut down even
further and index our enhanced categories on finite partially
ordered sets. As things stand, $\Posf^+$ is a large category, so one
has to check that our enhanced categories only admit a set of
isomorphism classes of functors between them (they do).

The first thing we prove is comparison with the standard theories of
enhancements. Our basis for comparison is given by complete Segal
spaces of Ch. Rezk \cite{rzk}, mostly because they are the closest
to intuition (quasicategories of Joyal and Lurie \cite{lu1} allow
for simpler proofs, but for our purposes this is immaterial). We
define a comparison functor from the homotopy category of complete
Segal spaces to the category of fibrations $\C \to \Posf^+$ with
small fibers, and isomorphism classes of cartesian functors between
them, and we prove that this comparison functor is fully
faithful. We then describe its essential image by imposing several
axioms on a fibration -- roughly five in total (or six if one wants
to look at fibrations over the whole $\Pos$ rather than just
$\Posf^+$). This defines an enhanced category. After that, we use
standard category theory techniques to bootstrap the whole theory --
an enhancement for the category of small enhanced categories,
enhanced categories of enhanced functors, enhanced limits and
colimits, the Yoneda embedding, the enhanced version of the
Grothendieck construction, etc.

It is perhaps instructive to mention how our theory describes
homotopy types. These correspond to enhanced groupoids, that is,
enhanced categories given by fibrations $\C \to \Posf^+$ whose
fibers are groupoids. In a sense, the whole gadget exhibits a sort
of an Eckmann-Hilton duality between the ideas of order (exemplified
by partially ordered sets $J \in \Posf^+$) and symmetry (exemplified
by the groupoids $\C_J$). In another sense, it restores the original
idea of ``symmetries between symmetries'', but in different
guise. There are no ``higher groupoids'', there are just groupoids
in the usual sense -- but a whole bunch of them (just as a scheme
can be thought of as a bunch of sets of its points over various
affine schemes).

\medskip

One feature of our approach to enhancements that we think is an
improvement over the existing alternatives is that it follows very
closely the usual categorical intuition and way of thinking. Thus
all the definitions, constructions and statements are quite simple
and natural, and can be explained in a reasonably concise way. The
proofs can not. Therefore on some reflection, we have decided to
split the exposition into two parts. The long and technical
\cite{big} contains all the proofs, and is completely independent of
the present paper. The present paper is an overview with all the
constructions and definitions, but no proofs at all; those are
replaced by precise references to \cite{big}. The idea is to give a
toolkit that is sufficient for practical applications, with all the
proofs safely hidden in the black box.

The paper is organized as follows. Section~\ref{cat.subs} contains a
very brief overview of standard category theory. This is needed
because there are no convenient textbooks, but it also serves a dual
purpose. Firstly, we fix notation and terminology, and spell out
precisely things that need to be spelled out precisely. Secondly, we
give a template that we will then repeat in the enhanced
setting. The short Section~\ref{pos.sec} contains the necessary
technical preliminaries about partially ordered sets, Segal spaces
and so on; this is kept to an absolute
minimum. Section~\ref{enh.sec} contains the main definitions and
structural results, including the comparison theorem mentioned
above, and Section~\ref{enh.cat.sec} describes the more advanced
parts of the theory (Grothendieck construction, the Yoneda package,
limits and Kan extensions). Finally, Section~\ref{app.sec} is an
appendix where we explain, to an interested reader, the main ideas
behind the proofs. Needless to say, it can be safely skipped.

\subsection*{Acknowledgements.}

This paper together with \cite{big} is the result of a project of
some duration; it is hopeless to try to mention here all the
colleagues who helped me along the way. I try to give at least a
partial list in \cite[Introduction]{big}. Here, let me just go back
to the beginning and express my deep gratitude to Vladimir Voevodsky
and to my Ph.D. advisor David Kazhdan who introduced me to
polycategories back in 1991, and generously tolerated my childish
attempts to dabble in the topic. Vladimir later on moved to other
things, but David remained consistently interested in the subject
through all these years, and I owe much of my understanding to
continued conversations with him. The same goes for Vladimir Hinich
and Sasha Beilinson. I am also very grateful to the referee for the
meticulous work on the manuscript and many suggested improvements.

\section{Category theory.}\label{cat.sec}

\subsection{Categories and functors.}\label{cat.subs}

We work in a minimal set-theoretic setup limited to small and large
categories; we do not assume the universe axiom. Even for a large
category, $\Hom$-sets are small. We use ``map'', ``morphism'' and
``arrow'' interchangeably. For any category $\C$, we write $c \in
\C$ as shorthand for ``$c$ is an object of $\C$''. For any category
$\C$, we denote by $\C^o$ the opposite category, and for any functor
$\gamma:\C_0 \to \C_1$, we denote by $\gamma^o:\C_0^o \to \C_1^o$
the opposite functor. A morphism $f:c \to c'$ defines a morphism
from $c'$ to $c$ in $\C^o$ that we denote by $f^o:c' \to c$. We
denote by $\ppt$ the point category (a single object, a single
morphism), and we denote by $[1]$ the ``single arrow category'' with
two objects $0,1 \in [1]$ and a single non-identity arrow $0 \to
1$. For any category $\C$ and object $c \in \C$, we denote by
$\eps(c):\ppt \to \C$ the functor with value $c$.  Morphisms in a
category $\C$ correspond to functors $[1] \to \C$, and commutative
squares in $\C$ correspond to functors $[1]^2 \to \C$. For any
categories $\C$, $\C'$, giving a functor $\gamma:\C \times [1] \to
\C'$ is equivalent to giving functors $\gamma_0,\gamma_1:\C \to \C'$
and a map $\gamma_0 \to \gamma_1$.

A {\em projector} in a category $\C$ is an endomorphism $p:c \to c$
of an object $c \in \C$ such that $p^2=p$; an {\em image} of a
projector $p:c \to c$ is an object $c' \in \C$ equipped with maps
$a:c' \to c$, $b:c \to c'$ such that $p=a \circ b$ and $b \circ a =
\id$. An image is unique up to a unique isomorphism if it exists, so
it's properly called ``the image'', and it is automatically
preserved by any functor $\C \to \C'$. A category is {\em
  Karoubi-closed} if all projectors have images. The {\em Karoubi
  completion} $\C'$ of a category $\C$ is the category of pairs
$\langle c,p \rangle$, $c \in \C$, $p:c \to c$ a projector, with
morphisms $\langle c,p \rangle \to \langle c',p' \rangle$ given by
morphisms $f:c \to c'$ such that $p' \circ f = f = f \circ p$. Then
$\C'$ is Karoubi-closed, we have a functor $\eps:\C \to \C'$, $c
\mapsto \langle c,\id \rangle$, and any functor $\gamma:\C \to \E$
to a Karoubi-closed category $\E$ factors through $\eps$, uniquely
up to a unique isomorphism.

A functor $\gamma:\C \to \C'$ is {\em full} resp.\ {\em faithful} if
it surjective resp.\ injective on $\Hom$-sets, and {\em essentially
  surjective} if it is surjective on isomorphism classes of
objects. A functor $\gamma$ is {\em conservative} if any map $f$
with invertible $\gamma(f)$ is itself invertible. A functor is an
{\em equivalence} if it is fully faithful and essentially surjective
--- or equivalently, if it is invertible up to an isomorphism ---
and we say that a functor is an {\em epivalence} if it is
essentially surjective, full and conservative (but maybe not
faithful). A category is {\em essentially small} if it is equivalent
to a small category, and a functor $\gamma:\C' \to \C$ is small if
the preimage $\gamma^{-1}(\C_0)$ of any essentially small full
subcategory $\C_0 \subset \C$ is essentially small. A {\em
  commutative square} of categories and functors is a square
\begin{equation}\label{cat.sq}
\begin{CD}
\C_{01} @>{\gamma_{01}^1}>> \C_1\\
@V{\gamma_{01}^0}VV @VV{\gamma_1}V\\
\C_0 @>{\gamma_0}>> \C
\end{CD}
\end{equation}
of categories and functors equipped with an isomorphism $\gamma_1
\circ \gamma_{01}^1 \cong \gamma_0 \circ \gamma_{01}^0$. We denote
by $\C_0 \times_\C \C_1$ the category of triples $\langle
c_0,c_1,\alpha \rangle$ of objects $c_l \in \C_l$, $l=0,1$ and an
isomorphism $\alpha:\gamma_0(c_0) \cong \gamma_1(c_1)$. We write
$\C_0 \times_\C \C_1 = \gamma_0^*\C_1$ when we want to emphasize the
dependence on $\gamma_0$. A square \eqref{cat.sq} induces a functor
\begin{equation}\label{cat.sq.eq}
\C_{01} \to \C_0 \times_\C \C_1,
\end{equation}
and a square is {\em cartesian} resp.\ {\em semicartesian} if
\eqref{cat.sq.eq} is an equivalence resp.\ an epivalence. A
commutative square \eqref{cat.sq} is {\em cocartesian} if for any
category $\C'$, functors $\gamma'_l:\C_l \to \C'$,
$l=0,1$, and an isomorphism $\gamma'_0 \circ \gamma^0_{01} \cong
\gamma_1' \circ \gamma^1_{01}$, there exists a functor $\gamma:\C
\to \C'$ and isomorphisms $\alpha_l:\gamma \circ \gamma_l \cong
\gamma_l'$, $l=0,1$, and the triple $\langle
\gamma,\alpha_0,\alpha_1 \rangle$ is unique up to a unique
isomorphism.

An {\em adjoint pair of functors} between categories $\C$, $\C'$ is
a pair of functors $\lambda:\C' \to \C$, $\rho:\C \to \C'$ equipped
with {\em adjunction maps} $a_\dg:\id \to \rho \circ \lambda$,
$a^\dg:\lambda \circ \rho \to \id$ such that the compositions
\begin{equation}\label{adj.eq}
\begin{CD}
\rho @>{a_\dg \circ \rho}>>  \rho \circ \lambda \circ \rho
@>{\rho(a^\dg)}>> \rho,\\
\lambda @>{\lambda(a_\dg)}>> \lambda \circ \rho \circ \lambda
@>{a^\dg \circ \lambda}>> \lambda
\end{CD}
\end{equation}
are both equal to the identity. The functor $\lambda$ is {\em
  left-adjoint} to $\rho$, and the functor $\rho$ is {\em
  right-adjoint} to $\lambda$. The adjunction maps induce
isomorphisms
\begin{equation}\label{adj.iso}
\Hom_\C(\lambda(c'),c) \cong \Hom_{\C'}(c',\rho(c)), \qquad c \in
\C,c' \in \C',
\end{equation}
functorial in $c$ and $c'$. In any adjoint pair, either one of the
functors $\lambda$, $\rho$ determines the other one and the
adjunction maps uniquely up to a unique isomorphism, so that having
an adjoint is a condition and not a structure. If we already have
$\rho$ and $\lambda$, then either of the maps $a_\dg$, $a^\dg$
uniquely defines the other one, so this is again a condition; if the
other map exists, we will say that $a_\dg$ resp.\ $a^\dg$ {\em
  defines an adjunction} between $\rho$ and $\lambda$. For brevity,
we will say that a functor $\phi:\C \to \C'$ is {\em left}
resp.\ {\em right-reflexive} if it admits a left-adjoint
$\phi_\dg:\C' \to \C$ resp.\ a right-adjoint $\phi^\dg:\C' \to
\C$. A full subcategory $\C' \subset \C$ is {\em left} resp.\ {\em
  right-admissible} if the embedding functor $\gamma:\C' \to \C$ is
left resp.\ right-reflexive. If we have a commutative square
\eqref{cat.sq} such that $\gamma_0$ and $\gamma_{01}^1$ are
right-reflexive, then the isomorphism $\gamma_0 \circ \gamma^0_{01}
\cong \gamma_1 \circ \gamma^1_{01}$ defines by adjunction a map
\begin{equation}\label{bc.eq}
\gamma^1_{01} \circ \gamma^{0\dg}_{01} \to \gamma_0^\dg \circ
\gamma_1,
\end{equation}
and dually for left-reflexive functors. We will call the maps
\eqref{bc.eq} the {\em base change maps}. 

A small category is {\em discrete} if all its maps are identity
maps, so that a discrete small category is the same thing as a
set. We denote by $\Sets$ resp.\ $\Cat$ the categories of sets
resp.\ small categories, and we have the full embedding $\Sets
\subset \Cat$ identifying discrete categories and sets. We denote by
$\iota:\Cat \to \Cat$ the involution $\C \mapsto \C^o$. A category
is {\em rigid} if its only isomorphisms are identity maps.
Conversely, a category is a {\em groupoid} if all its maps are
isomorphisms. Any functor between groupoids is trivially
conservative, so it is an epivalence as soon as it is essentially
surjective and full. A left or right-reflexive functor between
groupoids is an equivalence (because the adjunction maps are
automatically invertible). For any group $G$, we denote by $\ppt_G$
the groupoid with a single object with automorphism group $G$. The
{\em isomorphism groupoid} $\C_\Iso \subset \C$ of a category $\C$
has the same objects as $\C$, and those maps between them that are
invertible.

For any category $\C$ and essentially small category $I$, functors
$I \to \C$ form a category that we denote $\Fun(I,\C)$, and we
simplify notation by writing $I^o\C = \Fun(I^o,\C)$. If $\C$ is
equipped with a functor $\pi:\C \to I$, we define the category
$\Sec(I,\C)$ of
{\em sections} of the functor $\pi$ by the cartesian square
\begin{equation}\label{sec.sq}
\begin{CD}
\Sec(I,\C) @>>> \Fun(I,\C)\\
@VVV @VVV\\
\ppt @>{\eps(\id)}>> \Fun(I,I),
\end{CD}
\end{equation}
where the bottom arrow is the embedding onto $\id:I \to
I$. Explicitly, a section is a functor $s:I \to \C$ equipped with an
isomorphism $\alpha:\pi \circ s \to \id$, and $\Sec(I,\C)$ is the
category of pairs $\langle s,\alpha \rangle$. In general,
$\Fun(I,\C)$ comes equipped with the evaluation pairing $\ev:I
\times \Fun(I,\C) \to \C$, and for any category $I'$, a pairing
$\gamma:I \times I' \to \C$ factors as
\begin{equation}\label{fun.facto}
\begin{CD}
I \times I' @>{\id \times \wt{\gamma}}>> I \times \Fun(I,\C)
@>{\ev}>> \C
\end{CD}
\end{equation}
for a functor $\wt{\gamma}:I' \to \Fun(I,\C)$, unique up to a unique
isomorphism. In particular, any essentially small category $I$ comes
equipped with a $\Hom$-pairing $\Hom_I(-,-):I^o \times I \to \Sets$,
and \eqref{fun.facto} then gives rise to the fully faithful {\em
  Yoneda embedding}
\begin{equation}\label{yo.eq}
\Y:I \to I^o\Sets.
\end{equation}
For any functor $\gamma:I' \to I$ between essentially small $I$,
$I'$, we denote by $\gamma^*:\Fun(I,\C) \to \Fun(I',\C)$ the
pullback functor obtained by precomposition with $\gamma$. If
$I=[1]$, then $\Fun([1],\C) = \Ar(\C)$ is the arrow category of the
category $\C$: its objects are arrows $c_0 \to c_1$ in $\C$, and
morphisms are given by commutative squares. We let
$\sigma,\tau:\Ar(\C) \to \C$ be the functors sending an arrow $c_0
\to c_1$ to its source $c_0$ resp.\ target $c_1$, and we let
$\eta:\C \to \Ar(\C)$ send $c$ to $\id:c \to c$; then $\eta$ is
left-adjoint to $\sigma$ and right-adjoint to $\tau$.

For any category $\C$, we denote by $\C^<$ resp.\ $\C^>$ the
category obtained by adding the new initial resp.\ terminal object
$o$ to $\C$, and we let $\eps:\C \to \C^>$, $\eps:\C \to \C^<$ be
the embeddings, and $o:\C \to \C^>$, $o:\C \to \C^<$ be the constant
functors with value $o$. For any functor $\gamma:\C_0 \to \C_1$, the
{\em cylinder} $\Cyl(\gamma)$ and the {\em dual cylinder}
$\Cyl^o(\gamma)$ are defined by cartesian squares
\begin{equation}\label{cyl.def.sq}
\begin{CD}
\Cyl(\gamma) @>>> \C_0^>\\
@VVV @VV{\gamma^>}V\\
\C_1 \times [1] @>>> \C_1^>,
\end{CD}
\qquad\qquad
\begin{CD}
\Cyl^o(\gamma) @>>> \C_0^<\\
@VVV @VV{\gamma^<}V\\
\C_1 \times [1] @>>> \C_1^<,
\end{CD}
\end{equation}
where the bottom arrows correspond to tautological maps $\eps \to o$
resp.\ $o \to \eps$, and we note that we have cocartesian squares
\begin{equation}\label{cyl.sq}
\begin{CD}
\C_0 @>{\id \times t}>> \C_0 \times [1]\\
@V{\gamma}VV @VVV\\
\C_1 @>>> \Cyl(\gamma),
\end{CD}
\qquad\qquad
\begin{CD}
\C_0 @>{\id \times s}>> \C_0 \times [1]\\
@V{\gamma}VV @VVV\\
\C_1 @>>> \Cyl(\gamma),
\end{CD}
\end{equation}
where $s = \eps(0),t=\eps(1):\ppt \to [1]$ are the embedding onto
$0,1 \in [1]$. Explicitly, the collection of objects in
$\Cyl(\gamma)$ is the disjoint union of objects in $\C_0$ and
$\C_1$, and similarly for $\Cyl^o(\gamma)$, while morphisms are
given by
\begin{equation}\label{cyl.exp}
\begin{aligned}
\Cyl(\gamma)(c,c') &= \begin{cases} \C_l(c,c'), &\quad c,c' \in
  \C_l,l=0,1\\ \C_1(\gamma(c),c'), &\quad c \in \C_0,c' \in \C_1,\\
  \emptyset, &\quad c \in \C_1,c' \in \C_0,
\end{cases}\\
\Cyl^o(\gamma)(c,c') &= \begin{cases} \C_l(c,c'), &\quad c,c' \in
  \C_l,l=0,1\\ \C_1(c,\gamma(c')), &\quad c \in \C_1,c' \in \C_0,\\
  \emptyset, &\quad c \in \C_0,c' \in \C_1.
\end{cases}
\end{aligned}
\end{equation}  
The embeddings $s$, $t$ induce full embeddings $s:\C_0 \to
\Cyl(\gamma)$, $t:\C_1 \to \Cyl(\gamma)$, the embedding $t$ is
left-reflexive, with the adjoint functor $t_\dg$, and $\gamma \cong
t_\dg \circ s$. Dually, we also have full embeddings $s:\C_1 \to
\Cyl^o(\gamma)$, $t;\C_0 \to \Cyl^o(\gamma)$, $s$ is right-reflexive
with right-adjoint $s^\dg$, and $\gamma \cong s^\dg \circ
t$. Altogether, $\gamma$ has two canonical factorizations
\begin{equation}\label{cyl.facto}
\begin{CD}
\C_0 @>{s}>> \Cyl(\gamma) @>{t_\dg}>> \C_1 \qquad
\C_0 @>{t}>> \Cyl^o(\gamma) @>{s^\dg}>> \C_1.
\end{CD}
\end{equation}
For any functors $\rho:\C \to \C'$, $\lambda:\C' \to \C$, an
adjunction between $\rho$ and $\lambda$ defines an equivalence
$\alpha:\Cyl(\lambda) \cong \Cyl^o(\rho)$ equipped with isomorphisms
$a_s:\alpha \circ s \cong s$, $a_t:\alpha \circ t \cong t$, and
isomorphism classes of triples $\langle \alpha,a_s,a_t \rangle$
correspond bijectively to pairs of adjunction maps $a$, $a^\dg$ such
that \eqref{adj.eq} are the identity maps. If $\tau:\C \to \ppt$ is
the tautological projection, we have $\C^> \cong \Cyl(\tau)$, $\C^<
\cong \Cyl^o(\tau)$. We also have $s \cong \eps:\C \to \C^>$, $t
\cong \eps(o):\ppt \to \C^>$, and dually for $\C^<$.

For any functor $\pi:\C \to I$, the {\em left} resp.\ {\em right
  comma-category} $\C /_\pi I$ resp.\ $I \setminus_\pi \C$ are
defined by cartesian squares
\begin{equation}\label{comma.sq}
\begin{CD}
\C /_\pi I @>>> \Ar(I)\\
@V{\sigma}VV @VV{\sigma}V\\
\C @>{\pi}>> I,
\end{CD}
\qquad\qquad
\begin{CD}
I \setminus_\pi \C @>>> \Ar(I)\\
@V{\tau}VV @VV{\tau}V\\
\C @>{\pi}>> I,
\end{CD}
\end{equation}
and come equipped with functors
\begin{equation}\label{si.ta.co}
  \tau:\C /_\pi I \to I, \qquad \sigma:I \setminus_\pi \C \to I.
\end{equation}
Explicitly, objects in $\C /_\pi I$ are triples $\langle c,i,\alpha
\rangle$, $c \in \C$, $i \in I$, $\alpha:\pi(c) \to i$ a map, and
dually for $I \setminus_\pi \C$. The functor $\tau:I \setminus_\pi
\C$ in \eqref{comma.sq} has a right-adjoint $\eta:\C \to I
\setminus_\pi \C$ induced by $\eta:I \to \Ar(I)$, and dually,
$\sigma$ in \eqref{comma.sq} has a left-adjoint $\eta:\C \to \C
/_\pi I$. The functor $\pi:\C \to I$ then decomposes as
\begin{equation}\label{comma.facto}
\begin{CD}
\C @>{\eta}>> \C /_\pi I @>{\tau}>> I,\qquad
\C @>{\eta}>> I \setminus_\pi \C @>{\sigma}>> I,
\end{CD}
\end{equation}
where $\sigma$ and $\tau$ are the functors \eqref{si.ta.co}. For any
object $i \in I$, the {\em fiber} $\C_i$ of the functor $\pi$ is
given by $\C_i = \eps(i)^*\C$, where $\eps(i):\ppt \to I$ is the
embedding onto $i$, and the {\em left} resp.\ {\em right
  comma-fibers} are the fibers $\C /_\pi i = (\C /_\pi I)_i$, $i
\setminus_\pi \C)_i = (I \setminus_\pi \C)_i$ of the functors
\eqref{si.ta.co}. We will drop $\pi$ from notation when it is clear
from the context. In particular, if $\C = I$ and $\pi=\id$ is the
identity functor, then $I / I \cong I \setminus I \cong \Ar(I)$, and
for any $i \in I$, $I / i$ resp.\ $i \setminus I$ is the category of
objects $i' \in I$ equipped with a morphism $i' \to i$ resp.\ $i \to
i'$.

\subsection{Fibrations and cofibrations.}\label{fib.subs}

A functor $\pi:\C \to I$ is a {\em fibration} if the induced functor
$\pi\setminus\id:\Ar(\C) \cong \C \setminus \C \to I \setminus_\pi
\C$ admits a fully faithful right-adjoint
$(\pi\setminus\id)^\dg$. In this case, a map $f$ in $\C$ is {\em
  cartesian} if the corresponding object in $\Ar(\C)$ lies in the
essnetial image of the functor $(\pi\setminus\id)^\dg$. The
composition of cartesian maps is automatically cartesian, so that
all objects $c \in \C$ and cartesian maps between them form a
subcategory $\C_\flat \subset \C$. For any fibration $\pi:\C \to I$
and functor $\gamma:I' \to I$, the projection $\gamma^*\C \to I'$ is
a fibration, and for any functor $\gamma:\C \to I$, $\sigma$ in
\eqref{si.ta.co} is automatically a fibration. For any functor
$\gamma:\C_0 \to \C_1$, the projection $\Cyl^o(\gamma) \to \C_0
\times [1] \to [1]$ is a fibration, and it is clear from
\eqref{cyl.exp} that all fibrations over $[1]$ arise in this
way. Thus for any fibration $\pi:\C \to I$ and map $f:i \to i'$ in
$I$, with the corresponding functor $\eps(f):[1] \to I$, we have
$\eps(f)^* \cong \Cyl^o(f^*)$ for a unique functor $f^*:\C_{i'} \to
\C_i$ called the {\em transition functor} of the fibration $\C \to
I$. ``Unique'' here means ``unique up to a unique isomorphism'', so
that a composable pairs of maps $f:i \to i'$, $g:i' \to i''$ gives
rise to a canonical isomorphism $(g \circ f)^* \cong f^* \circ g^*$
satisfying a compatibility condition for composable triples; the
whole gadget is called a ``pseudofunctor''. Both the notion of a
fibration and a description of fibrations in terms of pseudofunctors
--- and vice versa --- were introduced by Grothendieck in
\cite{SGA}, and these days, this is known as the {\em Grothendieck
  construction}. A fibration $\C \to I$ is {\em constant} along a
map $f:i \to i'$ in $I$ if the transition functor $f^*$ is an
equivalence. A fibration $\C \to [1]^2$ is the same thing as a
commutative square \eqref{cat.sq}, and we say that a fibration $\C
\to I$ is {\em cartesian} resp.\ {\em semicartesian} over a
commutative square $\eps:[1]^2 \to I$ if so is the square
\eqref{cat.sq} corresponding to $\eps^*\C \to [1]^2$. We say that a
functor $\pi:\C \to I$ is a {\em family of groupoids} if
$\pi\setminus\id:\Ar(\C) \to I \setminus_\pi \C$ is an equivalence,
or equivalently, if $\pi$ is a fibration whose fibers $\C_i$, $i \in
I$ are groupoids. A functor $X:I^o \to \Sets \subset \Cat$ defines a
small fibration $IX \to X$ with discrete fibers, it is a family of
groupoids, and any small fibration with discrete fibers arises in
this way. One calls $IX$ the {\em category of elements} of the
functor $X$.

Dually, a functor $\pi:\C \to I$ is a {\em cofibration} if
$\pi^o:\C^o \to I^o$ is a fibration, or equivalently, if
$\id/\pi:\Ar(\C) \to \C/_\pi I$ has a fully faithful left-adjoint;
then $\tau$ in \eqref{si.ta.co} is a cofibration, cofibrations are
preserved by pullback, cofibrations over $[1]$ are given by
cylinders $\Cyl(\gamma)$, and for any cofibration $\C \to I$ and map
$f:i \to i'$ in $I$, we have the canonical transition functor
$f_!:\C_i \to \C_{i'}$. A functor $\C \to I$ is a {\em bifibration}
if it is both a fibration and a cofibration; in this case, $f_!$ is
left-adjoint to $f^*$ for any map $f$ in $I$, and converserly, a
fibration is a bifibration if all its transition functors $f^*$ are
left-reflexive.

For any fibration $\pi:\C \to I$, with fibers $\C_i$, $i \in I$, and
transition functors $f^*:\C_{i'} \to \C_i$ for morphisms $f:i \to
i'$ in $I$, one can construct the {\em transpose cofibration}
$\pi_\perp:\C_\perp \to I^o$ with the same fibers $(\C_\perp)_i =
\C_i$, $i \in I$, and transition functors $f^o_! = f^*$; the
isomorphisms $(f^o \circ g^o)_! \cong f^o_! \circ g^o_!$ are inverse
to the canonical isomorphisms $f^* \circ g^* \cong (g \circ
f)^*$. More invariantly, $\C_\perp$ is the category with the same
objects as $\C$, and morphisms from $c$ to $c'$ in $\C_\perp$ are
given by isomorphism classes of span diagrams
\begin{equation}\label{trans.dia}
\begin{CD}
  c @<{v}<< \wt{c} @>{f}>> c'
\end{CD}
\end{equation}
in $\C$ such that $\pi(v)$ is invertible, and $f$ is cartesian over
$I$. One checks that diagrams \eqref{trans.dia} have no non-trivial
automorphisms, and their composition is obtained by taking the
appropriate pullbacks in $\C$ that exist since $\pi:\C \to I$ is a
fibration. The functor $\pi_\perp$ sends $c \in \C_\perp$ to
$\pi(c)$, and a diagram \eqref{trans.dia} goes to $(\pi(v)^{-1}
\circ \pi(f))^o:\pi(c') \to \pi (c)$. The opposite functor
$\pi^o_\perp:\C^o_\perp \to I$ is a fibration with fibers $\C_i^o$
and transition functors $(f^*)^o$. If we have another fibration $\C'
\to I$, then a functor $\gamma:\C \to \C'$ cartesian over $I$
induces a functor $\gamma_\perp:\C_\perp \to \C'_\perp$ cocartesian
over $I^o$.

A fully faithful functor $\pi:\C' \to \C$ is a fibration if and only
if it is {\em left-closed}, in the sense that for any map $c \to c'$
in $\C$ such that $c' \in \C'$, we also have $c \in \C'$. For
example, the embedding $s = \eps(0):\ppt \to [1]$ onto $0 \in [1]$ is
left-closed, and this example is in fact universal: for any
left-closed full embedding $\C' \to \C$, we have $\C' \cong
\chi^*\ppt$ for a unique functor $\chi:\C \to [1]$. This is actually
the Grothendieck construction again: the fibers of a fully faithful
functor are either $\ppt$ or empty, and the full subcategory in
$\Cat$ spanned by $\ppt$ and $\emptyset$ is $[1]$. Dually, a fully
faithful $\pi:\C' \to \C$ is {\em right-closed} if $\pi^o$ is
left-closed, this happens iff $\pi$ is a cofibration, and the
universal right-closed full embedding is the embedding $t =
\eps(1):\ppt \to [1]$ onto $1 \in [1]$.

\begin{remark}
Technically, the definition of a fibration in \cite{SGA} is slightly
different: a functor $\pi:\C \to I$ is a fibration in our sense if
and only if $I \times_I \C \to I$ is a fibration in the original
sense of Grothendieck. Since we have a canonical equivalence $\C
\cong I \times_I \C$, this does not create any problems. We allow
ourselves to modify the definition since we want all categorical
notions to be invariant under equivalences.
\end{remark}

\begin{remark}
These days, what we call a cofibration is sometimes called an
``opfibration'', to avoid a clash of terminology with Quillen's
machinery of model categories. At the time of \cite{SGA}, this
machinery did not exist, so there was no problem. In our approach to
homotopical algebra, model categories are not much used either, so
we allow ourselves to go full circle and return to the original
terminology of \cite{SGA}.
\end{remark}

For any categories $\C_0$, $\C_1$ equipped with functors $\pi_l:\C_l
\to I$, $l=0,1$, a {\em lax functor} from $\C_0$ to $\C_1$ {\em over
  $I$} is a pair $\langle \gamma,\alpha \rangle$ of a functor
$\gamma:\C_0 \to \C_1$ and a morphism $\alpha:\pi_0 \to \pi_1 \circ
\gamma$. A {\em functor over $I$} is a lax functor $\langle
\gamma,\alpha \rangle$ over $I$ with invertible $\alpha$. We note
that there is a more general motion of a lax functor used in the
theory of $2$-categories, see e.g.\ \cite{ka.adj}, but we will not
need it. For any two lax functors $\langle \gamma,\alpha \rangle$,
$\langle \gamma',\alpha' \rangle$, a {\em morphism $\langle
  \gamma,\alpha \rangle \to \langle \gamma',\alpha' \rangle$ over
  $I$} is a morphism $b:\gamma \to \gamma'$ such that $\pi_1(b)
\circ \alpha = \alpha'$. Dually, for any categories $\C_0$, $\C_1$
equipped with functors $\phi_l:I \to \C_l$, a {\em lax functor} from
$\C_0$ to $\C_1$ {\em under $I$} is a pair $\langle \gamma,\alpha
\rangle$ of a functor $\gamma:\C_0 \to \C_1$ and a morphism
$\alpha:\phi_0 \to \phi_1 \circ \gamma$. A {\em functor under $I$}
is a lax functor $\langle \phi,\alpha \rangle$ under $I$ such that
$\alpha$ is an isomorphism. A {\em morphism $\langle \gamma,\alpha
  \rangle \to \langle \gamma',\alpha \rangle$ over $I$} between two
lax functors under $I$ is a morphism $b:\gamma \to \gamma'$ such
that $\alpha' \circ \phi_1(b) = \alpha$. If a functor $\gamma$ over
$I$ is left or right-reflexive, with some adjoint $\gamma^\dg$, then
we say that $\gamma$ is {\em left} resp.\ {\em right-reflexive over
  $I$} if the corresponding base change maps \eqref{bc.eq} are
isomorphisms, so that $\gamma^\dg$ is also a functor over $I$. If
$\pi_0$ and $\pi_1$ are fibrations, then a functor $\gamma:\C_0 \to
\C_1$ over $I$ is {\em cartesian} if it sends cartesian maps to
cartesian maps (or equivalently, the base change map \eqref{bc.eq}
for the adjoint functor $(\id\setminus\pi_l)^\dg$, $l=0,1$ is an
isomorphism). Dually, if $\pi_0$ and $\pi_1$ are cofibrations, then
$\gamma$ over $I$ is {\em cocartesian} if $\gamma^o$ is
cartesian. More generally, if we are given a commutative square
\eqref{cat.sq} such that $\gamma_1$ and $\gamma_{01}^0$ are
fibrations resp.\ cofibrations, then we say that $\gamma_{01}^1$ is
{\em cartesian} resp.\ {\em cocartesian over $\gamma_0$} if it sends
cartesian resp.\ cocartesian maps to cartesian resp.\ cocartesian
maps (or equivalently, if \eqref{cat.sq.eq} is cartesian
resp.\ cocartesian over $\C_0$).

As a slighly non-trivial application of the Grothendieck
construction, one can construct a relative version of the functor
categories $\Fun(-,-)$ of Subsection~\ref{cat.subs}. Assume given a
functor $I' \to I$, and some category $\E$. The {\em relative
  functor category} $\Fun(I'|I,\E)$ is a category over $I$ equipped
with an evaluation functor $I' \times_I \Fun(I'|I,\E) \to \E$ such
that for any category $\C$ over $I$, a functor $\gamma:I' \times_I
\C \to \E$ factors as
\begin{equation}\label{rel.fun.facto}
\begin{CD}
I' \times_I \C @>{\id \times \wt{\gamma}}>> I' \times_I \Fun(I'|I,\E)
@>{\ev}>> \E
\end{CD}
\end{equation}
for a functor $\wt{\gamma}:\C \to \Fun(I'|I,\C)$ over $I$, unique up to
a unique isomorphism. Then while it is not true that relative
functor categories exists for any functor $I' \to I$, it is true
when the functor is a small fibration or cofibration (\cite[Lemma
  1.4.1.2, Remark 1.4.1.5]{big}). In this case, $\Fun(I'|I,\E) \to
I$ is a cofibration resp.\ fibration with fibers $\Fun(I'|I,\E)_i
\cong \Fun(I'_i,\E)$, $i \in I$.

\subsection{Limits and Kan extensions.}\label{lim.subs}

For any categories $I$, $\E$, a {\em cone} of a functor $E:I \to \E$
is a functor $E_>:I^> \to \E$ equipped with an isomorphism
$\eps^*E_> \cong E$, and $E_>(o) \in \E$ is the {\em vertex} of the
cone. All cones for a fixed $E$ form a category, and a cone is {\em
  universal} if it is an initial object in this category. If the
universal cone exists, then $E$ {\em has a colimit}, and the colimit
$\colim_I E$ is the vertex of this universal cone. A functor $F:\E
\to \E'$ {\em preserves the colimit} $\colim_IE$ if the cone $F
\circ E_>$ for $F \circ E$ is universal. For example, we have $[1]^2
\cong \V^>$, where
\begin{equation}\label{V.eq}
  \V = \{0,1\}^<
\end{equation}
is the category with three objects $o,0,1$ and two non-identity maps
$o \to 0,1$; then any commutative square in a category $\E$ is the
cone of a functor $\V \to \E$, and the cone is universal iff the
square is cocartesian. {\em Coproducts} are colimits over discrete
categories.

More generally, for any functor $\gamma:I \to I'$, a {\em relative
  cone} of a functor $E:I \to \E$ is a functor $E_>:\Cyl(\gamma) \to
\E$ equipped with an isomorphism $s^*E_> \cong E$, where $s:I \to
\Cyl(\gamma)$ is the embedding \eqref{cyl.facto}. Again, a relative
cone is {\em universal} if it admits a unique map into any other
relative cone. If a universal relative cone $E_>$ exists, then $E$
{\em has a left Kan extension} with respect to $\gamma$, and this
left Kan extension is $\gamma_!E = t^*E_>$. A functor $F:\E \to \E'$
{\em preserves the Kan extension} $\gamma_!E$ if $F \circ E_>$ is a
universal relative cone. Sometimes all functors admits left Kan
extensions; for example, this happens if $\gamma$ is
right-reflexive, with adjoint $\gamma^\dg$, and then
\begin{equation}\label{kan.adj}
  \gamma_!E \cong \gamma^{\dg*}E,
\end{equation}
functorially with respect to $E$. Alternatively, one can impose
conditions on $E$. For example, a category $E$ is {\em cocomplete}
if any functor $I \to E$ from an essentially small $I$ has a
colimit; in this case, for any functor $\gamma:I \to I'$ between
essentially small categories, $\gamma_!:I' \to \E$ exists for any
$E:I \to \E$. If $I$ and $I'$ are small, and $\gamma_!E$ exists for
any $E:I \to \E$, then it is functorial in $E$ and defines a functor
$\gamma_!:\Fun(I,\E) \to \Fun(I',\E)$ left-adjoint to the pullback
functor $\gamma^*:\Fun(I',\E) \to \Fun(I,\E)$. For any cocomplete
category $\E$ and essentially small $I$, $\Fun(I,\E)$ is cocomplete;
in particular, $\Sets$ is cocomplete, and then so are
$\Fun(I,\Sets)$ and $I^o\Sets$.

In practice, to compute left Kan extensions, one can use the
following base change result. Assume given a diagram
\begin{equation}\label{bc.dia}
\begin{CD}
I'_0 @>{\gamma'}>> I'_1 @>{\pi'}>> I'\\
@V{\phi_0}VV @VV{\phi_1}V @VV{\phi}V\\
I_0 @>{\gamma}>> I_1 @>{\pi}>> I
\end{CD}
\end{equation}
of essentially small categories and functors such that the squares
are commutative and cartesian, $\pi$ and $\pi \circ \gamma$ are
cofibrations, and $\gamma$ is cocartesian over $I$. Then for any
target category $\E$ such that $\gamma_!:\Fun(I_0,\E) \to
\Fun(I_1,\E)$ and $\gamma'_!:\Fun(I'_0,\E) \to \Fun(I_1,\E)$ exist,
the base change map
\begin{equation}\label{kan.bc.eq}
  \gamma'_! \circ \phi_0^* \to \phi_1^* \circ \gamma_!
\end{equation}
of functors $\Fun(I_0,\E) \to \Fun(I',\E)$ is an isomorphism. In
particular, if one takes $I'=\ppt$ and $\pi=\id$, then
\eqref{kan.bc.eq} allows one to compute left Kan extensions
$\gamma_!E$ along a cofibration $\gamma$ pointwise; combining this
with \eqref{kan.adj} and using the decomposition
\eqref{comma.facto}, one obtains a functorial isomorphism
\begin{equation}\label{cokan.eq}
\gamma_!E(i) \cong \colim_{I' / i}E, \qquad i \in I,\quad E \in
\Fun(I',\E)
\end{equation}
for any functor $\gamma:I' \to I$ between essentially small
categories, and any target category $\E$ that admits all the
colimits in the right-hand side.

Dually, the {\em limit} $\lim E$ of a functor $E:I \to \E$ is
$\lim_I E = (\colim_{I^o}E^o)^o$, the {\em right Kan extension}
$\gamma_*E$ with respect to a functor $\gamma:I \to I'$ is given by
$\gamma_*E = (\gamma^o_!E^o)^o$, and both have the the dual versions
of all the properties of colimits and left Kan extensions. In
particular, we have base change isomorphisms for diagrams
\eqref{bc.dia} of fibrations and cartesian functors, and a dual
version of the functorial isomorphisms \eqref{kan.bc.eq}. If we are
given a functor $\gamma:I' \to I$ between essentially small
categories, and any functor $E:I' \to \E$, we then have isomorphisms
\begin{equation}\label{kan.eq}
\gamma_*E(i) \cong \lim_{i \setminus I'}E, \qquad i \in
I,\quad E \in \Fun(I',\E),
\end{equation}
a dual version of \eqref{cokan.eq}; if all the limits in the
right-hand side exists, then $f_*E$ also exists, and is given by
\eqref{kan.eq}. Limits over discrete categories are products, and
limits over $\V^o = \{0,1\}^<$ are cartesian squares. A category
$\E$ is {\em complete} if it has limits of all functors $I \to \E$
from an essentially small $I$, and in this case, $\Fun(I,\E)$ is
also complete, and all the right Kan extensions with target $\E$
exist. In particular, $\Sets$ is complete, and so is $I^o\Sets$ for
any essentially small $I$. The category $\Cat$ is also complete (in
fact, it is also cocomplete, but colimits in $\Cat$ are pathological
-- this is where the whole need for enhancements comes from).

\begin{remark}
It might happen that the right Kan extension $\gamma_*E$ of a
functor $E:I' \to \E$ along a functor $\gamma:I' \to I$ between
essentially small categories exists even while some of the limits in
\eqref{kan.eq} do not (see e.g.\ \cite[Remark 7.5.3.10]{big}). If
the target category $\E$ is essentially small, one can say that a
right Kan extension is {\em universal} if it is preserved by the
Yoneda embedding $\Y:\E \to \E^o\Sets$. Even if $\E$ is not
essentially small, $E$ together with $\gamma_*E$ factor through an
essentially small full subcategory $\E' \subset \E$; one says that
$\gamma_*E$ is {\em universal} if it is universal as Kan extension
of functors $I' \to \E'$, for any such $\E'$. Then all limits are
automatically universal, and a universal Kan extension $\gamma_*E$
exists if and only if so do all the limits in \eqref{kan.eq}. The
situation for colimits and left Kan extensions is dual.
\end{remark}

A category $\C$ with finite products is {\em cartesian-closed} if
for any $c \in \C$, the functor $c \times -:\C \to \C$ has a
right-adjoint $\Hhom(c,-)$. The categories $\Sets$, $\Cat$ are
cartesian-closed, and so are the categories $I^o\Sets$, $I^o\Cat$
for any essentially small $I$.

\subsection{Localization and the Yoneda embeddings.}\label{yo.subs}

A {\em localization} $h^W(\C)$ of a category $\C$ with respect to a
class of morphisms $W$ is defined by a cocartesian square
\begin{equation}\label{loc.sq}
\begin{CD}
W \times [1] @>>> \C\\
@VVV @VV{h}V\\
W @>>> h^W(\C),
\end{CD}
\end{equation}
if it exists (equivalently, the localization functor $h:\C \to
h^W(\C)$ inverts all maps $w \in W$ and is universal with this
property). By universality, the localization is unique up to a
unique isomorphism. Existence is usually non-trivial; if it holds,
one says that $\C$ is {\em localizable} with respect to the class
$W$. It always holds when $\C$ is small, but even then, it is
usually quite hard to describe the resulting category $h^W(\C)$. One
situation when it is easy is when $\C = \Cat$, and $W$ is the class
of all equivalences. In this case, $h^W(\Cat)$ is the category
$\Cat^0$ whose objects are small categories, and whose morphisms are
isomorphisms classes of functors. More generally, for any small
category $I$, we can localize $I^o\Cat$ with respect to pointwise
equivalences; the result is the category $(I^o\Cat)^0$ of functors
$I^o \to \Cat$ and pointwise functors between them considered modulo
pointwise isomorphisms. Note that this is {\em not} the same as
considering functors to $\Cat^0$ -- we always have a tautological
functor
\begin{equation}\label{loc.cat.eq}
(I^o\Cat)^0 \to I^o\Cat^0,
\end{equation}
but it is only an equivalence if $I$ is discrete. Sometimes it is an
epivalence (e.g. when $I=[1]$ or $I=\V$, see \cite[Lemma
  4.1.1.9]{big}).

For any category $I$, we denote by $\Cat \bb I$ the category of
small categories $\C$ equipped with a functor $\pi:\C \to I$, with
morphisms given by lax functors over $I$, and we let $\Cat \bbi I
\subset \Cat \bb I$ be the subcategory with the same objects and
morphisms given by functors over $I$. The forgetful functor
\begin{equation}\label{cat.bb}
  \Cat \bb I \to \Cat, \qquad \langle \C,\pi \rangle \mapsto \C
\end{equation}
is a fibration with fibers $(\Cat \bb I)_\C \cong \Fun(\C,I)$, and
$\Cat \bbi I \subset \Cat \bb I$ is the subcategory of maps
cartesian over $\Cat$. Dually, if $I$ is essentially small, we
denote by $I \bbd \Cat$ the category of small categories $\C$
equipped with a functor $\phi:I \to \C$, with morphisms given by lax
functors under $I$, and we let $I \bbdi \Cat \subset I \bbd \Cat$ be
the subcategory with the same objects and morphisms given by
functors under $I$. The forgetful functor $I \bbd \Cat \to \Cat$ is
a cofibration with fibers $(I \bbd \Cat)_\C \cong \Fun(I,\C)$, and
$I \bbdi \Cat \subset I \bbd \Cat$ is the subcategory of maps
cocartesian over $\Cat$. We let $\Cat_\idot = \ppt \bbd \Cat$;
explicitly, this is the category of pairs $\langle \C,c \rangle$, $c
\in \C \in \Cat$, and the cofibration $\Cat_\idot \to \Cat$ has
fibers $(\Cat_\idot)_\C \cong \C$. For any category $I$, we have
\begin{equation}\label{cat.bb.fib}
\Cat \bb I \cong \Fun(\Cat_\idot|\Cat,I),
\end{equation}
where $\Fun(\Cat_\idot|\Cat,-)$ is the relative functor category of
Subsection~\ref{fib.subs} (indeed, both categories are fibered over
$\Cat$, with fiber $\Fun(\C,I)$ over $\C \in \Cat$). If $I$ is
essentially small, we have a full embedding
\begin{equation}\label{yo.cat.eq}
\Y:\Cat \bb I \to I^o\Cat
\end{equation}
sending $\langle \C,\pi \rangle$ to the functor $i \mapsto i
\setminus_\pi \C$; this is an extended version of the Yoneda
embedding (and restricts to the usual Yoneda embedding \eqref{yo.eq}
on the fiber $(\Cat \bb I)_{\ppt} \cong \Fun(\ppt,I) \cong I$ of the
forgetful fibration \eqref{cat.bb.fib}).

To understand the Yoneda Lemma in this language, it is useful to
distinguish an even smaller subcategory $\Cat \mm I \subset \Cat
\bbi I$ of fibrations and cartesian functors over $I$ (and dually,
the subcategory $\Cat \mc I \subset \Cat \bbi I$ of cofibrations and
cocartesian functors). Altogether, we have a circular diagram
\begin{equation}\label{3.dia}
\begin{CD}
\Cat \mm I @>{a}>> \Cat \bbi I @>{b}>> \Cat \bb I @>{y}>> \Cat
\mm I,
\end{CD}
\end{equation}
where $a$ and $b$ are the embeddings, and $y$ sends $\pi:\C \to I$
to the fibration $\sigma:I \setminus_\pi \C \to I$ of
\eqref{si.ta.co}. For any $\pi:\C \to I$, the decomposition
\eqref{comma.facto} provides a functor $\eta:\C \to I \setminus_\pi
\C$ over $I$, its left-adjoint $\tau:I \setminus_\pi \C \to \C$ is
canonically a lax functor over $I$, and if $\pi$ is a fibration,
then $\eta$ also has a right-adjoint $\eta^\dg = \sigma \circ
(\pi\setminus\id)^\dg:I \setminus_\pi \C \to \C$ cartesian over
$I$. In terms of \eqref{3.dia}, this defines maps $\id \to a \circ y
\circ b$, $\id \to b \circ a \circ y$, $y \circ b \circ a \to \id$
that provide a sort of a $2$-categorical adjunction between $a$ and
$y \circ b$ and between $a \circ y$ and $b$.

To make it into a genuine adjunction, one can localize with respect
to equivalences. Just as in the absolute case $I=\ppt$, all the
categories in \eqref{3.dia} admits localizations $(\Cat \mm I)^0$,
$(\Cat \bbi I)^0$, $(\Cat \bb I)^0$; these have the same objects as
the categories we localize, and morphisms are cartesian functors
resp.\ functors resp.\ lax functors considered up to an isomorphism
over $I$. Then \eqref{3.dia} induces a diagram of localized
categories where $a$ is genuinely left-adjoint to $y \circ b$, and
$a \circ y$ is genuinely right-adjoint to $b$. It is a purely formal
consequences of this fact (see \cite[Lemma 1.2.7]{big}) that
\begin{equation}\label{yo.loc}
y:(\Cat \bb I)^0 \to (\Cat \mm I)^0
\end{equation}
is fully faithful. The Grothendieck construction provides an
equivalence $(I^o\Cat)^0 \cong (\Cat \mm I)^0$, and under this
identification, \eqref{yo.loc} is induced by the extended Yoneda
embedding \eqref{yo.cat.eq}.

The localized categories $\Cat^0$, $(I^o\Cat)^0 \cong (\Cat \mm
I)^0$ are still cartesian-closed, but they are no longer complete,
nor cocomplete, and neither are the categories $(\Cat \bb I)^0$,
$(\Cat / I)^0$ (to observe this, note that already cartesian squares
\eqref{cat.sq} of small categories are {\em not} cartesian
squares in $\Cat^0$ -- the universal property needs the actual
isomorphism $\alpha:\gamma_0 \circ \gamma^0_{01} \to \gamma_1 \circ
\gamma^1_{01}$, not just the fact that some isomorphism exists, and
passing to $\Cat^0$ forgets $\alpha$). Analogously, localizing the
category $I \bbd \Cat$ with respect to equivalences gives the
category $(I \bbd \Cat)^0$ of functors $\phi:I \to \C$, with
morphisms given by lax funders under $I$ considered up to an
isomorphism under $I$, and this category is neither complete nor
cocomplete either. Moreover, the projection $(\Cat \bb I)^0 \to
\Cat^0$ is not a fibration, and the projection $(I \bbd \Cat)^0 \to
\Cat^0$ is no longer a cofibration. If we take $I = \ppt$,
$\Cat_\idot^0 = (\ppt \bbd \Cat)^0$, then we have a commutative
square
\begin{equation}\label{cat.loc}
\begin{CD}
\Cat_\idot @>>> \Cat_\idot^0\\
@VVV @VVV\\
\Cat @>>> \Cat^0,
\end{CD}
\end{equation}
and this square is semicartesian, but not cartesian. To see this
effect explicitly, consider a group $G$, with the corresponding
groupoid $\ppt_G \in \Cat$. Then we have an exact
sequence
\begin{equation}\label{G.Z.dia}
\begin{CD}
1 @>>> Z @>>> G @>{a'}>> \Aut(G) @>>> G' @>>> 1,
\end{CD}
\end{equation}
where $Z \subset G$ is the center of the group $G$, $\Aut(G)$ is the
automorphism group of $G$, and $a$ is the adjoint action map. In
terms of \eqref{G.Z.dia}, we have $\Aut_{\Cat}(\ppt_G) = \Aut(G)$
and $\Aut_{\Cat^0}(\ppt_G) = G'$.  Moreover, since $\ppt_G$ has a
unique object, it canonically defines an object in $\Cat_\idot$, and
we have $\Aut_{\Cat_\idot}(\ppt_G) \cong G \rtimes \Aut(G)$,
$\Aut_{\Cat_\idot^0}(\ppt_G) \cong (G \rtimes \Aut(G))/G \cong
\Aut(G)$. Then \eqref{cat.loc} induces a semicartesian square
$$
\begin{CD}
\ppt_{G \rtimes \Aut(G)}@>>> \ppt_{\Aut(G)}\\
@VVV @VVV\\
\ppt_{\Aut(G)} @>>> \ppt_{G'},
\end{CD}
$$
and this square is not cartesian unless $Z$ is trivial.

\subsection{Accessible categories.}\label{acc.subs}

As we have mentioned in Subsection~\ref{cat.subs}, even in a large
category, $\Hom$-sets are assumed to be small, and this has the
well-known unpleasant side effect: functors between two large
categories may not form a category (there may be too many morphisms
between them). To circumvent this difficulty, one can only consider
functors that preserve some sort of colimits, and the standard way
to do it is to use {\em filtered colimits} and {\em accessible
  categories}. This theory was originated by Grothendieck in
\cite{SGA4} and developed by Gabriel and Ulmer \cite{GU}; in its
present from, the theory appears in \cite{MP}, and then in
\cite{AR}. The latter is a wonderful book that is the standard
reference on the subject, so we only recall the bare details (an
overview that uses the same notation and technology as in this paper
is available in \cite{anan}).

As usual, a {\em cardinal} is an isomorphism class of sets. For any
set $S$, $|S|$ is its cardinality, and we write $|S| \leq |S'|$ if
there is an injective map $S \to S'$. A {\em chain} in a small
category $I$ is a diagram
\begin{equation}\label{i.dot.dia}
\begin{CD}
i_0 @>>> \dots @>>> i_n
\end{CD}
\end{equation}
of some length $n \geq 0$, and a chain is {\em non-degenerate} if
none of the maps $i_l \to i_{l+1}$ is an identity map. We write
$\Ch(I)$ for the set of all non-degenerate chains in a small
category $I$, and we denote $|I|=|\Ch(I)|$. For any cardinal
$\kappa$, we let $\Sets_\kappa \subset \Sets$ resp.\ $\Cat_\kappa
\subset \Cat$ be the full subcategory spanned by sets $S$
resp.\ small categories $I$ such that $|S| < \kappa$ resp.\ $|I| <
\kappa$. A category $\E$ is {\em $\kappa$-complete} resp.\ {\em
  $\kappa$-cocomplete} if $\lim E$ resp.\ $\colim E$ exists for any
functor $E:I \to E$ from a small category $I$ such that $|I| <
\kappa$. A cardinal $\kappa$ is {\em regular} if $\Sets_\kappa$ is
$\kappa$-cocomplete. A category $\E$ is {\em $\kappa$-filtered} if
any functor $E:I \to \E$ from a small category $I$ with $|I| <
\kappa$ admits a cone $E_>$. A category $\E$ is {\em
  $\kappa$-filtered-cocomplete} if $\colim E$ exists for any functor
$E:I \to \E$ from a $\kappa$-filtered $I$ (that is, $\E$ has all
$\kappa$-filtered colimits). An object $e \in \E$ in a
$\kappa$-filtered-cocomplete category $\E$ is {\em $\kappa$-compact}
if $\Hom(e,-):\E \to \Sets$ preserves $\kappa$-filtered colimits,
and we denote by $\Comp_\kappa(\E) \subset \E$ the full subcategory
of $\kappa$-compact objects. For any category $\C$, the {\em
  $\kappa$-inductive completion} $\Ind_\kappa(\C)$ is the category
of pairs $\langle I,c \rangle$ of a $\kappa$-filtered small category
$\I$ and a functor $c:I \to \C$, with morphisms from $\langle I,c
\rangle$ to $\langle I',c' \rangle$ given by
\begin{equation}\label{ind.eq}
\Hom(\langle I,c\rangle,\langle I',c' \rangle) = \lim_{i \in
  I^o}\colim_{i' \in I'}\Hom(c(i),c'(i')).
\end{equation}
We have a fully faithful embedding $\C \to \Ind_\kappa(\C)$, the
category $\Ind_\kappa(\C)$ is $\kappa$-filtered-cocomplete, and it
is universal with this property: any functor $\C \to \E$ to a
$\kappa$-filtered-cocomplete category $\E$ extends to a functor
$\Ind_\kappa(\C) \to \E$ that preserves $\kappa$-filtered colimits,
and such an extension is unique up to a unique isomorphism. In
particular, if a category $\C$ is already $\kappa$-filtered
cocompete, we have a unique comparison functor
\begin{equation}\label{ind.comp.eq}
\Ind_\kappa(\Comp_\kappa(\C)) \to \C
\end{equation}
that preserves $\kappa$-filtered colimits. It is easy to see from
\eqref{ind.eq} that this functor is automatically fully faithful.

\begin{defn}
For any regular cardinal $\kappa$, a category $\C$ is {\em
  $\kappa$-accessible} if it is $\kappa$-filtered-cocomplete,
$\Comp_\kappa(\C)$ is essnetially small, and the fully faithful
embedding \eqref{ind.comp.eq} is an equivalence. A
$\kappa$-accessible category is {\em $\kappa$-presentable} if it is
cocomplete. A functor $\C \to \C'$ between $\kappa$-accessible
categories is {\em $\kappa$-accessible} if it preserves
$\kappa$-filtered colimits. A category $\C$ is {\em accessible}
resp.\ {\em presentable} if it is $\kappa$-accessible
resp.\ $\kappa$-presentable for some regular cardinal $\kappa$, and
a functor $\C \to \C'$ between accessible categories is {\em
  accessible} if it is $\kappa$-accessible for some $\kappa$ such
that both $\C$ and $\C'$ are $\kappa$-accessible.
\end{defn}

If a category $\C$ is $\kappa$-presentable, then it is also
$\mu$-presentable for any regular cardinal $\mu > \kappa$. For
accessible categories, this is not necessarily true, but for any
regular $\kappa$, there exists a cofinal collection of cardinals
$\mu$ such that the same conclusion holds (for the precise condition
on $\mu$, see e.g.\ \cite[Definition 2.1.6.4]{big}). Thus one can
``raise the accessibility index'' indefinitely, and in any small
collection of accessible categories and functors between them, one
can choose a single $\kappa$ that works for all of them. Moreover,
we have the following two fundamental facts.
\begin{enumerate}
\item For any cartesian square \eqref{cat.sq} such that $\C$,
  $\C_0$, $\C_1$, $\gamma_0$, $\gamma_1$ are accessible, so are
  $\C_{01}$ and $\gamma^l_{01}$, $l=0,1$.
\item Accessible functors $\C \to \E$ between two accessible
  categories $\C$, $\E$ form a well-defined category
  $\Fun(\C,\E)$. We have
\begin{equation}\label{fun.acc}
\Fun(\C,\E) = \bigcup_\kappa\Fun_\kappa(\C,\E),
\end{equation}
where the union is over all regular cardinals $\kappa$ such that
$\C$ and $\E$ are $\kappa$-accessible, and $\Fun_\kappa(\C,\E)
\subset \Fun(\C,\E)$ is the full subcategory spanned by
$\kappa$-accessible functors. All the categories
$\Fun_\kappa(\C,\E)$ are also accessible.
\end{enumerate}
We note that both in \thetag{i} and \thetag{ii} it is essential that
we are allowed to raise the accessibility index (in \thetag{i}, it
suffices to take the successor cardinal $\kappa^+$ to the given
cardinal $\kappa$).

Every small category becomes accessible after Karoubi
completetion. Most of the large categories one works with -- such as
$\Sets$, $\Cat$, the categories of schemes, smooth manifolds,
algebras, modules, and so on and so forth -- are accessible. One is
tempted to conclude that all categories ``in nature'' are
accessible, at least up to the Karoubi completion, and simply
working in the accessible world solves all the problems. However,
there is a caveat: the category $\Sets^o$ opposite to the category
of $\Sets$ is {\em not} accessible, and in general, only rarely
passing to the opposite category preserves accessibility. Thus in
the accessible world, simply taking the opposite category is not
allowed. Whether it is a bug or a feature depends on the point of
view.

\section{Orders and simplices.}\label{pos.sec}

\subsection{Partially ordered sets.}

We denote the category of partially ordered sets by $\Pos$. Every
partially ordered set is a small category in the standard way, so we
have a full embedding $\Pos \to \Cat$. The subcategory $\Pos \subset
\Cat$ contains $\Sets \subset \Cat$, and is closed under limits and
coproducts (but not colimits). All partially ordered sets are rigid
as categories, so that a commutative square \eqref{cat.sq} of
partially ordered sets is simply a commutative square in $\Pos$, and
it is cartesian iff it is cartesian in $\Pos$. For any $J \in \Pos$,
the opposite category $J^o$ is also a partially ordered set, and we
denote by $\iota:\Pos \to \Pos$ the involuton $J \mapsto J^o$. The
categories $J^>$, $J^<$ are also partially ordered sets, and for any
small category $\C \in \Cat$, $\Fun(\C,J)$ is a partially ordered
set. For any map $f:J_0 \to J_1$ in $\Pos$, the comma-categories
$J_0 /_f J_1$, $J_1 \setminus_f J_0$ and the cylinder and the dual
cylinder $\Cyl(f)$, $\Cyl^o(f)$ are also partially ordered sets. We
say that a map $f:J' \to J$ is a {\em full embedding}, or {\em
  left-closed}, or {\em right-closed} when this holds in $\Cat$.

Since partially ordered sets are rigid, the Grothendieck
construction identifies cofibrations $J' \to J$ in $\Pos$ with
honest functors $J \to \Pos$. Thus we have an equivalence
\begin{equation}\label{arc.eq}
\Pos \bb \Pos \cong \Arc(\Pos),
\end{equation}
where $\Arc(\Pos) \subset \Ar(\Pos)$ stands for the subcategory
spanned by cofibrations $f:J' \to J$, with maps from $f_0:J'_0 \to
J_0$ to $f_1:J'_1 \to J_1$ given by commutative squares
\begin{equation}\label{arc.sq}
\begin{CD}
J_0' @>{g'}>> J_1'\\
@V{f_0}VV @VV{f_1}V\\
J_0 @>{g}>> J_1
\end{CD}
\end{equation}
such that $g'$ is cocartesian over $g$. The fibration $\Pos \bb \Pos
\to \Pos$ then corresponds to $\tau:\Arc(\Pos) \to \Pos$ that
becomes a fibration after restriction to $\Arc(\Pos) \subset
\Ar(\Pos)$.

We denote by $\N$ the partially ordered set of integers $n \geq 0$,
with the usual order, and for any $n \in \N$, we denote by $[n] =
\N/n$ the partially ordered set of integers $l$, $0 \leq l \leq
n$. We have $[0] = \ppt$, and for any $n \geq 0$, we let $e:[n] \to
[0]$ be the tautological map. For any $n \geq m \geq 0$, we denote
by $s:[m] \to [n]$ resp.\ $t:[m] \to [n]$ the unique left-closed
resp.\ right-closed full embedding (explicitly, $s(l) = l$ and $t(l)
= l + n - m$ for any $l \in [m] = \{0,\dots,m\}$). Note that $[1]$
is the single arrow category, with embeddings $s,t:[0] \to [1]$ onto
$0$ resp.\ $1$, so the notation is consistent. For any $n > l > 0$,
we have a commutative square
\begin{equation}\label{seg.sq}
\begin{CD}
[0] @>{s}>> [l]\\
@V{t}VV @VV{t}V\\
[n-l] @>{s}>> [n].
\end{CD}
\end{equation}
The product $[1]^2$ is also a partially ordered set, and so is $\V =
\{0,1\}^<$.

The squares \eqref{seg.sq} are both cartesian and cocartesian both
in $\Pos$ and in $\Cat$, and stay cartesian and cocartesian after
applying the forgetful functor $\Pos \to \Sets$. In general, $\Pos$
has all colimits, but they often behave badly (in particular, they
need no be preserved by the embedding $\Pos \to \Cat$ or the
forgetful functor $\Pos \to \Sets$). Another example of a
well-behaved colimit is the following: if we have left-closed full
embeddings $J_{01} \to J_0$, $J_{01} \to J_1$, then there exists a
cocartesian square
\begin{equation}\label{st.sq}
\begin{CD}
J_{01} @>>> J_0\\
@VVV @VVV\\
J_1 @>>> J
\end{CD}
\end{equation}
in $\Pos$. The square is also cocartesian in $\Cat$ and cartesian
--- in effect, we have two left-closed full subsets $J_0,J_1 \subset
J$, and $J_{01}=J_0 \cap J_1$. We call such squares {\em standard
  pushout squares}. If $J_{01}=[0]$ and $J_0=J_1=[1]$, then $J=[1]
\copr_{[0]} [1] \cong \V$, and this example is universal: for every
standard pushout square \eqref{st.sq}, there exists a unique
characteristic map $\chi:J \to \V$ such that $J_l = [1] \times_\V J
= J/l$, $l=0,1$, and $J_{01} = [0] \times_\V J = J_o$.

For some of our purposes, the category $\Pos$ is too big, and we
will need to restrict our attention to smaller full subcategories
that are still large enough. The precise meaning of ``large enough''
is as follows.

\begin{defn}\label{amp.def}
A full subcategory $\I \subset \Pos$ is {\em ample} if it is closed
under finite limits and standard pushouts \eqref{st.sq},
contains any subset $J' \subset J$ of some set $J \in \I$, and
contains all finite partially ordered sets.
\end{defn}

In particular, Definition~\ref{amp.def} implies that for any
cofibration $J' \to J$ in an ample $\I \subset \Pos$, all the fibers
$J'_j \subset J'$, $j \in J$ are also in $\I$, so that
\eqref{arc.eq} induces a full embedding
\begin{equation}\label{arc.I.eq}
\Arc(\I) \to \I \bb \I.
\end{equation}
Ample categories that we will need are described by using the
following notion of ``dimension''.

\begin{defn}\label{dim.def}
A partially ordered set $J$ has {\em chain dimension $\leq n$} for
some integer $n \geq 0$ if for any injective map $f:[m] \to I$, $m
\geq 0$, we have $m \leq n$. A partially ordered set $J$ has {\em
  finite chain dimension} if it has chain dimension $\leq n$ for
some $n \geq 0$. In this case, $\dim J$ is the smallest such integer
$n$.
\end{defn}

The full subcategory spanned by $J \in \Pos$ of finite chain
dimension is ample, and we denote it by $\Posf \subset \Pos$. We say
that $J \in \Pos$ is {\em left-finite} resp.\ {\em left-bounded} if
for any $j \in J$, the comma-set $J/j = \{j' \in J| j' \leq j \}$ is
finite resp.\ has finite chain dimension. We denote by
$\Posf^\pm,\Posf^+ \subset \Pos$ the full subcategories spanned by
left-finite resp.\ left-bounded sets; both are ample. The involution
$\iota:\Pos \to \Pos$, $J \mapsto J^o$ sends $\Posf \subset \Pos$
into itself --- in fact, we have $\dim J^o = \dim J$ --- but does
not preserve either $\Posf^\pm$ or $\Posf^+$. The simplest example
of a left-finite set that is not of finite chain dimension is $\N$.

\subsection{Standard squares.}

Our enhanced categories are defined as fibrations over $\Pos$ or its
ample subcategories that are cartesian or semicartesian over certain
commutative squares including but not limited to \eqref{st.sq}; let
us describe these additional squares. Firstly, the {\em barycentric
  subdivision} $B(J)$ of a partially ordered set $J$ is the set of
all finite non-empty totally ordered subsets $S \subset J$, ordered
by inclusion, and we define $\oB(J)$ by the cartesian square
\begin{equation}\label{BJ.sq}
\begin{CD}
\oB(J) @>{\xi'}>> \overline{J}\\
@V{i'}VV @VV{i}V\\
B(J) @>{\xi}>> J,
\end{CD}
\end{equation}
where $\overline{J}$ is $J$ with the discrete order, $i:\overline{J}
\to J$ is the tautological embedding, and $\xi:BJ \to J$ sends $S
\subset J$ to its maximal element. We note that $B(J) \cong B(J^o)$
canonically, $B(J)$ is always left-finite, and if $J$ has finite
chain dimension, then so does $B(J)$, and $\dim B(J) = \dim J$.

Next, take some set $S \in \Sets \subset \Pos$, and let $S^\hush =
B(S^>)^o$. Explicitly, $S^\hush$ fits into a standard pushout square
\eqref{st.sq} of the form
\begin{equation}\label{BS.eq}
\begin{CD}
  S @>>> S^>\\
  @V{\id \times s}VV @VVV\\
  S \times [1] @>>> S^\hush,
\end{CD}
\end{equation}
where the top arrow is the canonical embedding onto $S \cong S^>
\ssetminus \{o\}$ --- in other words, $S^\hush$ is obtained by
taking the product $S \times [1]$, and adding a single new element
$o$ with the order relations $s \times 0 \leq o$, $s \in S$.  The map
$\xi$ of \eqref{BJ.sq} then fits into a cartesian cocartesian square
\begin{equation}\label{S.sq}
\begin{CD}
S^> @>>> \ppt\\
@VVV @VVV\\
S^\hush @>{\xi^o}>> S^<,
\end{CD}
\end{equation}
where the arrow on the right is the embedding onto $o \in S^<$ (and
the arrow on the left then appears in the standard pushout square
\eqref{BS.eq}).

Now consider $\N$ with the standard order, and look at the embedding
$\N \to B(\N)$ sending an even integer $2n \in \N$ to the
single-element subset $\{n\} \subset\N)$, and an odd integer $2n+1$
to the two-element subset $\{n,n+1\} \subset \N$. The embedding is
of course not order-preserving, so $\N$ inherits a new partial
order; denote $\N$ with this new order by $\ZZ_\infty \subset
B(\N)$. Explicitly, the order relations in $\ZZ_\infty$ are $0 \leq
1 \geq 2 \leq 3 \dots$, so that $\ZZ_\infty$ can be visualized as an
infinite zigzag. The map $\xi:B(\N) \to \N$ of \eqref{BJ.sq}
restricts to a map $\zeta:\ZZ_\infty \to \N$, $\zeta(2n-l)=n$, $n
\geq 0$, $l=0,1$ that fits into a cartesian cocartesian square
\begin{equation}\label{N.sq}
\begin{CD}
\oN \times [1] @>>> \oN\\
@VVV @VV{q}V\\
\ZZ_\infty @>{\zeta}>> \N,
\end{CD}
\end{equation}
where as in \eqref{BJ.sq}, $\oN$ is $\N$ with discrete order, and
$q$ sends $n$ to $n+1$.

Moreover, for any $m \geq 0$, we can let $\ZZ_m = \{0,\dots,m \}
\subset \ZZ_\infty$, and then $\zeta$ induces a map $\zeta_m:\ZZ_m
\to [n] \subset \N$ for any $m = 2n-l$, $l=0,1$. In particular, we
have a map $\zeta_3:\ZZ_3 \to [2]$, and \eqref{N.sq} induces a
cartesian cocartesian square
\begin{equation}\label{Z.sq}
\begin{CD}
[1] @>>> \ppt\\
@VVV @VVV\\
\ZZ_3 @>{\zeta_3}>> [2],
\end{CD}
\end{equation}
where the arrow on the right is the embedding onto $1 \in [2] =
\{0,1,2\}$.

All the squares \eqref{S.sq}, \eqref{N.sq}, \eqref{Z.sq} lie in
$\Posf^+ \subset \Pos$. The square \eqref{Z.sq} also lies both in
$\Posf$ and in $\Posf^\pm$, the square \eqref{N.sq} lies in
$\Posf^\pm$, and the square \eqref{S.sq} lies in $\Posf$ but not in
$\Posf^\pm$ (because $S^>$ is not left-finite as soon as $S$ is
infinite). The latter is somewhat inconvenient, so we modify
\eqref{BS.eq} by defining a partially ordered set $S^\hhush$ via the
standard pushout square
\begin{equation}\label{BSS.eq}
\begin{CD}
  S @>>> B(S^>)\\
  @V{\id \times s}VV @VVV\\
  S \times [1] @>>> S^\hhush = S,
\end{CD}
\end{equation}
where the top arrow $S = B(S) \to B(S^>)$ is the left-closed
embedding induced by the standard embedding $\eps:S \to S^>$. Then
$S^\hhush$ is left-finite, and we have a cartesian cocartesian
square
\begin{equation}\label{SS.sq}
\begin{CD}
B(S^>) @>>> \ppt\\
@VVV @VVV\\
S^\hhush @>>> S^<,
\end{CD}
\end{equation}
where the map on the left is the composition of $\id \copr
\xi:S^\hhush \to S^\hush$ and the map $\xi^o:S^\hush \to S^<$ of
\eqref{S.sq}. The square \eqref{SS.sq} lies both in $\Posf$ and in
$\Posf^\pm$.

Now, a useful property of the cartesian cocartesian squares
\eqref{S.sq}, \eqref{N.sq}, \eqref{Z.sq} and \eqref{SS.sq} is that
they produce many other cartesian cocartesian squares by
pullbacks. Namely, for any partially ordered set $J$ equipped with a
map $J \to S^<$, with the fiber $J_o$ over $o \in S^<$, the squares
\eqref{S.sq} and \eqref{SS.sq} produce cartesian squares
\begin{equation}\label{J.S.sq}
\begin{CD}
J_o \times S^> @>>> J_o\\
@VVV @VVV\\
J \times_{S^<} S^\hush @>>> J,
\end{CD}
\qquad\qquad
\begin{CD}
J_o \times B(S^>) @>>> J_o\\
@VVV @VVV\\
J \times_{S^<} S^\hhush @>>> J,
\end{CD}
\end{equation}
and both are also cocartesian. Analogously, for any $J \in \Pos$
equipped with a map $J \to \N$, with the corresponding cofibration
$\tau:J / \N \to \N$, \eqref{N.sq} produces a cartesian cocartesian
square
\begin{equation}\label{J.N.sq}
\begin{CD}
\coprod_{n \geq 1} (J/n) \times [1] @>>> \coprod_{n \geq 1} J/n\\
@VVV @VVV\\
\zeta^*(J / \N) @>>> J / \N.
\end{CD}
\end{equation}
For the square \eqref{Z.sq}, one has to be more careful: simply
taking some map $J \to [2]$ would produce a cartesian square that is
not necessarily cocartesian. However, one can consider the arrow set
$\Ar([1]) \cong [1]/[1]$, and it is easy to see that it is totally
ordered, thus can be canonically identified with $[2]$; under this
identification, the projection $\tau:[1] / [1] \to [1]$ sends $0 \in
[2]$ to $0 \in [1]$ and $1,2 \in [2]$ to $1 \in [1]$. Then for any
partially ordered set $J$ equipped with a map $\chi:J \to [1]$, we
have the comma-set $J/[1]$ equipped with the map $\chi/\id:J/[1] \to
[1]/[1] \cong [2]$, and \eqref{Z.sq} produces a cartesian
cocartesian square
\begin{equation}\label{J.Z.sq}
\begin{CD}
J_0 \times [1] @>>> J_0\\
@VVV @VVV\\
\zeta_3^*(J/[1]) @>>> J/[1],
\end{CD}
\end{equation}
where $J_0 = \chi^{-1}(0) \subset J$ is the fiber of the map $\chi$
over $0 \in [1]$.

\subsection{Nerves and Segal spaces.}\label{del.subs}

As usual, we let $\Delta \subset \Pos$ be the full subcategory
spanned by $[n]$, $n \geq 0$. Equivalently, $\Delta$ is the category
of all finite non-empty totally ordered sets. A {\em simplicial
  object} in a category $\C$ is a functor $\Delta^o \to \C$. For any
category $I$, we can consider the preimage $\Delta \bb I \subset
\Cat \bb I$ of $\Delta \subset \Cat$ under the fibration $\Cat \bb I
\to \Cat$, and we have the subcategory $\Delta \bbi I \subset \Delta
\bb I$ formed by maps cartesian over $\Delta$. Then $\Delta \bbi I
\to \Delta$ is a family of groupoids; explicitly, its fiber over
$[n] \in \Delta$ is the groupoid of functors $[n] \to I$ and
isomorphisms between these functors. To recover $I$ from $\Delta
\bbi I$, say that a map $f:[m] \to [n]$ in $\Delta$ is {\em special}
if $f(m) = n$ -- or equivalently, if $f$ is left-reflexive -- and
for any fibration $\pi:\C \to \Delta$, let $+$ be the class of maps
$f$ in $\C$ with special $\pi(f)$. Then $I \cong h^+(\Delta \bbi I)$
(for a proof of this, see e.g.\ \cite[Lemma 4.2.1.1]{big}).

One can also characterize families of groupoids $\C \to \Delta$ that
arise as $\Delta \bbi I$ for some $I$. To do this, say that $\C$ is a
{\em Segal family} if it is cartesian over the square \eqref{seg.sq}
for any $n > l > 0$. Then let $\mu_l:[1] \to [3]$, $l = 0,1$ be the
embeddings given by $\mu_l(a) = 2a+l$, consider the commutative
square
\begin{equation}\label{compl.sq}
\begin{CD}
\C_{[0]} @>{\delta}>> \C_{[0]} \times \C_{[0]}\\
@V{e^*}VV @VV{e^* \times e^*}V\\
\C_{[3]} @>{\mu_0^* \times \mu_1^*}>> \C_{[1]} \times \C_{[1]},
\end{CD}
\end{equation}
where $\delta$ is the diagonal embedding, and say that the family
$\C$ is {\em complete} if the square \eqref{compl.sq} is
cartesian. In these terms, a family of groupoids $\C \to \Delta$ is
a complete Segal family if and only if $\C \cong \Delta \bbi I$ for
some category $I$, and if this holds, then we have a canonical
equivalence $I \cong h^+(\C)$ (for a proof, see
e.g.\ \cite[Proposition 4.2.3.3]{big}).

\begin{remark}
The square \eqref{compl.sq} is essentially induced by a commutative
square
\begin{equation}\label{compl.co.sq}
\begin{CD}
[1] \copr [1] @>{\mu_0 \copr \mu_1}>> [3]\\
@V{e \copr e}VV @VV{e}V\\
[0] \copr [0] @>>> [0]
\end{CD}
\end{equation}
of finite partially ordered sets, and this is an example of a
cocartesian square in $\Pos$ that is not cocartesian in
$\Cat$. Neither is it a cocartesian commutative square of categories
in the sense of Subsection~\ref{cat.subs}; to turn it into such a
cocartesian square, one has to replace $\ppt = [0]$ in the bottom
right conner with an equivalent category $e(\{0,1\})$ (where for any
set $S$, we let $e(S)$ be the small category whose objects are
elements $s \in S$, and that has exactly one morphism between any
two objects).
\end{remark}

If the category $I$ is small, one can also construct a small
discrete fibration $\Delta I \to \Delta$ called the {\em simplicial
  replacement} for $I$. To do this, define the {\em nerve functor}
$N$ by
\begin{equation}\label{nerve.eq}
  N = \phi_*\Y:\Cat \to \Delta^o\Sets,
\end{equation}
where $\Y$ is the Yoneda embedding \eqref{yo.eq} for $\Delta$, and
$\phi:\Delta \to \Cat$ is the standard full embedding, and let
$\Delta I = \Delta N(I)$. Then explicitly, if one computes the right
Kan extension in \eqref{nerve.eq} by \eqref{kan.eq}, one observes
that $N(I)([n])$ is the set of functors $[n] \to I$, and then
$\Delta I \to \Delta$ is the discrete fibration whose fibers
$(\Delta I)_{[n]}$ are these sets. We have a natural embedding
$\Delta I \to \Delta \bbi I$ that is bijective on objects, but
$\Delta \bbi I$ has more morphisms: we allow isomorphisms between
functors $[n] \to I$. We still have $I = h^+(\Delta I)$, so passing
to $\Delta I$ loses no information -- in fact, \eqref{nerve.eq} is a
fully faithful embedding -- but $\Delta \bbi I$ has better
functoriality with respect to equivalences: an equivalence $\gamma:I
\to I'$ induces an equivalence $\Delta \bbi I \to \Delta \bbi I'$,
while the induced functor $\Delta I \to \Delta I'$ is only an
equivalence if $\gamma$ is an isomorphism. The price to pay is that
$\Delta \bbi I$ is only a family of groupoids, and does not
correspond to a simplicial set. To remedy this, one can observe that
$\Delta \bbi I$ is actually rigid enough to correspond to a
simplicial object in $\Cat$, and then again take the nerve. The
result is a {\em bisimplicial set}, that is, a functor
$N^2(I):\Delta^o \to \Delta^o\Sets$, given by
\begin{equation}\label{binerve.eq}
N^2(I)([n]) = N(\Fun([n],I)_\Iso), \qquad [n] \in \Delta,
\end{equation}
where $\Fun([n],I)_\Iso \subset \Fun([n],I)$ is the isomorphism
groupoid of the small category $\Fun([n],I)$. Then $N^2 I \in
\Delta^o\Delta^o\Sets \cong (\Delta \times \Delta)^o\Sets$ defines a
discrete fibration $\pi:\Delta^2 I \to \Delta^2 = \Delta \times
\Delta$, and $I \cong h^E(\Delta^2 I)$, where $E=+ \times \Tot$ is
the class of maps $f$ such that $\pi(f) = f_0 \times f_1$ is the
product of a special map $f_0$ in the first factor $\Delta$, and any
map $f_1$ in the second one.

A beautiful idea of Ch.\ Rezk \cite{rzk} is that, roughly speaking,
it is the double nerve \eqref{binerve.eq} that should be generalized
to a homotopical setting. Namely, $\Delta^o\Sets$ has the standard
model structure of \cite{qui.ho} where in particular, nerves of
small groupoids are fibrant, and a functor $\gamma:I \to I'$ between
small groupoids is an equivalence if and only if $N(\gamma):N(I) \to
N(I')$ is a weak equivalence. Moreover, since $\Delta$ is a Reedy
category (\cite{reedy}, \cite{hovey}), $\Delta^o\Delta^o\Sets$
acquires a Reedy model structure, with weak equivalences given by
pointwise weak equivalences. Recall that a commutative square in a
model category is {\em homotopy cartesian} if it is weakly
equivalent to a cartesian square of fibrant object and
fibrations. Recall also that $\Delta^o\Delta^o\Sets \cong (\Delta
\times \Delta)^o\Sets$ is cartesian-closed, and for any $X \in
\Delta^o\Delta^o\Sets$ and partially ordered set $J$, denote
\begin{equation}\label{X.n.2}
X(J)^{(2)} = \Hhom_{\Delta^o\Delta^o\Sets}(N^2(J),X) \in
\Delta^o\Delta^o\Sets.
\end{equation}
Note that since partially ordered sets $J$ are rigid, their double
nerves $N^2(J)$ are constant along the second factor $\Delta$, so
that $N^2(J) \cong \pi_1^{o*}N(J)$, where $\pi_1:\Delta \times
\Delta \to \Delta$ is the projection onto the first factor. In
particular, if $X$ is constant along the first factor, then
$X([n])^{(2)} \cong X([0])^{(2)}$ for any $n \geq 0$.

\begin{defn}\label{css.def}
A {\em Segal space} is a fibrant bisimplicial set $X:\Delta^o \to
\Delta^o\Sets$ such that for any $n > l > 0$, the square
$$
\begin{CD}
X([n]) @>>> X([l])\\
@VVV @VVV\\
X([n-l]) @>>> X([0])
\end{CD}
$$
induced by \eqref{seg.sq} is homotopy cartesian. A Segal space
$X$ is {\em complete} if the square
$$
\begin{CD}
X([3])^{(2)} @>>> X([1])^{(2)} \times X([1])^{(2)}\\
@VVV @VVV\\
X([0])^{(2)} @>>> X([0])^{(2)} \times X([0])^{(2)}
\end{CD}
$$
induced by \eqref{compl.co.sq} is homotopy cartesian as well.
\end{defn}

Then in particular, for any small category $I$, its double nerve
$N^2(I)$ of \eqref{binerve.eq} is a complete Segal space, and a
functor $\gamma:I \to I'$ is an equivalence iff $N^2(\gamma):N^2(I)
\to N^2(I')$ is a weak equivalence of bisimplicial sets. More
generally, for any complete Segal space $X$, one can consider the
discrete fibration $\Delta^2 X \to \Delta^2$, and define its {\em
  truncation} by $h^{+ \times \Tot}(\Delta^2 X)$. The space $X$
itself -- considered up to a weak equivalence, of course -- provides
a homotopical enhancement for its truncation, and as far as
enhancements based on model category techniques are concerned, this
is as good as it gets.

\section{Enhanced categories.}\label{enh.sec}

\subsection{Reflexive families.}

We are now ready to give our main definitions. We start with the
notion of a ``reflexive family'' -- this is a somewhat technical
gadget that nevertheless allows one to prove something. Fix a full
subcategory $\I \subset \Pos$ ample in the sense of
Definition~\ref{amp.def}. By a {\em family of categories} over $\I$
we will understand a fibration $\C \to \I$. For any such family $\C$
and $J \in \I$, we have a functor
\begin{equation}\label{nondege.eq}
\prod_{j \in J}\eps(j)^*:\C_J \to \prod_{j \in J}\C_\ppt,
\end{equation}
where $\eps(j):\ppt \to J$ is the embedding onto $j \in J$. Say that
the family $\C$ is {\em non-degenerate} if all the functors
\eqref{nondege.eq} are conservative. Moreover, for any $J$, let
$s,t:J \to J \times [1]$ be the embeddings onto $J \times \{0\}$, $J
\times \{1\}$, and let $e:J \times [1] \to J$ be the projection onto
the first factor, so that $e \circ s = e \circ t = \id$.

\begin{defn}\label{relf.def}
A family of categories $\C \to \I$ is {\em reflexive} if for any $J
\in \I$, the isomorphism $s^* \circ e^* \cong \id$ induced by $e
\circ s = \id$ defines an adjunction between the transition functors
$e^*:\C_J \to \C_{J \times [1]}$ and $s^*:\C_{J \times [1]} \to
\C_J$.
\end{defn}

For any reflexive family of categories $\C \to \I$ and $J \in \I$,
$c \in \C_{J \times [1]}$, the second adjunction map $a(c):c \to
e^*s^*c$ is functorial with respect to $c$, and sending $c$ to the
arrow $t^*a(c)$ in $\C_J$ provides a functor
\begin{equation}\label{nu.J.eq}
\nu_J:\C_{J \times [1]} \to \Ar(\C_J).
\end{equation}
The family $\C$ is {\em separated} if \eqref{nu.J.eq} is an
epivalence for any $J \in \I$.

For any category $\E$, we denote $\Unf(\I,\E) = \iota(\I) \bb \E$,
where $\iota(\I) \subset \Pos$ is the image of $\I \subset \Pos$
under the involution $\iota:\Pos \to \Pos$, $J \mapsto
J^o$. Explicitly, objects in the category $\Unf(\I,\E)$ are pairs
$\langle J,\alpha \rangle$ of $J \in \I$ and a functor $\alpha:J \to
\E^o$, with morphisms $\langle J,\alpha \rangle \to \langle
J',\alpha' \rangle$ given by pairs of a map $f:J \to J'$ and a
morphism $f^*\alpha' \to \alpha$. Then the forgetful functor
$\Unf(\I,\E) \to \I$, $\langle J,\alpha \rangle \mapsto J$ is a
fibration with fibers $\Unf(\I,\E)_J \cong J^o\E$, so that
$\Unf(\I,\E)$ is a family of categories over $\I$, and this family
is non-degenerate, reflexive and separated (for the latter, note
that \eqref{nu.J.eq} for $\Unf(\I,\E)$ is an equivalence). A functor
$\gamma:\E \to \E'$ to some category $\E'$ induces a functor
$\Unf(\gamma):\Unf(\I,\E) \to \Unf(I,\E')$, cartesian over $\I$.

\begin{prop}\label{trunc.prop}
For any reflexive family $\C \to \I$, there exists a {\em truncation
  functor}
\begin{equation}\label{trunc.eq}
\unf(\C):\C \to \Unf(\I,\C_\ppt),
\end{equation}
cartesian over $\I$, such that for any category $\E$ and
functor $\gamma:\C \to \Unf(\I,\E)$ cartesian over $I$, we have a
unique isomorphism $\gamma \cong \Unf(\gamma_\ppt) \circ \unf(\C)$
that restricts to $\id$ over $\ppt \in \I$.
\end{prop}

\proof{} This is part of \cite[Proposition 7.1.2.1]{big}. \endproof

In particular, the fiber of the truncation functor \eqref{trunc.eq}
over some $J \in \I$ is a functor
\begin{equation}\label{tru.J}
  \C_J \to J^o\C_\ppt,
\end{equation}
and we also have the functors \eqref{nu.J.eq}. One can actually
define an even more general family of functors that includes both
\eqref{nu.J.eq} and \eqref{tru.J}. Namely, let $\E$ be the category
$\I$ itself. Then \eqref{arc.I.eq} provides a fully faithful
embedding
\begin{equation}\label{arf.eq}
\Arf(\I) \to \Unf(\I,\I),
\end{equation}
where $\Arf(\I) \subset \Ar(\I)$ consists of fibrations $J' \to J$,
with maps given by squares \eqref{arc.sq} such that $g':J_0' \to
J_1'$ is cartesian over $g:J_0 \to J_1$. The projection
$\tau:\Arf(\I) \to \I$ is a fibration that turns $\Arf(\I)$ into a
reflexive family of $\I$, and the full embedding \eqref{arf.eq} is
cartesian over $\I$. However, we also have the projection
$\sigma:\Arf(\I) \to \I$.

\begin{lemma}\label{C.refl.le}
For any reflexive family $\pi:\C \to \I$, the category $\sigma^*\C$
with the projection $\tau \circ \sigma^*(\pi):\sigma^*\C \to
\Arf(\I) \to \I$ is a reflexive family.
\end{lemma}

\proof{} This is \cite[Lemma 7.1.2.6]{big} combined with
\cite[Example 7.1.1.8]{big}. \endproof

By virtue of Lemma~\ref{C.refl.le}, we can apply
Proposition~\ref{trunc.prop} to the family $\sigma^*\C$ and consider
the corresponding truncation functor \eqref{trunc.eq} and its
components \eqref{tru.J}. Explicitly, for any fibration $J' \to J$
in $\I$, we obtain a functor
\begin{equation}\label{tru.J.J}
\C_{J'} \to \Sec(J,(j_\idot^*\C)_\perp),
\end{equation}
where $j_\idot:J^o \to \I$ corresponds to $J' \to J$ under
\eqref{arf.eq}, and $(j_\idot^*\C)_\perp \to J$ is the transpose
cofibration. When $J' = J$, this gives \eqref{tru.J}, and the
fibration $J \times [1] \to [1]$ recovers \eqref{nu.J.eq}.

It is also useful to restrict the family $\sigma^*\C$ of
Lemma~\ref{C.refl.le} to various subcategories in $\Arf(\I)$. In
particular, for any $J \in \I$, we have the embedding $j:\I \to
\Arf(\I)$ sending $J' \in \I$ to the projection $J' \times J \to J'$
(in terms of \eqref{arf.eq}, this corresponds to equipping $J'$ with
the constant functor $\alpha:{J'}^o \to \I$ with value $J$). If
we let $J^o_h\C = j^*\sigma^*\C$, then we have the following result.

\begin{lemma}\label{J.C.le}
For any $J \in \I$ and reflexive family $\C \to \I$, $J^o_h\C \to
\I$ is a reflexive family over $\I$, and there exists a functor
\begin{equation}\label{j.eq}
\ev:\Unf(\I,J^o) \times_\I J^o_h\C \to \C,
\end{equation}
cartesian over $\I$, with the following universal property: for any
reflexive family $\C' \to \I$, a functor $\gamma:\Unf(\I,J^o)
\times_\I \C' \to \C$ cartesian over $\I$ factors as
\begin{equation}\label{J.o.C}
\begin{CD}
\Unf(\I,J^o) \times_\I \C' @>{\id \times \wt{\gamma}}>>
\Unf(\I,J^o) \times_\I J^o_h\C @>{\ev}>> \C,
\end{CD}
\end{equation}
for a certain functor $\wt{\gamma}:\C' \to J^o_h\C$, cartesian over
$\I$ and unique up to a unique isomorphism.
\end{lemma}

\proof{} This is \cite[Corollary 7.1.2.7]{big} (with
\cite[Example 7.1.1.8]{big}). \endproof

\subsection{Enhanced categories and functors.}

Now recall that we have a collection of cartesian cocartesian
squares \eqref{J.S.sq}, \eqref{J.N.sq}, \eqref{J.Z.sq} in $\Pos$. At
the insistence of the referee, and for the convenience of the
reader, here they are again:
$$
\begin{aligned}
&\begin{CD}
J_o \times S^> @>>> J_o\\
@VVV @VVV\\
J \times_{S^<} S^\hush @>>> J,
\end{CD}
\qquad\qquad
\begin{CD}
J_o \times B(S^>) @>>> J_o\\
@VVV @VVV\\
J \times_{S^<} S^\hhush @>>> J,
\end{CD}\\[4mm]
&\begin{CD}
\coprod_{n \geq 1} (J/n) \times [1] @>>> \coprod_{n \geq 1} J/n\\
@VVV @VVV\\
\zeta^*(J / \N) @>>> J / \N,
\end{CD}
\qquad
\begin{CD}
J_0 \times [1] @>>> J_0\\
@VVV @VVV\\
\zeta_3^*(J/[1]) @>>> J/[1],
\end{CD}
\end{aligned}
$$
To distinguish between the two squares on the left resp.\ right in
\eqref{J.S.sq}, let us denote them by \eqrefi{J.S.sq}{i}
resp.\ \eqrefi{J.S.sq}{ii}. Recall that we also have standard
pushout squares \eqref{st.sq}.

\begin{defn}\label{exci.def}
Assume given a family of categories $\C \to \I$ over an ample full
subcategory $\I \subset \Pos$. Then $\C$ is {\em additive} if for
any $J \in \I$ equipped with a map $J \to S$ to a discrete $S \in
\Sets \subset \Pos$, the product
$$
\C_J \to \prod_{s \in S}\C_{J_s}
$$
of transition functors $\C_J \to \C_{J_s}$ for the embeddings $J_s
\to J$ is an equivalence. The family $\C$ is {\em semiexact} if it
is semicartesian over any standard pushout square \eqref{st.sq} in
$\I$. The family $\C$ {\em satifies excision} resp.\ {\em modified
  excision} if if it is cartesian over any square
\eqrefi{J.S.sq}{i} resp.\ \eqrefi{J.S.sq}{ii} in
$\I$. The family $\C$ is {\em semicontinuous} resp.\ {\em satisfies
  the cylinder axiom} if it is cartesian over any square
\eqref{J.N.sq} resp.\ \eqref{J.Z.sq} in $\I$.
\end{defn}

\begin{defn}\label{enh.def}
An {\em enhanced category} is a non-degenerate separated reflexive
family of categories $\C \to \Posf^+$ that is additive, semiexact,
semicontinuous and satisfies excision and the cylinder axiom. An
enhanced category $\C$ is {\em small} if so is the fibration $\C \to
\Posf^+$, and a small enhanced category is {\em $\kappa$-bounded},
for a regular cardinal $\kappa$, if $|\C_J| < \kappa$ for any finite
$J \in \Posf^+$. An {\em enhanced functor} between enhanced
categories $\C$, $\C'$ is a functor $\gamma:\C \to \C'$ cartesian
over $\Posf^+$, and an {\em enhanced map} between enhanced functors
$\gamma$, $\gamma'$ is a map $\gamma \to \gamma'$ over $\Posf^+$.
\end{defn}

For any enhanced categories $\C$, $\C'$, their coproduct $\C \copr
\C'$ is an enhanced category, and so is their product over $\Posf^+$
that we denote by
\begin{equation}\label{enh.prod.eq}
\C \times^h \C' = \C \times_{\Posf^+} \C'.
\end{equation}
For any category $\E$, the reflexive family $\Unf(\E) =
\Unf(\Posf^+,\E)$ is an enhanced category (\cite[Example
  7.2.1.16]{big}). We recall that explicitly, $\Unf(\E) \to \Posf^+$
is a fibration with fibers $\Unf(\E)_J \cong J^o\E$, $J \in
\Posf^+$. For any functor $\gamma:\E' \to \E$, the functor
$\Unf(\gamma):\Unf(\E') \to \Unf(\E)$ is an enhanced functor
(\cite[Example 7.2.2.2]{big}), and all enhanced functors $\Unf(\E')
\to \Unf(\E)$ are of this form, for some $\gamma:\E \to \E'$ unique
up to a unique isomorphism (\cite[Proposition
  7.1.2.1]{big}). Enhanced maps $\Unf(\gamma) \to \Unf(\gamma')$
correspond bijectively to maps $\gamma \to \gamma'$. If we have a
cartesian square
\begin{equation}\label{U.prod}
\begin{CD}
\C' @>{\wt{\gamma}}>> \C\\
@V{\pi'}VV @VV{\pi}V\\
\Unf(\E') @>{\Unf(\gamma)}>> \Unf(\E),
\end{CD}
\end{equation}
where $\C$ is an enhanced category and $\pi$ is an enhanced functor,
then $\C'$ is an enhanced category, and $\pi'$, $\wt{\gamma}$ are
enhanced functors (\cite[Lemma 7.2.2.4]{big}). In general, we think
of an enhanced category $\C$ as providing an ``enhancement'' for its
``truncation'' $\C_\ppt$; if $\pi$ in \eqref{U.prod} is the
truncation functor \eqref{trunc.eq}, then we see that any category
$\E'$ equipped with a functor $\E' \to \E = \C_\ppt$ inherits an
enhancement. In particular, say that a class $v$ of morphisms in a
category $I$ is {\em closed} if it is closed under compositions and
contains all the identity maps, and for any closed class $v$, denote
by $I_v \subset I$ the subcategory with the same objects as $I$ and
only those morphisms that are in $v$. Then for any closed class $v$
of morphisms in $\C_\ppt$, \eqref{U.prod} provides an enhanced
category
\begin{equation}\label{h.v.eq}
\C_{hv} = \Unf(\C_{\ppt, v}) \times_{\Unf(\C_\ppt)} \C.
\end{equation}
For any enhanced category $\C$, the reflexive family $\sigma^*\C \to
\Posf^+$ of Lemma~\ref{C.refl.le} is an enhanced category
(\cite[Lemma 7.2.4.4]{big}); since we have $(\sigma^*\C)_\ppt \cong
\C$, this means that $\C$ itself treated simply as a category also
inherits a canonical enhancement. The terminal enhanced category is
$\Unf(\ppt) \cong \Posf^+$, and we denote it by $\ppt^h =
\Unf(\ppt)$. An {\em enhanced object} in an enhanced category is an
object $c \in \C_\ppt$ in the truncation $\C_\ppt$; these correspond
bijectively to enhanced functors $\eps^h(c):\ppt^h \to \C$. An {\em
  enhanced morphism} between enhanced objects $c$, $c'$ is an
enhanced functor $f:\Unf([1]) \to \C$ equipped with isomorphisms
$\Unf(s)^*f \cong c$, $\Unf(t)^*f \cong c'$. Alternatively, for any
$J \in \Posf^+$ and enhanced category $\C$, the reflexive family
$J^o_h\C$ of Lemma~\ref{J.C.le} is an enhanced category
(\cite[Corollary 7.2.4.6]{big}); we define the {\em enhanced arrow
  category} $\Ar^h(\C)$ as $\Ar^h(\C) \cong [1]^o_h\C$, with the
enhanced functors
\begin{equation}\label{s.t.e.ar.h}
\eta = e^*:\C \to \Ar^h(\C), \qquad \sigma = t^*,\tau =
s^*:\Ar^h(\C) \to \C
\end{equation}
induced by the projection $e:[1] \to [0] = \ppt$ and embeddings
$s,t:[0] \to [1]$. Then enhanced morphisms in $\C$ correspond to
enhanced objects in $\Ar^h(\C)$, or equivalently, to objects in
$\C_{[1]} \cong \Ar^h(\C)_\ppt$. Since $\C$ is separated,
\eqref{nu.J.eq} for $J=\ppt$ is an epivalence, so that isomorphism
classes of enhanced morphisms $c \to c'$ are the same thing as
morphisms $c \to c'$ in $\C_\ppt$; however, enhanced morphisms form
a groupoid, not a set, so two of them cannot be equal.

An enhanced functor $\C' \to \C$ is {\em fully faithful} if it is
fully faithful as a functor. Then $\C'_\ppt \to \C_\ppt$ is also
fully faithful, and $\C'$ fits into a cartesian square
\eqref{U.prod}, where $\pi$ is the truncation functor
\eqref{trunc.eq}, and $\gamma:\C'_\ppt \to \C_\ppt$ is the embedding
(\cite[Corollary 7.2.2.18]{big}). In words, full enhanced
subcategories $\C' \subset \C$ correspond bijectively to full
subcategories $\C'_\ppt \subset \C_\ppt$, or equivalently, to
collections of enhanced objects in $\C$. In other words, the concept
of a full embedding is simply inherited from the usual category
theory, with no changes. Another such concept is that of adjunction.

\begin{defn}\label{enh.refl.def}
An enhanced functor $\gamma:\C \to \C'$ is {\em left} resp.\ {\em
  right-reflexive} if it is left resp.\ right-reflexive over
$\Posf^+$, and the adjoint functor $\gamma_\dg$ resp. $\gamma^\dg$
is also cartesian over $\Posf^+$.
\end{defn}

Explicitly, an enhanced adjunction between two given enhanced
functors $\lambda:\C \to \C'$, $\rho:\C' \to \C$ is defined by
enhanced maps $\lambda \circ \rho \to \id$, $\id \to \rho \circ
\lambda$ such that the compositions \eqref{adj.eq} are identity
maps. If this happens, $\lambda$ is left-reflexive, $\rho$ is
right-reflexive, and we have canonical isomorphism $\rho \cong
\lambda_\dg$ and $\lambda \cong \rho^\dg$. For example, for any
enhanced category $\C$, $\sigma$ resp.\ $\tau$ in \eqref{s.t.e.ar.h}
is right resp.\ left-adjoint to $\eta$, with adjunction defined by
isomorphisms $\tau \circ \eta \cong \id \cong \sigma \circ \eta$; by
adjunction, $\eta$ is a fully faithful embedding, and the other
adjunction maps $\eta \circ \sigma \to \id \to \eta \circ \tau$
define an enhanced map
\begin{equation}\label{ar.map}
  \eta \circ \sigma \to \eta \circ \tau
\end{equation}
between enhanced functors $\eta \circ \sigma,\eta \circ
\tau:\Ar^h(\C) \to \C \to \Ar^h(\C)$.

Since for any non-degenerate reflexive family $\C \to \Posf^+$, the
truncation functor \eqref{trunc.eq} is conservative, an enhanced
category $\C \to \Posf^+$ is a family of groupoids if and only if
$\C_\ppt$ is a groupoid. In this case, we say that $\C$ is an {\em
  enhanced groupoid}. A family of groupoids $\C \to \Posf^+$ is
trivially non-degenerate, and it is reflexive iff it is constant
over $e:J \times [1] \to [1]$ for any $J \in \Posf^+$, so it is
automatically separated. The cylinder axiom of
Definition~\ref{exci.def} is also automatic: a reflexive family of
groupoids $\C \to \Posf^+$ is an enhanced groupoid iff it is
additive, semiexact, semicontinuous, and satisfies excision.  For
any enhanced category $\C$, \eqref{h.v.eq} for the closed class $v =
\Iso$ of all invertible maps in $\C_\ppt$ gives the {\em enhanced
  isomorphism groupoid}
\begin{equation}\label{h.iso.eq}
\C_{h\Iso} = \Unf(\C_{\ppt,\Iso}) \times_{\Unf(\C_\ppt)} \C
\end{equation}
of the enhanced category $\C$. The truncation $\C_{h\Iso,\ppt} \cong
\C_{\ppt,\Iso}$ is the isomorphism groupoid of the truncation
$\C_\ppt$, and the whole enhanced groupod \eqref{h.iso.eq} has the
same universal property: a enhanced functor $\E \to \C$ from an
enhanced groupoid $\E$ factors through $\C_{h\Iso} \to \C$, uniquely
up to a unique isomorphism.

An enhanced object $o:\Posf^+ \to \C$ in an enhanced category $\C$
is {\em initial} resp.\ {\em terminal} if $\eps^h(o)$ is right
resp.\ left-reflexive in the sense of
Definition~\ref{enh.refl.def}. Note that the embedding functors $s =
\eps^h(0),t = \eps^h(1):\ppt^h \to \Unf([1])$ corresponding to
objects $0,1 \in [1]$ admit right-adjoints $s^\dg,t^\dg:\Unf([1])
\to \ppt^h$ sending $\langle J,\alpha \rangle \in \Unf([1])$ to
$s^*J$ resp.\ $t^*J$, and then for any enhanced category $\C$, we
can define families of categories
\begin{equation}\label{h.gt.lt}
\C^{h>} = s^{\dg*}\C \to \Unf([1]) \to \Posf^+, \qquad \C^{h<} =
t^{\dg*}\C \to \Unf([1]) \to \Posf^+.
\end{equation}
Both are enhanced categories (\cite[Example 7.2.2.16]{big}), we have
identifications $\C^{h>}_\ppt \cong (\C_\ppt)^>$, $\C^{h<}_\ppt
\cong (\C_\ppt)^<$, and the enhanced object $o$ in $\C^{h>}$
resp.\ $\C^{h<}$ corresponding to $o \in \C^>_\ppt$
resp.\ $\C^<_\ppt$ is terminal resp.\ initial. The adjunction maps
$s \circ s^\dg,t \circ t^\dg \to \id$ then induce enhanced functors
\begin{equation}\label{h.gt.lt.1}
\C \times^h \Unf([1]) \to \C^{h>}, \qquad \C \times^h \Unf([1]) \to
\C^{h<},
\end{equation}
an enhanced version of the bottom arrows in \eqref{cyl.sq}.

\subsection{Extensions and restrictions.}

The main reason we have chosen $\Posf^+$ as the domain of definition
for enhanced categories is that this minimizes the axioms one need
to impose in Definition~\ref{enh.def}. However, it is actually
sufficient to define enhanced categories over $\Posf^\pm \subset
\Posf^+$, and one can canonically extend them to families of
categories over the whole $\Pos$.

\begin{defn}
A family of categories over $\Pos$ is {\em bar-invariant} if it is
cartesian over the square \eqref{BJ.sq} for any $J \in \Pos$.
\end{defn}

\begin{prop}\label{lf.prop}
The restriction of an enhanced category $\C \to \Posf^+$ to
$\Posf^\pm \subset \Posf^+$ satisfies modified excision in the sense
of Definition~\ref{exci.def}. Conversely, any non-degenerate
separated reflexive family $\C \to \Posf^+$ that is additive,
semiexact, semicontinuous, and satisfies the cylinder axiom and
modified excision extends to a bar-invariant family of categories
$\C^\Tot \to \Pos$ that is non-degenerate, reflexive, separated,
additive, semiexact, semicontinuous, and satisfies excision and the
cylinder axiom, so that its restriction to $\Posf^+ \subset \Pos$ is
an enhanced category. Moreover, for any two such families $\C_0,\C_1
\to \Posf^\pm$ with such extensions $\C_0^\Tot,\C_1^\Tot \to \Pos$,
a functor $\gamma:\C_0 \to \C_1$ cartesian over $\Posf^\pm$ extends
to a functor $\gamma^\Tot:\C_0^\Tot \to \C_1^\Tot$ cartesian over
$\Pos$, and a map $f:\gamma_0 \to \gamma_1$ over $\Posf^\pm$ between
two such functors uniquely extends to a map $f^\Tot:\gamma^\Tot_0
\to \gamma^\Tot_1$ over $\Pos$.
\end{prop}

\proof{} This is \cite[Proposition 7.2.3.7]{big} and \cite[Lemma
  7.2.3.6]{big}. \endproof

One immediate application of canonical extensions of
Proposition~\ref{lf.prop} is the definition of an opposite enhanced
category.

\begin{prop}\label{iota.prop}
For any enhanced category $\C \to \Posf^+$, with canonical extension
$\C^\Tot \to \Pos$, the fibration $\iota^*(\C^\Tot)^o_\perp \to
\Pos$ is a canonical extension of an enhanced category $\C^\iota \to
\Posf^+$.
\end{prop}

\proof{} This is \cite[Lemma 7.2.3.11]{big}. \endproof

\begin{defn}
The enhanced category $\C^\iota$ of Proposition~\ref{iota.prop} is
the {\em enhanced opposite} of the enhanced category $\C$.
\end{defn}

By virtue of the functoriality part of Proposition~\ref{lf.prop}, an
enhanced functor $\gamma:\C_0 \to \C_1$ induces a canonical {\em
  enhanced-opposite functor} $\gamma^\iota:\C_0^\iota \to
\C_1^\iota$, and any enhanced map $\gamma_0 \to \gamma_1$ induces an
enhanced map $\gamma_1^\iota \to \gamma_0^\iota$. An enhanced
functor $\gamma$ is left resp.\ right-reflexive iff its
enhanced-opposite $\gamma^\iota$ is right resp.\ left-reflexive, and
$(\gamma_\dg)^\iota \cong (\gamma^\iota)^\dg$, $(\gamma^\dg)^\iota
\cong (\gamma^\iota)_\dg$. We have $\Unf(\E)^\iota \cong \Unf(\E^o)$
for any $\E$, and $\Unf(\gamma)^\iota \cong \Unf(\gamma^o)$ for any
functor $\gamma:\E \to \E'$.

While Proposition~\ref{lf.prop} requires a proof, the extension
procedure itself is very straightforward and dictated by
bar-invariance: since for any partially ordered set $J$, both $B(J)$
and $\oB(J)$ are left-finite, and so is $\overline{J}$, $\C^\Tot$ is
simply given by the cartesian square
\begin{equation}\label{C.tot}
\begin{CD}
\C^\Tot @>>> S^*\C\\
@VVV @VVV\\
B^*\C @>>> \oB^*\C,
\end{CD}
\end{equation}
where the functor $S:\Pos \to \Sets \subset \Pos$ sends $J$ to the
underlying discrete set $\overline{J}$. One can then describe the
enhanced-opposite category $\C^\iota$ purely in terms of $\C$: it is
given by the cartesian square
\begin{equation}\label{C.iota}
\begin{CD}
\C^\iota @>>> (S \circ \iota)^*\C^o_\perp\\
@VVV @VVV\\
(B \circ \iota)^*\C^o_\perp @>>> (\oB \circ \iota)^*\C^o_\perp
\end{CD}
\end{equation}
induced by \eqref{BJ.sq} for $J^o$.

Another extension result (\cite[Proposition 7.1.5.2]{big}) says that
an enhanced category $\C \to \Posf^+$ is completely determined by
the subcategory $\C_\flat \subset \C$ with the same objects and maps
cartesian over $\Posf^+$. By definition, $\C_\flat \to \Posf^+$ is a
family of groupoids. If $\C$ is an enhanced groupoid, then $\C_\flat
= \C$. In general, $\C \supset \C_\flat$ is bigger, but one can
still recover the whole $\C$ from the family of groupoids
$\C_\flat$, and there is also a characterization of families of
groupoids that appear in this way (in particular, they are additive,
semiexact, semicontinuous, and satisfy excision and the cylinder
axiom in the sense of Definition~\ref{exci.def}). We do not
reproduce this result here since it is technical rather than
genuinely useful; however, here is the corresponding reconstruction
result for enhanced functors.

\begin{lemma}\label{fl.le}
For any enhanced categories $\C,\C' \to \Posf^+$, with the
underlying families of groupoids $\C_\flat,\C'_\flat \to \Posf^+$ of
cartesian maps, a functor $\gamma_\flat:\C_\flat \to \C'_\flat$ over
$\Posf^+$ extends to an enhanced functor $\gamma:\C \to \C'$,
uniquely up to a unique isomorphism.
\end{lemma}

\proof{} This is \cite[Lemma 7.2.2.9]{big}. \endproof

Finally, we note that the extension $\C^\Tot$ given by
Proposition~\ref{lf.prop} is really canonical, in that it actually
recovers functors up to a unique isomorphism. If one is only
interested in small enhanced categories and isomorphism classes of
functors, then there is a further extension result that allows one
to only consider families over $\Posf$ and drop one more axiom.

\begin{defn}\label{restr.enh.def}
A {\em restricted enhanced category} is an additive semiexact
non-degenerate reflexive separated small fibration $\C \to \Posf$
that satisfies excision and the cylinder axiom. A {\em restricted
  enhanced functor} between restricted enhanced categories $\C$,
$\C'$ is a functor $\C \to \C'$ cartesian over $\Posf$.
\end{defn}

\begin{prop}\label{pl.prop}
Let $i:\Posf \to \Posf^+$ be the embedding. Then for any restricted
enhanced category $\C \to \Posf$, there exists a small enhanced
category $\C^+$ and an equivalence $\C \cong \C^+$. Moreover, for
any two small enhanced categories $\C^+_0$, $\C^+_1$ and a
restricted enhanced functor $\gamma:i^*\C^+_0 \to \C^+_1$, there
exists an enhanced functor $\gamma^+;\C^+_0 \to \C^+_1$ and an
isomorphism $i^*(\gamma^+) \cong \gamma$, and $\gamma^+$ is unique
up to an isomorphism.
\end{prop}

\proof{} This is \cite[Corollary 7.3.1.5]{big}. \endproof

We also have a version of Proposition~\ref{pl.prop} for enhanced
groupoids, where we ask that $\C$ is a family of groupoids, and drop
the conditions that are automatic: a small reflexive family of
groupoids $\C \to \Posf$ that is additive, semiexact and satisfies
excision extends to an enhanced groupoid $\C^+ \to \Posf^+$, and any
restricted enhanced functor $\gamma:\C_0 \to \C_1$ extends to an
enhanced functor $\gamma^+:\C_0^+ \to \C_1^+$, uniquely up to an
isomorphism that need not be unique. This is \cite[Lemma
  7.2.1.10]{big} (combined with \cite[Corollary 7.3.1.5]{big}).

\subsection{Examples.}\label{exa.subs}

Trivial examples of enhanced categories are those of the form
$\Unf(\E)$, for a category $\E$. The first non-trivial example is
enhancement for the category $\Cat^0$ of small categories and
isomorphism classes of functors. Define a category $\cCat$ as
follows. Its objects are pairs $\langle J,\C \rangle$ of $J \in
\Posf^+$ and a small category $\C$ equipped with a fibration $\C \to
J$. Morphisms from $\langle J',\C' \rangle$ to $\langle J,\C
\rangle$ are commutative squares
\begin{equation}\label{ccat.sq}
\begin{CD}
\C' @>{\phi}>> \C\\
@VVV @VVV\\
J' @>{f}>> J,
\end{CD}
\end{equation}
where $f$ is a map in $\Posf^+$, and $\phi$ is cartesian over
$f$. Squares \eqref{ccat.sq} are considered modulo an
isomorphism. More generally, assume given a category $\E$, and
define a category $\cCat \bb \E$ as follows. Object are triples
$\langle J,\C,\alpha \rangle$ of some $\langle J,\C \rangle \in
\cCat$ and a $\alpha:\C_\perp \to \E$, where $\C_\perp \to J^o$ is
the transpose cofibration to the fibration $\C \to J$. Morphisms
from $\langle J',\C',\alpha' \rangle$ to $\langle J,\C,\alpha
\rangle$ are triples $\langle f,\phi,a \rangle$ of a square
\eqref{ccat.sq} defining a morphism in $\cCat$, and a map $a:\alpha'
\to \alpha \circ \phi_\perp$; these triples are considered modulo
isomorphisms $b:\phi' \cong \phi$ such that $a = a' \circ
\alpha(b_\perp)$. We have the forgetful functors
\begin{equation}\label{ccat.bb}
\cCat \bb \E \to \cCat, \qquad \langle J,\C,\alpha \rangle \mapsto
\langle J,\C \rangle
\end{equation}
and
\begin{equation}\label{ccat.pos}
\cCat \to \Posf^+, \qquad \langle J,\C \rangle \mapsto J.
\end{equation}
We note that while \eqref{ccat.pos} and the composition $\cCat \bb
\E \to \Posf^+$ of \eqref{ccat.bb} and \eqref{ccat.pos} are
fibrations, \eqref{ccat.bb} by itself is not (already over $\ppt \in
\Posf^+$). We have $\cCat_\ppt \cong \Cat^0$, $(\cCat \bb \E)_\ppt
\cong (\Cat \bb \E)^0$, and we can also define the extended versions
$\cCat^\Tot,\cCat^\Tot \bb \E \to \Pos$ of the families $\cCat,\cCat
\bb \E \to \Posf^+$ by allowing $J$ to be an arbitrary partially
ordered set.

Moreover, if the category $\E$ is essentially small, we can also
construct an enhancement $\E \bbd \cCat$ for the category $\E \bbd
\Cat$. Its objects are triples $\langle J,\C,\alpha \rangle$ of
$\langle J,\C \rangle \in \Cat$ and a functor $\alpha:\E \times J^o
\to \C_\perp$ over $J^o$, and morphisms from $\langle J',\C',\alpha'
\rangle$ to $\langle J,\C,\alpha \rangle$ are triples $\langle
f,\phi,a \rangle$ of a square \eqref{ccat.sq} defining a morphism in
$\cCat$, and a map $a:\alpha \circ \phi_\perp \circ (\id \times f^o)
\to \alpha'$; these triples are considered modulo isomorphisms
$b:\phi' \cong \phi$ such that $a = \alpha(b_\perp) \circ a'$. We
again have the forgetcul functor
\begin{equation}\label{ccat.bbd}
\E \bbd \cCat \to \cCat, \qquad \langle J,\C,\alpha \rangle \mapsto
\langle J,\C \rangle
\end{equation}
whose composition with \eqref{ccat.pos} is a fibration. We denote
$\cCat_\idot = \ppt \bbd \cCat$, so that $(\cCat_\idot)_\ppt \cong
\Cat^0_\idot$ is the category of \eqref{cat.loc}.

\begin{prop}\label{ccat.prop}
The family $\cCat \to \Posf^+$ is an enhanced category, and so is
the family $\cCat \bb \E \to \Posf^+$ for any category $\E$, while
$\cCat^\Tot$ resp.\ $\cCat^\Tot \bb \E$ are their canonical
bar-invariant extensions provided by Proposition~\ref{lf.prop}. If
the category $\E$ is essentially small, then the family $\E \bbd
\cCat \to \Posf^+$ is also an enhanced category.
\end{prop}

\proof{} As explained in \cite[Subsection 7.3.6]{big}, the statement
for $\cCat$ is actually a part of \cite[Proposition 7.3.6.1]{big}
complemented by \cite[Lemma 7.3.6.2]{big}, while the statements for
the lax functor categories are then parts of \cite[Propositions
  7.3.2.1,7.3.7.4]{big}, \endproof

The enhancement $\cCat \bb \E$ for $(\Cat \bb \E)^0$ given by
Proposition~\ref{ccat.prop} induces enhancements for the other two
categories of \eqref{3.dia} by cartesian squares
\begin{equation}\label{3.enh.dia}
\begin{CD}
\cCat \mm \E @>{a}>> \cCat \bbi \E @>{b}>> \cCat \bb \E\\
@VVV @VVV @VVV\\
\Unf((\Cat \mm \E)^0) @>{\Unf(a)}>> \Unf((\Cat \bbi \E)^0) @>{\Unf(b)}>>
\Unf((\Cat \bb \E)^0),
\end{CD}
\end{equation}
where the vertical arrow on the right is the truncation functor
\eqref{trunc.eq}, and the categories on the left only make sense
when $\E$ is essentially small. Explicitly, $\cCat \bbi \E \subset
\cCat \bb \E$ is the subcategory of triples $\langle J,\C,\alpha
\rangle$ such that $\alpha:\C_\perp \to \E$ is cocartesian over the
projection $J^o \to \ppt$, and morphisms are given by triples
$\langle f,\phi,a \rangle$ with invertible $a$. Note that in such a
situation, $\alpha \cong \beta_\perp$ for a unique functor $\beta:\C
\to \E$ cartesian over $J \to \ppt$, so that $\cCat \bbi \E$ is also the
category of triples $\langle J,\C,\beta \rangle$, $\beta:\C \to \E$
cartesian over $J \to \ppt$, with morphisms given by commutative
squares
\begin{equation}\label{ccat.1.sq}
\begin{CD}
\C' @>{\phi}>> \C\\
@V{\pi' \times \beta'}VV @VV{\pi \times \beta}V\\
J' \times \E @>{f \times \id}>> J \times \E
\end{CD}
\end{equation}
such that $\phi$ is cartesian over $f \times \id$, considered up to
an isomorphism over $J \times \E$. In these terms, $\cCat \mm \E
\subset \cCat \bbi \E$ is the subcategory of triples $\langle
J,\C,\beta \rangle$ such that $\pi \times \beta:\C \to J \times \E$
is a fibration, with morphisms given by squares \eqref{ccat.1.sq}
such that $\phi$ is cartesian over $f \times \id$ (for a proof, see
\cite[Lemma 7.4.4.1]{big}).

Alternatively, one can also construct $\cCat$ by considering the
fibration $\pi:\Unf(\Cat) \to \Posf^+$, and localizing $\Unf(\Cat)$
with respect to fiberwise equivalences --- that is, maps $f$ such
that $\pi(f) = \id_J$ for some $J \in \Posf^+$, and $f$ is a
pointwise equivalence of functors $J^o \to \Cat$. The same procedure
applied to the categories \eqref{3.dia} gives the enhanced
categories of \eqref{3.enh.dia}, and for $\E \bbd \Cat$, in
particular $\Cat_\idot$, we obtain $\E \bbd \cCat$, in particular
$\cCat_\idot$.

For another example of an enhanced category obtained by
localization, assume given a model category $\C$ (that is, a closed
model category in the original sense of \cite{qui.ho}). Then $\C$ is
localizable with respect to the class $W$ of weak equivalences, and
we can construct a natural enhancement for the localization
$h^W(\C)$. By definition, $\C$ is finitely complete, and for any
left-finite $J \in \Posf^\pm$, the functor category $J^o\C$ carries
a Reedy model structure. If $\C$ is complete, then the same holds
for any left-bounded $J \in \Posf^+$. Consider the fibration
$\pi:\Unf(\C) \to \Posf^+$, and let $\Unf(W)$ be the class of maps
$f$ in $\Unf(\C)$ such that $\pi(f) = \id_J$ for some $J \in \Pos$,
and $f$ is a pointwise weak equivalence in $\Unf(\C)_J \cong J^o\C$.

\begin{prop}\label{mod.prop}
The category $\Unf(\Posf^\pm,\C)$ is localizable with respect to the
class $\Unf(W)$, and the localization
$h^{\Unf(W)}(\Unf(\Posf^\pm,\C)) \to \Posf^\pm$ satisfies the
assumptions of Proposition~\ref{lf.prop}, thus defines an enhanced
category $\Hh^W(\C)$ such that $\Hh^W(\C)_\ppt \cong h^W(\C)$. If
$\C$ is complete, then $\Unf(\C)$ is also localizable with respect
to $\Unf(W)$, and $\Hh^W(\C) \cong h^W(\Unf(\C))$. A Quillen-adjoint
pair of functors $\lambda:\C \to \C'$, $\rho:\C' \to \C$ between
model categories $\C$, $\C'$ induces an adjoint pair of enhanced
functors $\Hh(\lambda):\Hh^W(\C) \to \Hh^W(\C')$,
$\Hh(\rho):\Hh^W(\C') \to \Hh^W(\C)$.
\end{prop}

\proof{} This is \cite[Lemma 7.2.4.1]{big}. \endproof

For yet another example, consider the category $C_\idot(\A)$ of
chain complexes in an additive category $\A$, and let $\Ho(\A)$ be
the corresponding homotopy category (objects are complexes,
morphisms are chain-homotopy equivalence classes of
maps). Alternatively, $\Ho(\A)$ is obtained by localizing
$C_\idot(\A)$ with respect to the class $W$ of chain-homotopy
equivalences. To construct an enhancement for $\Ho(\A)$, note that
for any $J \in \Posf^\pm$, the category $J^o\A$ is also additive,
and we have $C_\idot(J^o\A) \cong J^oC_\idot(\A)$. We also have the
homotopy category $\Ho(J^o\A)$, but this turns out to be the wrong
thing to consider. To get it right, one has to either impose an
additional condition on complexes in $J^o\A$, or -- which is simpler
-- to enlarge the class $W$. Thus we consider the fibration
$\pi:\Unf(\Posf^\pm,C_\idot(\A)) \to \Posf^\pm$, and we let $W$ be a
class of maps $f$ in $\Unf(\Posf^\pm,C_\idot(\A))$ such that $\pi(f)
= \id_J$ for some $J \in \Posf^\pm$, and $f$ is a pointwise
chain-homotopy equivalence in $J^oC_\idot(\A)$ (that is, becomes a
chain-homotopy equivalence after evaluation at any $j \in J$).

\begin{prop}\label{ho.prop}
The category $\Unf(\Posf^\pm,C_\idot(\A))$ is localizable with
respect to the class $W$, and the localization
$h^W(\Unf(\Posf^\pm,C_\idot(\A)))$ satisfies the assumptions of
Proposition~\ref{lf.prop}, thus defines an enhanced category
$\Hho(\A)$ such that $\Hho(\A)_\ppt \cong \Ho(\A)$.
\end{prop}

\proof{} This is \cite[Lemma 7.2.4.2]{big}. \endproof

Finally, if we restrict our attention to small enhanced categories,
then a huge supply of those is provided by complete Segal spaces of
Subsection~\ref{del.subs}. Namely, recall that any complete Segal
space $X \in \Delta^o\Delta^o\Sets$ defines a small category $h(X) =
h^{+ \times \Tot}(\Delta^2X)$, where $\Delta^2X \to \Delta^2$ is the
corresponding discrete fibration. If $X = N^2(I)$ for small category
$I$, then $h(X) \cong I$. Now consider the embedding
\begin{equation}\label{rho.eq}
\rho = ((\iota \circ \phi) \times \id \times \id) \circ (\delta
\times \id):\Delta \times \Delta \to \Delta \times \Delta \times
\Delta \to \Posf^+ \times \Delta \times \Delta,
\end{equation}
where $\iota:\Delta \to \Delta$ is the involution $[n] \mapsto
[n]^o$, $\delta:\Delta \to \Delta \times \Delta$ is the diagonal
embedding, and $\phi:\Delta \to \Posf^+$ is the standard
embedding. Then for any $X \in \Delta^o\Delta^o\Sets$, let
\begin{equation}\label{unf.X}
\Unf(X) = h^{\id \times + \times \Tot}((\Posf^+ \times \Delta \times
\Delta)\rho^o_*X),
\end{equation}
where $\pi:(\Posf^+ \times \Delta \times \Delta)\rho^o_*X \to
\Posf^+ \times \Delta \times \Delta$ is the discrete fibration
corresponding to the right Kan extension $\rho^o_*X:(\Posf^+ \times
\Delta \times \Delta)^o \to \Sets$, and $\id \times + \times \Tot$
is the class of maps $f$ such that $\pi(f) = f_0 \times f_1 \times
f_2$, with $f_0 = \id_J$ for some $J \in \Posf^+$, special $f_1:[n]
\to [m]$ (that is, $f_1(n)=m$), and any $f_2$. By definition,
$\Unf(X)$ comes with a fibration $\Unf(X) \to \Posf^+$, and if one
computes $\phi^o_*X$ by \eqref{kan.eq}, then for any $J \in
\Posf^+$, one obtains an equivalence
\begin{equation}\label{unf.X.J}
\Unf(X)_J \cong h(X(J^o)^{(2)}),
\end{equation}
where $X(J^o)^{(2)}$ is an in \eqref{X.n.2}. In particular, we have
$\Unf(X)_\ppt \cong h(X)$. If $X = N^2(I)$ is the double nerve of a
small category $I$, then $\Unf(X) \cong \Unf(I)$.

\begin{prop}\label{seg.prop}
For any complete Segal space $X$, the category $\Unf(X)$ of
\eqref{unf.X} with its fibration $\Unf(X) \to \Posf^+$ is a small
enhanced category, and for any map $f:X \to X'$ of complete Segal
spaces, $\Unf(f):\Unf(X) \to \Unf(X')$ is an enhanced functor. If
$X:\Delta^o \times \Delta^o \to \Sets$ is constant along the first
factor $\Delta^o$, then $\Unf(X)$ is an enhanced groupoid.
\end{prop}

\proof{} This is \cite[Lemma 7.3.1.1]{big}. \endproof

\subsection{Representability.}

As it happens, the last example of Subsection~\ref{exa.subs}
provided by Proposition~\ref{seg.prop} is universal. Namely, we have
the following representability theorem.

\begin{theorem}\label{repr.thm}
For any small enhanced category $\C$, there exists a complete Segal
space $X$ and an equivalence $\C \cong \Unf(X)$. For any complete
Segal spaces $X$, $X'$ and an enhanced functor $\gamma:\Unf(X) \to
\Unf(X')$, there exists a map $f:X \to X'$ and an isomorphism
$\gamma \cong \Unf(f)$, and two maps $f,f':X \to X'$ define functors
$\Unf(f)$, $\Unf(f')$ isomorphic over $\Posf^+$ if and only if they
are homotopic --- that is, define the same map in
$h^W(\Delta^o\Delta^o\Sets)$. If $\C$ is $\kappa$-bounded, for a
regular cardinal $\kappa$, then one can choose $X \in
\Delta^o\Delta^o\Sets_\kappa$. If $\C$ is an enhanced groupoid, one
can choose $X:\Delta^o \times \Delta^o \to \Sets$ constant along the
first factor $\Delta^o$.
\end{theorem}

\proof{} This is \cite[Theorem 7.3.1.3]{big}. \endproof

As an immediate corollary of Theorem~\ref{repr.thm}, we see that for
any two small enhanced categories $\C$, $\C'$, isomorphism classes
of enhanced functors $\C \to \C'$ form a set (\cite[Corollary
  7.3.1.4]{big}), so that we have a well-defined category of small
enhanced categories and isomorphism classes of enhanced functors
that we denote by $\Cat^h$. Theorem~\ref{repr.thm} then provides an
equivalence
\begin{equation}\label{cat.h.eq}
\Cat^h \cong h^W_{css}(\Delta^o\Delta^o\Sets),
\end{equation}
where $h^W_{css}(\Delta^o\Delta^o\Sets) \subset
h^W(\Delta^o\Delta^o\Sets)$ is the full subcategory spanned by
complete Segal spaces of Definition~\ref{css.def}. Moreover,
\eqref{cat.h.eq} restricts to an equivalence
\begin{equation}\label{sets.h.eq}
  \Sets^h \cong h^W(\Delta^o\Sets),
\end{equation}
where $\Sets^h \subset \Cat^h$ is the full subcategory spanned by
small enhanced groupoids -- in enhanced context, an enhanced
groupoid is really an enhanced version of a set -- and the full
subcategory $h^W(\Delta^o\Sets) \subset
h^W_{css}(\Delta^o\Delta^o\Sets)$ is spanned by $X:\Delta^o \times
\Delta^o \to \Sets$ constant along the first factor $\Delta^o$. By
Proposition~\ref{pl.prop}, the same category $\Cat^h$ can be
described as the category of restricted small enhanced categories of
Definition~\ref{restr.enh.def} and isomorphism classes of enhanced
functors. However, there is more: as less direct corollary of
Theorem~\ref{repr.thm}, we also have the following result.

\begin{defn}\label{enh.cocart.def}
A commutative square \eqref{cat.sq} of enhanced categories and
enhnaced functors is {\em enhanced-cocartesian} if for any enhanced
category $\C'$, enhanced functors $\gamma_l:\C_l \to \C'$, $l =
0,1$, and isomorphism $\gamma'_0 \circ \gamma^0_{01} \cong
\gamma'_1 \circ \gamma^1_{01}$, there exists an enhanced functor
$\gamma:\C \to \C'$ and isomorphisms $a_l:\gamma \circ \gamma_l
\cong \gamma_l'$, $l=0,1$, and for any enhanced functors
$\gamma,\gamma':\C \to \C'$ and enhanced maps $a_l:\gamma \circ
\gamma_l \to \gamma' \circ \gamma_l$, $l=0,1$ such that
$\gamma^{0*}_{01}(a_0) = \gamma^{1*}_{01}(a_1)$, there exists a
unique enhanced map $a:\gamma \to \gamma'$ such that $\gamma_l^*(a)
= a_l$, $l=0,1$.
\end{defn}

\begin{prop}\label{enh.loc.prop}
For any small enhanced category $\C$, there exists a partially
ordered set $J \in \Posf^+$, a set $W$ of maps $w:[1] \to J$,
and an enhanced functor $\phi:\Unf(J^o) \to \C$ that fits into a
commutative square
\begin{equation}\label{pre.sq}
\begin{CD}
W \times \Unf([1]^o) @>{q}>> \Unf(J^o)\\
@V{p}VV @VV{\phi}V\\
W \times \ppt^h @>{i}>> \C
\end{CD}
\end{equation}
of enhanced categories and enhanced functors that is
enhanced-cocartesian in the sense of
Definition~\ref{enh.cocart.def}.
\end{prop}

\proof{} Immediately follows from the more precise \cite[Proposition
  7.3.3.4]{big}, combined with \cite[Example
  7.3.4.5]{big}. \endproof

The square \eqref{pre.sq} is an enhanced version of the square
\eqref{loc.sq} (in particular, $W \times \Unf([1]^o)$ resp.\ $W
\times \ppt^h$ are coproduct of copies of $\Unf([1]^o)$
resp.\ $\ppt^h \cong \Posf^+$ numbered by elements $w \in W$, $q$ is
the coproduct of the embeddings $\Unf(w^o):\Unf([1]^o) \to
\Unf(J^o)$, and $p$ is the coproduct of projections $\Unf([1]^o) \to
\ppt^h$). Informally, Proposition~\ref{enh.loc.prop} says that every
small enhanced category can be obtained as an enhanced localization
of a partially ordered set. In particular, this immediately implies
that for any small enhanced category $\C$ and any enhanced category
$\E$, enhanced functors $\C \to \E$ and enhanced morphisms between
them form a well-defined category $\Fun^h(\C,\E)$ -- effectively,
for any $J \in \Posf^+$, we have $\Fun^h(\Unf(J^o),\E) \cong \E_J$
by Lemma~\ref{J.C.le}, and then a cocartesian square \eqref{pre.sq}
for $\C$ induces a cartesian square
\begin{equation}\label{fun.h.sq}
\begin{CD}
\Fun^h(\C,\E) @>>> \E_J\\
@VVV @VVV\\
\prod_{w \in W}\E_\ppt @>>> \prod_{w \in W}\E_{[1]}.
\end{CD}
\end{equation}
In particular, an enhanced map $\gamma \to \gamma'$ between two
enhanced functors $\gamma,\gamma':\C \to \E$ is always induced by
the universal map \eqref{ar.map} via an enhanced functor $\C \to
\Ar^h(\E)$. As another immediate corollary, for any enhanced
category $\E$, we have a well-defined category $\Cat^h \bb^h \E$
whose objects are small enhanced categories $\C$ equipped with
enhanced functors $\alpha:\C \to \E$, with morphisms from $\langle
\C',\alpha' \rangle$ to $\langle \C,\alpha \rangle$ defined by pairs
$\langle \phi,a \rangle$ of an enhanced functor $\phi:\C' \to \C$
and an enhanced map $a:\alpha' \to \alpha \circ \phi$, considered up
to an enhanced isomorphism $b:\phi' \cong \phi$ such that $a = a'
\circ \alpha(b)$. We also have the subcategory $\Cat^h \bbi^h \E
\subset \Cat^h \bb^h \E$ with the same objects, and morphisms given
by pairs $\langle \phi,a \rangle$ with invertible $a$, and if $\E$
is small, we have the category $\E \bbd^h \Cat^h$ of small enhanced
categories $\C$ equipped with enhanced functors $\alpha:\E \to \C$,
with morphisms from $\langle \C',\alpha' \rangle$ to $\langle
\C,\alpha \rangle$ defined by pairs $\langle \phi,a \rangle$ of an
enhanced functor $\phi:\C' \to \C$ and an enhanced map $a:\alpha'
\to \alpha \circ \phi$, considered up to an enhanced isomorphism
$b:\phi' \cong \phi$ such that $a = \alpha(b) \circ a'$. As in
\eqref{cat.loc}, we denote $\Cat^h_\idot = \ppt^h \bbd^h
\Cat^h$. There are two other less immediate but very useful
corollaries.

\begin{corr}\label{fun.h.corr}
For any small enhanced category $\C$ and enhanced category $\E$,
there exists an enhanced category $\fFun^h(\C,\E)$ and en enhanced
functor
\begin{equation}\label{ev.C.E}
  \ev:\C \times^h \fFun^h(\C,\E) \to \E,
\end{equation}
where $- \times^h -$ is the enhanced product \eqref{enh.prod.eq},
such that for any enhanced category $\C'$, any enhanced functor
$\gamma:\C \times^h \C' \to \E$ factors as
\begin{equation}\label{facto.C.E}
\begin{CD}
\C \times^h \C' @>{\id \times \wt{\gamma}}>> \C \times^h
\fFun^h(\C,\E) @>{\ev}>> \E
\end{CD}
\end{equation}
for an enhanced functor $\wt{\gamma}:\C' \to \fFun^h(\C,\E)$, and
$\wt{\gamma}$ and the factorization \eqref{facto.C.E} are unique up
to a unique enhanced isomorphism.
\end{corr}

\proof{} This is \cite[Corollary 7.3.3.5]{big}. \endproof

In particular, the universal property of Corollary~\ref{fun.h.corr}
for $\C' = \ppt^h$ provides an identification between the fiber
$\fFun^h(\C,\E)_\ppt$ and the category $\Fun^h(\C,\E)$ of
\eqref{fun.h.sq}, so that Corollary~\ref{fun.h.corr} provides a
canonical enhancement for the functor category
$\Fun^h(\C,\E)$. Explicitly, objects in $\fFun^h(\C,\E)$ are pairs
$\langle J,\gamma \rangle$, $J \in \Posf^+$, $\gamma:\C \times^h
\Unf(J^o) \to \C$ an enhanced functor, with morphisms from $\langle
J',\gamma' \rangle$ to $\langle J,\gamma \rangle$ given by pairs
$\langle f,\alpha \rangle$ of a map $f:J' \to J$ and an enhanced map
$\alpha:\gamma' \to \gamma \circ (\id \times^h \Unf(f^o))$. If $\E$
itself is small, then Corollary~\ref{fun.h.corr} implies that
$\Cat^h$ is cartesian-closed. As in the unenhanced case, we simplify
notation by writing $\C^\iota\E = \fFun^h(\C^\iota,\E)$ for any
enhanced category $\E$ and small enhanced category $\C$, and we
further simplify to $I^o_h\E = \Unf(I)^\iota\E$ when $\C = \Unf(I)$
for some essentially small category $I$.

\begin{corr}\label{semi.corr}
\begin{enumerate}
\item For any enhanced categories $\C$, $\C_1$, $\C_0$ and enhanced
  functors $\gamma_l:\C_l \to \C$, $l=0,1$ such that $\C_0$ and
  $\gamma^1$ is small, there exists a semicartesian commutative
  square \eqref{cat.sq} of enhanced categories and enhanced
  functors, with small $\gamma^0_{01}$ and $\C_{01}$.
\item Assume given a semicartesian commutative squares
  \eqref{cat.sq} of enhanced categories and enhanced functors, and
  another semicartesian commutative square
\begin{equation}\label{cat.1.sq}
\begin{CD}
\C'_{01} @>{{\gamma'}_{01}^1}>> \C_1\\
@V{{\gamma'}_{01}^0}VV @VV{\gamma_1}V\\
\C_0 @>{\gamma_0}>> \C
\end{CD}
\end{equation}
of enhanced categories and enhanced functors such that $\C'_{01} \to
\Posf^+$ is small. Then there exists an enhanced functor
$\phi:\C'_{01} \to \C_{01}$ and enhanced isomorphisms
$a_l:\gamma^l_{01} \circ \phi \cong {\gamma'}^l_{01}$, $l = 0,1$,
and the triple $\langle \phi,a_0,a_1 \rangle$ is unique up to an
enhanced isomorphism.
\end{enumerate}
\end{corr}

\proof{} This is \cite[Lemma 7.3.37]{big} and \cite[Corollary
  7.3.3.6]{big}. \endproof

If $\C$ and $\C_0$ in Corollary~\ref{semi.corr} are of the form
$\Unf(\E)$, $\Unf(\E')$ for some categories $\E$, $\E'$, then the
semicartesian square is the cartesian square \eqref{U.prod}. The
square is actually cartesian for any $\C_0$, as long as $\C \cong
\Unf(\E)$, or either of $\gamma_l$, $l=0,1$, is fully faithful. But
in general, taking a cartesian square does not work: $\C_0 \times_\C
\C_1$ is not an enhanced category (what breaks is
semiexactness). What Corollary~\ref{semi.corr} shows is that for
small enhanced categories this can be corrected, and in a somewhat
functorial way. In particular, by
Corollary~\ref{semi.corr}~\thetag{ii}, the semicartesian square
provided by Corollary~\ref{semi.corr}~\thetag{i} is unique up to a
unique equivalence. We call the corresponding category $\C_{01}$ the
{\em semicartesian product} of $\C_0$ and $\C_1$ over $\C$, and
denote it by $\C_0 \times^h_\C \C_1$ (that reduces to
\eqref{enh.prod.eq} for $\C = \ppt^h$). One has to remember that
while the semicartesian product is unique up a unique equivalence,
the equivalence itself is only unique up to a {\em non-unique}
isomorphism. This is not a bug but a feature; this is precisely
where the whole homotopy theory comes from.

\section{Enhanced category theory.}\label{enh.cat.sec}

\subsection{Categories of categories.}

Let us now describe the basics of enhanced category theory, in
parallel to basic category theory of Section~\ref{cat.sec}. We start
with constructing an enhancement for the category \eqref{cat.h.eq}
of small enhanced categories, and the categories $\Cat^h \bb^h \E$
for enhanced categories $\E$.

\begin{defn}
For any enhanced category $\C$ and partially ordered set $I$, an
{\em $I$-augmentation} for $\C$ is an enhanced functor $\pi:\C \to
\Unf(I)$ that is a fibration, and an {\em $I$-coaugmentation} is an
enhanced functor $\pi:\C \to \Unf(I)$ such that $\pi^\iota:\C^\iota
\to \Unf(I)^\iota \cong \Unf(I^o)$ is an $I^o$-augmentation. For any
map $f:I' \to I$ and commutative square
\begin{equation}\label{ccat.h.sq}
\begin{CD}
\C' @>{\phi}>> \C\\
@V{\pi'}VV @VV{\pi}V\\
\Unf(J') @>{\Unf(f)}>> \Unf(J)
\end{CD}
\end{equation}
of enhanced categories and enhanced functors such that $\pi$, $\pi'$
are augmentations, $\phi$ is {\em augmented over $f$} if it it
cartesian over $\Unf(f)$, and if $\pi$, $\pi'$ are coaugmentations,
then $\pi$ is {\em coaugmented over $f$} if $\pi^\iota$ is augmented
over $f^o$. In particular, for any $I$-augmentations
resp.\ $I$-coaugmentations $\C,\C' \to \Unf(I)$, an enhanced functor
$\gamma:\C' \to \C$ over $I$ is {\em augmented} resp.\ {\em
  coaugmented} if it is augmented resp.\ coaugmented over $\id_I$.
\end{defn}

For any $I \in \Posf^+$, $I$-augmented small enhanced categories and
isomorphism classes of $I$-augmented enhanced functors between them
form a category that we denote by $\Cat^h(I)$. If $\C = \Unf(\E)$,
$\pi=\Unf(\gamma)$ for some $\gamma:\E \to I$, then $\pi$ is an
$I$-augmentation resp.\ $I$-coaugmentation iff $\gamma$ is a
fibration resp.\ cofibration (\cite[Lemma 7.4.3.1]{big}), and in the
former case, $\pi^\iota \cong \Unf(\gamma^o)$ is the coaugmentation
corresponding to the opposite cofibration $\E^o \to I^o$. We also
have an enhanced version of the transpose cofibration $\E_\perp \to
I^o$ and the transpose-opposite fibration $\E^o_\perp \to
I$. Namely, if a partially ordered set $J$ is equipped with a map
$\alpha:J \to I^o$, then $J^o$ of course does not admit a natural
map to $I^o$, but both $B(J^o) \cong B(J)$ and $\overline{J^o} \cong
\overline{J}$ do --- we have maps $\alpha \circ \xi:B(J) \to I^o$
and $\alpha \circ i:\overline{J} \to I^o$, where $\xi$ and $i$ are
as in \eqref{BJ.sq}. Therefore, we can extend the functors $B \circ
\iota$ and $S \circ \iota$ in \eqref{C.iota} to functors
$B_I^\iota,S_I^\iota:\Unf(I) \to \Unf(I)$ by setting
$B_I^\iota(\langle J,\alpha \rangle) = \langle B(J),\alpha \circ \xi
\rangle$ and $S_I^\iota(\langle J,\alpha \rangle) = \langle
\overline{J},\alpha \circ i \rangle$. Moreover, if we consider the
square \eqref{BJ.sq} for $J^o$, with the corresponding maps
$i':\oB(J^o) \to B(J^o) \cong B(J)$, $\xi':\oB(J^o) \to
\overline{J^o} \cong \overline{J}$, then $\alpha \circ \xi \circ i'
\geq \alpha \circ i \circ \xi'$ pointwise. Thus if we define
$\oB^\iota_I(\langle J,\alpha \rangle)=\langle \oB(J^o),\alpha \circ
i \circ \xi'\rangle$, then $\xi'$ is a map over $I^o$, and $i'$ is
still a map in $\Unf(I)$. Therefore, for any $I$-augmented enhanced
category $\C \to \Unf(I)$, we can define a category
$\C^\iota_{h\perp}$ fibered over $\Unf(I)$ by the cartesian square
\begin{equation}\label{C.I.iota}
\begin{CD}
\C^\iota_{h\perp} @>>> (S_I^\iota)^*\C^o_\perp\\
@VVV @VV{{\xi'}^*}V\\
(B_I^\iota)^*\C^o_\perp @>{{i'}^*}>> (\oB_I^\iota)^*\C^o_\perp.
\end{CD}
\end{equation}
One checks (\cite[Subsection 7.2.2]{big}) that $\C^\iota_{h\perp}$
is an enhanced category, it is $I$-augmented by the fibration
$\C^\iota_{h\perp} \to \Unf(I)$, and gives the {\em opposite
  coaugmentation} $\C_{h\perp} \to \Unf(I^o)$.

Now as in Subsection~\ref{exa.subs}, define a category $\cCat^h$ as
follows. Objects are pairs of a partially ordered set $I \in \Posf^+$
and an $I$-augmented small enhanced category $\C$. Morphisms from
$\langle I',\C' \rangle$ to $\langle I,\C \rangle$ are defined by
pairs $\langle f,\phi \rangle$ of a map $f:I' \to I$ and a
commutative square \eqref{ccat.h.sq} such that $\phi$ is an enhanced
functor augmented over $f$.  Two pairs $\langle f,\phi\rangle$,
$\langle f',\phi'\rangle$ define the same morphism if $f=f'$, and
$\phi \cong \phi'$ over $\Posf^+$. Moreover, for any enhanced
category $\E$, define a category $\cCat^h \bb^h \E$ as
follows. Objects are triples $\langle I,\C,\alpha \rangle$, where
$\langle I,\C \rangle$ gives an object in $\cCat^h$, and
$\alpha:\C_{h\perp} \to \E$ is an enhanced functor. Morphisms
$\langle I',\C',\alpha' \rangle \to \langle I,\C,\alpha \rangle$ are
defined by triples $\langle f,\phi,a \rangle$ of a commutative
square \eqref{ccat.h.sq} defining a morphism in $\Cat^h$, and an
enhanced map $a:\alpha' \to \alpha \circ \phi_{h\perp}$. Two triples
$\langle f,\phi,a \rangle$, $\langle f',\phi',a' \rangle$ define the
same morphism if $f=f'$, and there exists an isomorphism $b:\phi'
\cong \phi$ over $\Unf(I)$ such that $a = a' \circ
\alpha(b_{h\perp})$. As in Subsection~\ref{exa.subs}, we have
a forgetful functor
\begin{equation}\label{ccat.h.bb}
\cCat^h \bb^h \E \to \cCat^h, \qquad \langle I,\C,\alpha \rangle
\mapsto \langle I,\C \rangle,
\end{equation}
an enhanced version of \eqref{ccat.bb}, and a forgetful functor
\begin{equation}\label{ccat.h.pos}
\cCat^h \to \Posf^+, \qquad \langle I,\C \rangle \mapsto I,
\end{equation}
an enhanced version of \eqref{ccat.pos}. Just as in
Subsection~\ref{exa.subs}, both \eqref{ccat.h.pos} and the
composition of \eqref{ccat.h.pos} and \eqref{ccat.h.bb} are
fibrations, with fibers $\cCat^h_\ppt \cong \Cat^h$, $\cCat^h_I
\cong \Cat^h(I)$, $(\cCat^h \bb^h \E)_\ppt \cong \Cat^h \bb^h \E$,
but \eqref{ccat.h.bb} by itself is not a fibration. We then
construct fibrations $\cCat^{h\Tot},\cCat^{h\Tot} \bb^h \E \to \Pos$
by allowing $I$ to be an arbitrary partially ordered set.

Moreover, if $\E$ is small, define a category $\E \bbd^h \cCat^h$ as
follows. Objects are triples $\langle I,\C,\alpha \rangle$ of $I \in
\Posf^+$, an object $\langle I,\C \rangle \in \cCat^h$, and an
enhanced functor $\alpha:\E \times^h \Unf(I^o) \to \C_{h\perp}$ over
$\Unf(I^o)$. Morphisms $\langle I',\C',\alpha' \rangle \to \langle
I,\C,\alpha \rangle$ are defined by triples $\langle f,\phi,a
\rangle$ of a commutative square \eqref{ccat.h.sq} defining a
morphism in $\Cat^h$, and an enhanced map $a:\alpha' \to \alpha
\circ \phi_{h\perp}$. Two triples $\langle f,\phi,a \rangle$,
$\langle f',\phi',a' \rangle$ define the same morphism if $f=f'$,
and there exists an isomorphism $b:\phi' \cong \phi$ over $\Unf(I)$
such that $a = a' \circ \alpha(b_{h\perp})$. We again have a
forgetful functor
\begin{equation}\label{ccat.h.bbd}
\E \bbd^h \cCat^h \to \cCat^h, \qquad \langle I,\C,\alpha \rangle
\mapsto \langle I,\C \rangle,
\end{equation}
an enhanced version of \eqref{ccat.bbd}, whose composition with
\eqref{ccat.h.pos} is a fibration, so that $\E \bbd^h \cCat^h$ is a
family over $\Posf^+$. We have $(\E \bbd^h \cCat^h)_\ppt \cong \E
\bbd^h \Cat^h$, and we denote $\cCat^h_\idot = \ppt^h \bbd^h
\cCat^h$.

\begin{prop}\label{ccat.h.prop}
The family $\cCat^h \to \Posf^+$ is an enhanced category, and so are
the family $\cCat^h \bb^h \E \to \Posf^+$ for any enhanced category
$\E$, and the family $\E \bbd^h \cCat^h \to \Posf^+$ for any small
enhanced category $\E$. The fibrations $\cCat^{h\Tot},\cCat^{h\Tot}
\bb^h \E \to \Pos$ are the canonical bar-invariant extensions of the
enhanced categories $\cCat^h$, $\cCat^h \bb^h \E$. Moreover, the
equivalence \eqref{cat.h.eq} extends to an enhanced equivalence
\begin{equation}\label{ccat.h.eq}
  \cCat^h \cong \Hh_{css}^W(\Delta^o\Delta^o\Sets),
\end{equation}
where $\Hh^W(\Delta^o\Delta^o\Sets)$ is the enhanced localization of
Proposition~\ref{mod.prop}, and $\Hh^W_{css}(\Delta^o\Delta^o\Sets)
\subset \Hh^W(\Delta^o\Delta^o\Sets)$ is the full enhanced
subcategory corresponding to $h^W_{css}(\Delta^o\Delta^o\Sets)
\subset h^W(\Delta^o\Delta^o\Sets)$.
\end{prop}

\proof{} The first claim is \cite[Propositions
  7.3.6.1,7.3.7.1,7.3.7.4]{big}, the second claim is
\cite[Lemma 7.3.6.2]{big}, and \eqref{ccat.h.eq} is
\cite[\thetag{7.3.6.2}]{big}. \endproof

We also note, see \cite[Subsection 7.3.7]{big}, that we have full
embeddings $\Pos \subset \Cat^0 \subset \Cat^h$, and these extend to
enhanced full embeddings $\pPos \subset \cCat \subset \cCat^h$,
where the enhancement $\pPos \cong \Unf(\Pos)$ for $\Pos$ is trivial
-- since partially ordered sets are rigid, there is nothing to
enhance. By the same rigidity, the projection \eqref{ccat.h.bb} is a
fibration over $\pPos \subset \cCat^h$. More precisely, if we define
a category $\pPos \bb^h \E$ by the cartesian square
\begin{equation}\label{ppos.E}
\begin{CD}
  \pPos \bb^h \E @>>> \cCat^h \bb^h \E\\
  @VVV @VVV\\
  \pPos @>>> \cCat^h,
\end{CD}
\end{equation}
then the vertical arrow on the left is a fibration, and $\pPos \bb^h
\E$ is an enhanced category. Explicitly, $\Pos \bb^h \E = (\pPos
\bb^h \E)_\ppt$ is naturally identified with $\iota^*\E^\Tot$, where
$\E^\Tot \to \Pos$ is the canonical bar-invariant extension of the
enhanced category $\E$, and the equivalence $\iota:\Pos \bb^h \E
\cong \iota^*\E^\Tot \to \E^\Tot$ extends to an equivalence $\pPos
\bb^h \E \cong \sigma^*\E$, where $\sigma^*\E$ is the enhanced
category of Lemma~\ref{C.refl.le}. Another useful full subcategory
in $\Cat^h$ is the subcategory $\Sets^h \subset \Cat^h$ of enhanced
groupoids; it inherits an enhancement $\sSets^h \subset \cCat^h$,
and \eqref{sets.h.eq} extends to an enhanced equivalence
\begin{equation}\label{ssets.h.eq}
\sSets^h \cong \Hh^W(\Delta^o\Sets).
\end{equation}
For any enhanced category $\E$, we then also have the full enhanced
subcategory $\sSets^h \bb^h \E \subset \cCat^h \bb^h \E$. One easily
checks that the intersection $\Pos \cap \Sets^h \subset \Cat^h$ is
the category $\Sets$ of usual discrete sets, so that the enhanced
category $\E \cong (\cCat^h \bb^h \E)_\ppt \subset \cCat^h \bb^h \E$
is contained both in $\pPos \bb^h \E$ and in $\sSets^h \bb^h \E$.

\subsection{Enhanced fibrations and cofibrations.}

Next, we construct enhanced versions of cylinders and
comma-categories of Subsection~\ref{cat.subs}. For any enhanced
functor $\gamma:\C_0 \to \C_1$ between small enhanced categories
$\C_0$, $\C_1$, the enhanced cylinder and dual cylinder are defined
by semicartesian squares
\begin{equation}\label{enh.cyl.sq}
\begin{CD}
\Cyl_h(\gamma) @>>> \C_0^{h>}\\
@VVV @VV{\gamma^{h>}}V\\
\Unf([1]) \times^h \C_1 @>>> \C_1^{h>},
\end{CD}
\qquad\qquad
\begin{CD}
\Cyl_h^\iota(\gamma) @>>> \C_0^{h<}\\
@VVV @VV{\gamma^{h<}}V\\
\Unf([1]) \times^h \C_1 @>>> \C_1^{h<},
\end{CD}
\end{equation}
where $\C_l^{h<}$, $\C_l^{h>}$, $l=0,1$ are the enhanced categories
\eqref{h.gt.lt}, and the bottom arrows are the enhanced functors
\eqref{h.gt.lt.1}. As in the unenhanced case, we have enhanced full
embeddings $s:\C_0 \to \Cyl_h(\gamma)$, $t:\C_1 \to \Cyl_h(\gamma)$,
the latter is left-reflexive with an adjoint enhanced functor
$t_\dg$, and $\gamma \cong t_\dg \circ s$. Dually, we have enhanced
full embeddings $s:\C_1 \to \Cyl^\iota_h(\gamma)$, $t:\C_0 \to
\Cyl_h(\gamma)$, the former is right-reflexive with an adjoint
enhanced functor $s^\dg$, and we have $\gamma \cong s^\dg \circ
t$. Altogether, $\gamma$ admits an enhanced version of the two
factorizations \eqref{cyl.facto}, while both the cylinder
$\Cyl_h(\gamma)$ and the dual cylinder $\Cyl^\iota_h(\gamma)$ are
enhanced categories under $\C_0 \copr \C_1$, with respect to the
enhanced functor $s \copr t$. Moreover, the dual cylinder
$\Cyl^\iota_h(\gamma)$ resp.\ cylinder $\Cyl_h(\gamma)$ is
$[1]$-augmented resp.\ $[1]$-coaugmented by the projection
$\Cyl_h(\gamma),\Cyl_h^\iota(\gamma) \to \Unf([1]) \times^h \C_1 \to
\Unf([1])$. Conversely, any $[1]$-augmented enhanced category $\C$,
with enhanced fibers $\C_0$, $\C_1$, is of the form $\C \cong
\Cyl^\iota_h(\gamma)$ for a functor $\gamma:\C_1 \to \C_0$, unique
up to an isomrphism, and dually, any $[1]$-coaugmented enhanced
category $\C$ is of the form $\C \cong \Cyl_h(\gamma)$.

\begin{lemma}\label{adj.cyl.le}
Assume given a pair $\lambda:\C_0 \to \C_1$, $\rho:\C_1 \to \C_0$ of
enhanced functors between small enhanced categories. Then pairs of
enhanced maps $\lambda \circ \rho \to \id$, $\id \to \rho \circ
\lambda$ defining an adjunction between $\lambda$ and $\rho$
correspond bijectively to isomorphism classes of equivalences
$\Cyl(\lambda) \cong \Cyl^o(\rho)$ under $\C_0 \copr \C_1$.
\end{lemma}

\proof{} This is part of \cite[Lemma 7.4.1.5]{big}. \endproof

\begin{corr}\label{adj.pb.corr}
If we have enhanced functors $\C_1,\C_1' \to \C$ between small
enhanced categories, and a pair of enhanced functors $\lambda:\C_1
\to \C_1'$, $\rho:\C_1' \to \C_1$ adjoint over $\C$, then for any
small enhanced $\C_0$ and any enhanced functor $\gamma:\C_0 \to \C$,
the enhanced functors $\gamma^*(\lambda):\C_0 \times^h_\C \C_1 \to
\C_0 \times^h \C_1'$, $\gamma^*(\rho):\C_0 \times^h_\C \C'_1 \to
\C_0 \times^h \C_1$ provided by Corollary~\ref{semi.corr} form an
adjoint pair.
\end{corr}

\proof{} This is \cite[Corollary 7.4.1.6]{big}. \endproof

For any enhanced functor $\pi:\C \to \E$ between small enhanced
categories, the enhanced left and right comma-categories are defined
by semicartesian squares
\begin{equation}\label{enh.comma.sq}
\begin{CD}
\C /^h_\pi \E @>>> \Ar^h(\E)\\
@V{\sigma}VV @VV{\sigma}V\\
\C @>{\pi}>> \E,
\end{CD}
\qquad\qquad
\begin{CD}
\E \setminus^h_\pi \C @>>> \Ar^h(\E)\\
@V{\tau}VV @VV{\tau}V\\
\C @>{\pi}>> \E,
\end{CD}
\end{equation}
where $\Ar^h(\E)$, $\sigma$ and $\tau$ in the right-hand side are as
in \eqref{s.t.e.ar.h}. The enhanced functor $\tau$ in the left-hand
side of \eqref{enh.comma.sq} has a right-adjoint $\eta:\C \to \E
\setminus^h_\pi \C$ induced by $\eta$ of \eqref{s.t.e.ar.h}, and
dually, $\sigma$ in the left-hand side of \eqref{enh.comma.sq} has a
left-adjoint $\eta:\C \to \C /^h_\pi \E$. The enhanced functor
$\pi:\C \to \E$ decomposes as
\begin{equation}\label{enh.comma.facto}
\begin{CD}
\C @>{\eta}>> \C /^h_\pi \E @>{\tau}>> \E,\qquad
\C @>{\eta}>> \E \setminus^h_\pi \C @>{\sigma}>> \E,
\end{CD}
\end{equation}
an enhanced version of the factorizations \eqref{comma.facto}. We
also have commutative squares
$$
\begin{CD}
\C @<{\sigma}<< \Ar^h(\C) @>{\tau}>> \C\\
@V{\pi}VV @VV{\Ar^h(\pi)}V @VV{\pi}V\\
\E @<{\sigma}<< \Ar^h(\E) @>{\tau}>> \E,
\end{CD}
$$
and by Corollary~\ref{semi.corr}~\thetag{ii}, these induce enhanced
functors
\begin{equation}\label{pi.id}
\id /^h \pi:\Ar^h(\C) \to \C /^h_\pi \E, \qquad \pi \setminus^h
\id:\Ar^h(\C) \to \E \setminus^h_\pi \C.
\end{equation}
For any enhanced object $e \in \E_\ppt$ of the enhanced category
$\E$, with corresponding enhanced functor $\eps^h(e):\ppt^h \to \E$,
the {\em enhanced fiber} $\C_e$ of the functor $\pi$ is given by
$\C_e = \ppt^h \times^h_\E \C$, and the {\em enhanced left
  resp.\ right comma-fibers} are the enhanced fibers $\C /^h_\pi e =
(\C /^h_\pi \E)_e$, $e \setminus^h_\pi \C = (\E \setminus^h_\pi
\C)_e$ of the projections $\tau$ resp.\ $\sigma$ of
\eqref{enh.comma.facto}. As in the unenhanced setting, we drop $\pi$
from notation when it is clear from the context.

\begin{defn}\label{enh.fib.def}
An enhanced functor $\pi:\C \to \E$ between enhanced categories
$\C$, $\E$ is an {\em enhanced fibration} if there exists a
left-admissible enhanced full subcategory $\E \setminus^h_\pi \C
\subset \Ar^h(\C)$ such that the commutative square
\begin{equation}\label{enh.fib.sq}
\begin{CD}
\E \setminus^h_\pi \C @>{\tau}>> \C\\
@V{\Ar^h(\pi)}VV @VV{\pi}V\\
\Ar^h(\E) @>{\tau}>> \E
\end{CD}
\end{equation}
is semicartesian. For any two enhanced fibrations $\C,\C' \to \E$,
an enhanced functor $\gamma:\C \to \C'$ over $\E$ is {\em cartesian}
if $\Ar^h(\gamma)$ sends $\E \setminus \C \subset \Ar^h(\C)$ into $E
\setminus \C' \subset \Ar^h(\C')$. Dually, an enhanced functor $\pi$
is an {\em enhanced cofibration} if $\pi^\iota$ is an enhanced
fibration, and a functor $\gamma:\C \to \C'$ over $\E$ for two
enhanced cofibrations $\C,\C' \to \E$ is {\em cocartesian} if
$\gamma^\iota$ is cartesian. An enhanced functor $\pi$ is an {\em
  enhanced bifibration} if it is both an enhanced fibration and an
enhanced cofibration.
\end{defn}

Note that the enhanced categories in Definition~\ref{enh.fib.def}
are not required to be small. Still, if $\pi:\C \to \E$ is an
enhanced fibration, then the enhanced full subcategory $\E
\setminus_\pi \C \subset \Ar^h(\C)$ is unique (\cite[Corollary
  7.4.2.3]{big}), so being an enhanced fibration is a condition and
not a structure. If $\C$ and $\E$ are small, then by
Corollary~\ref{semi.corr}~\thetag{ii}, $\E \setminus^h_\pi \C$ is
the right enhanced comma-category \eqref{enh.comma.sq}, and this
explains our notation. The enhanced functor $\Ar^h(\C) \to \E
\setminus^h_\pi \C$ left-adjoint to the embedding is then the
enhanced functor $\pi \setminus^h \id$ of \eqref{pi.id}, so that an
enhanced functor $\pi:\C \to E$ is an enhanced fibration iff $\pi
\setminus^h \id$ is right-reflexive.  For an enhanced cofibration
$\pi:\C \to \E$, we similary denote $\C /^h_\pi \E = (\E^\iota
\setminus^h_{\pi^\iota} \C^\iota)^\iota$; if $\C$ and $\E$ are
small, this is the left enhanced comma-category of
\eqref{enh.comma.sq}, and $\pi$ is an enhanced cofibration iff $\id
/^h \pi$ is left-reflexive.

\begin{exa}
Say that an enhanced functor $\pi:\C \to \E$ is an {\em enhanced
  family of groupoids} if the commutative square
$$
\begin{CD}
\Ar^h(\C) @>{\tau}>> \C\\
@V{\Ar^h(\pi)}VV @VV{\pi}V\\
\Ar^h(\E) @>{\tau}>> \E
\end{CD}
$$
is semicartesian. Then an enhanced family of groupoids is
taulotogically an enhanced fibration, with $\E \setminus^h \C =
\Ar^h(\C)$. If $\E = \ppt^h$, then $\C \to \ppt^h$ is an enhanced
family of groupoids iff $\C$ is an enhanced groupoid. In general,
for any enhanced fibration $\pi:\C \to \E$, we can let $\flat$ be
the closed class of enhanced morphisms in $\C$ that are cartesian
over $\E$; then \eqref{h.v.eq} provides an enhanced category
\begin{equation}\label{h.flat.sq}
  \C_{h\flat} = \Unf(\C_{\ppt,\flat}) \times_{\Unf(\C_\ppt)} \C,
\end{equation}
the induced enhanced functor $\C_{h\flat} \to \E$ is an enhanced
family of groupoids, the embedding $\C_{h\flat} \to \C$ is cartesian
over $\E$, and for any enhanced family of groupoids $\C' \to \E$, an
enhanced functor $\C' \to \C$ cartesian over $\E$ factors through
$\C_{h\flat} \to \C$, uniquely up to a unique isomorphism.
\end{exa}

For any enhanced functor $\pi:\C \to \E$ between small enhanced
categories, the enhanced functor $\sigma$ resp.\ $\tau$ of
\eqref{enh.comma.facto} is an enhanced fibration resp.\ cofibration
(\cite[Corollary 7.4.2.11]{big}). For any semicartesian square
\eqref{cat.sq} of enhanced categories and enhanced functors, if
$\gamma_1$ is a fibration or a cofibration, then so is
$\gamma^0_{01}$ (\cite[Lemma 7.4.2.15]{big}). If we have a
commutative square \eqref{cat.sq} such that $\gamma_1$ and
$\gamma^0_{01}$ are enhanced fibrations resp.\ cofibrations, then we
say that $\gamma^1_{01}$ is cartesian resp.\ cocartesian over
$\gamma_0$ if the induced functor $\C_{01} \to \C_0 \times^h_\C
\C_1$ is cartesian resp.\ cocartesian. If $\E = \Unf(I)$ for a small
category $I$, then semicartesian squares \eqref{enh.comma.sq} are in
fact cartesian, and $\pi:\C \to \E$ is an enhanced fibration if and
only it is a fibration in the usual sense (\cite[Lemma
  7.4.3.1]{big}). In particular, if $I$ is a partially ordered set,
then an enhanced functor $\C \to \Unf(I)$ is an enhanced fibration
resp.\ cofibration if and only if it is an $I$-augmentation
resp.\ $I$-coaugmentation. Thus for any small enhanced fibration
$\pi:\C \to \E$ and enhanced morphism $f:e \to e'$ in its base $\E$,
with the corresponding enhanced functor $\eps^h(f):\Unf([1]) \to
\E$, we have $\eps^h(f)^*\C \cong \Cyl^\iota_h(f^*)$ for a unique
  {\em transition functor} $f^*:\C_{e'} \to \C_e$ between enhanced
  fibers $\C_e$ and $\C_{e'}$. If we have another small enhanced
  fibration $\C' \to \E$ and a functor $\nu:\C' \to \C$, then for
  any $f:e \to e'$, we have a morphism
\begin{equation}\label{fu.fib.eq}
\nu_e \circ f^* \to f^* \circ \nu_{e'},
\end{equation}
see \cite[(7.4.3.3)]{big}, and $\nu$ is cartesian over $\E$ iff all
the maps \eqref{fu.fib.eq} are invertible. If this holds, we then
have commutative squares
\begin{equation}\label{fu.fib.sq}
\begin{CD}
\C'_{e'} @>{\nu_{e'}}>> \C_{e'}\\
@V{f^*}VV @VV{f^*}V\\
\C'_e @>{\nu_e}>> \C_e.
\end{CD}
\end{equation}
Dually, if we are given a small enhanced cofibration $\pi:\C \to
\E$, then for any enhanced morphism $f:e \to e'$, we have
$\eps^h(f)^*\C \cong \Cyl_h(f_!)$ for a unique transition functor
$f_!:\C_e \to \C_{e'}$, and for any enhanced functor $\nu:\C' \to
\C$ over $\E$ between small enhanced cofibrations $\C',\C \to \E$,
we have a morphism $f_!  \circ \nu_e \to \nu_{e'} \circ f_!$. The
enhanced functor $\nu$ is cocartesian over $\E$ iff all these
morphisms are invertible.

\begin{lemma}
Let $\pi:\C \to \E$ be a small enhanced fibration, with enhanced
fibers $\C_e$, $e \in \E_\ppt$ and transition functors $f^*$.
\begin{enumerate}
\item Assume given an enhanced full subcategory $\C' \subset \C$
  such that for any enhanced morphism $f:e \to e'$ in $\E$,
  $f^*:\C_{e'} \to \C_e$ sends $\C'_{e'}$ into $\C'_e$. Then the
  induced functor $\C' \to \E$ is an enhanced fibration, and the
  embedding enhanced functor $\nu:\C' \to \C$ is cartesian over
  $\E$. Conversely, if we have an enhanced fibration $\C' \to \E$
  and an enhanced functor $\nu:\C' \to \C$ cartesian over $\E$, then
  $\nu$ is fully faithful resp.\ an equivalence if and only if so
  are its fibers $\gamma_e:\C'_e \to \C_e$ for all enhanced objects
  $e \in \E_\ppt$.
\item Assume given a small enhanced fibration $\C' \to \E$, and an
  enhanced functor $\gamma:\C' \to \C$ cartesian over $\E$. Then if
  $\gamma_e:\C'_e \to \C_e$ is left-reflexive for any $e \in
  \E_\ppt$, the enhanced functor $\gamma$ is itself left-reflexive
  over $\E$, and the maps \eqref{fu.fib.eq} for the adjoint enhanced
  functor $\gamma^\dg$ are the base change maps for the
  corresponding commutative squares \eqref{fu.fib.sq}.
\item The functor $\pi$ is an enhanced bifiration if and only if
  $f^*$ is left-reflexive for any $f:e \to e'$.
\end{enumerate}
\end{lemma}

\proof{} \thetag{i} is \cite[Lemma 7.4.3.2]{big}, \thetag{ii} is
\cite[Lemma 7.4.3.3]{big}, and \thetag{iii} is \cite[Corollary
  7.4.3.5]{big}. \endproof

For any small enhanced category $\E$, the forgetful functors
\eqref{ccat.h.bb} are enhanced fibrations; for any $\C \in \Cat^h$,
the enhanced fiber $(\cCat^h \bb^h \E)_\C$ is the enhanced functor
category $\fFun^h(\C,\E)$ (\cite[Proposition 7.4.3.6]{big}). Note
that while the enhanced categories $\cCat^h \bbi^h \E$, $\cCat^h
\bb^h \E$ themselves are not small, the forgetful functors
\eqref{ccat.h.bb} are small as soon as so is $\E$. The forgetful
functors \eqref{ccat.h.bbd} are enhanced cofibrations, with enhanced
fibers $(\E \bbd^h \cCat^h)_\C \cong \fFun^h(\E,\C)$ (this is also
in \cite[Proposition 7.4.3.6]{big}). In particular, we have an
enhanced cofibration
\begin{equation}\label{ccat.h.dot}
\cCat^h_\idot \to \cCat^h
\end{equation}
with enhanced fiber $(\cCat^h_\idot)_\C \cong \C$ for any $\C \in
\cCat^h$. If we restrict our attention to $\cCat \subset \cCat^h$,
then \eqref{ccat.h.dot} induces an enhanced cofibration $\cCat_\idot
\to \cCat$, and then \eqref{cat.loc} extends to a semicartesian
square of enhanced cofibrations. Thus if we consider the
localization $\Cat^0$ as an enhanced category, as opposed to simply
a category, it behaves well.

\subsection{The Yoneda package.}

Assume given an enhanced category $\E$, and let $\Cat^h \bbi^h \E
\subset \Cat^h \bb^h \E$ be the subcategory with the same objects,
and morphisms given by isomorphism classes of enhanced functors over
$\E$. Moreover, if $\E$ is small, let $\Cat^h \mmh \E \subset \Cat^h
\bbi^h \E$ be the subcategory of enhanced fibrations $\pi:\C \to
\E$, with morphisms given by isomorphism classes of
enhanced-cartesian functors over $\E$. Then as in \eqref{3.enh.dia},
the enhancement $\cCat^h \bb^h \E$ for the category $\Cat^h \bb \E$
of Proposition~\ref{ccat.h.prop} induces enhancements for the
categories $\Cat^h \mmh \E \subset \Cat^h \bbi^h \E \subset \Cat^h
\bb^h \E$ by cartesian squares
\begin{equation}\label{3.h.enh.dia}
\begin{CD}
\cCat^h \mmh \E @>{a}>> \cCat^h \bbi^h \E @>{b}>> \cCat^h \bb^h \E\\
@VVV @VVV @VVV\\
\Unf(\Cat^h \mmh \E) @>>> \Unf(\Cat^h \bbi^h \E) @>>> \Unf(\Cat^h \bb^h
\E),
\end{CD}
\end{equation}
where again, the vertical arrow on the right is the truncation
functor \eqref{trunc.eq}, and the categories on the left are only
considered if $\E$ is small. Explicitly --- see \cite[Lemma
  7.4.4.1]{big} for a proof --- $\cCat^h \bbi^h \E \subset \cCat^h
\bb^h \E$ is the subcategory of objects $\langle I,\C,\alpha
\rangle$ such that $\alpha$ is coaugmented over the projection $I^o
\to \ppt$, and with morphisms given by triples $\langle f,\phi,a
\rangle$ with invertible $a$. Any such $\alpha$ is of the form
$\alpha \cong \beta_{h\perp}$ for a unique enhanced functor
$\beta:\C \to \E$ augmented over $I \to \ppt$, so that
alternatively, $\cCat^h \bbi^h \E$ is the category of triples $\langle
I,\C,\beta \rangle$, $\beta:\C \to \E$ an enhanced functor augmented
over $I \to \ppt$, and morphisms are isomorphism classes of
commutative squares
\begin{equation}\label{ccat.h.1.sq}
\begin{CD}
\C' @>{\phi}>> \C\\
@V{\pi' \times \beta'}VV @VV{\pi \times \beta}V\\
\Unf(I') \times^h \E @>{\Unf(f) \times \id}>> \Unf(I) \times^h \E
\end{CD}
\end{equation}
such that $\phi$ is augmented over $f$. If $\E$ is small, then
$\cCat^h \mmh \E \subset \cCat^h \bbi^h \E$ is the subcategory of
triples $\langle I,\C,\beta \rangle$ such that $\pi \times \beta:\C
\to \Unf(I) \times^h \E$ is an enhanced fibration, with morphisms
given by isomorphism classes of squares \eqref{ccat.h.1.sq} such
that $\phi$ is enhanced-cartesian over $\Unf(f) \times \id$.

Assume that $\E$ is small, and consider the enhanced functor
category $\E^\iota\cCat^h = \fFun^h(\E^\iota,\cCat^h)$ of
Corollary~\ref{fun.h.corr}, Explicitly, objects in $\E^\iota\cCat^h$
are pairs $\langle I,\gamma \rangle$, $I \in \Posf^+$,
$\gamma:\E^\iota \times^h \Unf(I^o) \to \cCat^h$ an enhanced
functor. Note that the involution $\C \mapsto \C^\iota$ is
functorial with respect to $\C \in \Cat^h$, and the equivalence
$\iota:\Cat^h \to \Cat^h$, $\C \mapsto \C^\iota$ has a natural
enhancement
\begin{equation}\label{iota.enh}
\iota_h:\cCat^h \to \cCat^h, \qquad \langle I,\C \rangle \mapsto
\langle I,\C^\iota_{h\perp} \rangle.
\end{equation}
Then for every $\langle I,\gamma \rangle \in \fFun^h(\E^\iota,\C)$,
one can construct a semicartesian square
\begin{equation}\label{G.pr.sq}
\begin{CD}
\C(\gamma) @>>> \cCat^h_\idot\\
@V{\pi(\gamma)}VV @VVV\\
\E^\iota \times^h \Unf(I^o) @>{\iota_h \circ \gamma}>> \cCat^h,
\end{CD}
\end{equation}
where the arrow on the right is the enhanced cofibration
\eqref{ccat.h.dot}, and then $\pi(\gamma)^\iota:\C(\gamma)^\iota \to
\E \times^h \Unf(I)$ is an enhanced fibration, thus defines an
object in $\cCat^h \mmh \E$. This is obviously functorial with
respect to maps cartesian over $\Posf^+$, and we can then construct
an enhanced functor
\begin{equation}\label{G.pr.eq}
\fFun(\E^\iota,\cCat^h) \to \cCat^h \mmh \E
\end{equation}
by applying our general extension result, Lemma~\ref{fl.le}.

\begin{prop}\label{groth.prop}
For any small enhanced category $\E$, the enhanced functor
\eqref{G.pr.eq} is an equivalence.
\end{prop}

\proof{} This is \cite[Proposition 7.4.4.2]{big}. \endproof

Proposition~\ref{groth.prop} is an enhanced version of the
Grothendieck construction, and we see that the statement is somewhat
cleaner: there is no need to introduce ``pseudofunctors'' (this is
already embedded into the enhanced functor formalism and the
definition of the enhanced category $\cCat^h$). Dually, one can also
consider the enhanced subcategory $\cCat^h \mch \E \subset \cCat^h
\bbi^h \E$ of enhanced cofibrations and functors enhanced-cocartesian
over $\E$; then \eqref{G.pr.eq} for $\E^\iota$ immediately provides
an equivalence
\begin{equation}\label{G.iota.eq}
\cCat^h \mch \E \cong \fFun^h(\E,\cCat^h).
\end{equation}
Moreover, enhanced functors on the right-hand side of
\eqref{G.pr.eq} can be composed with the involution
\eqref{iota.enh}, and this provides a non-trivial operation on the
left-hand side: for any enhanced fibration $\C \to \E$, we obtain a
transpose-opposite enhanced fibration $\C^\iota_\perp$ and transpose
enhanced cofibration $\C_{h\perp} \to \E$. If $\E = \Unf(I)$ for a
partially ordered set $I$, we recover the transpose-opposite
augmentation and transpose coaugmentation.

We can now use the equivalence \eqref{G.pr.eq} to construct an enhanced
version of the Yoneda embedding \eqref{yo.eq}. What we need for this
is an enhanced version
\begin{equation}\label{3.h.dia}
\begin{CD}
\cCat^h \mmh \E @>{a}>> \cCat^h \bbi^h \E @>{b}>> \cCat^h \bb^h \E
@>{y}>> \cCat^h \mmh \E
\end{CD}
\end{equation}
of the diagram \eqref{3.dia} of Subsection~\ref{yo.subs}. We already
have all the categories and the embeddings $a$, $b$, so it remains
to construct $y$. This is provided by the following.

\begin{prop}\label{yo.prop}
For any small enhanced category $\E$, there exists a fully faithful
enhanced functor $y:\cCat^h \bb^h \E \to \cCat^h \mmh \E$ and
functorial maps $y \circ b \circ a \to \id$, $\id \to a \circ y
\circ b$, $b \circ a \circ y \to \id$ that define an adjunction both
between $a \circ y$ and $b$ and between $a$ and $y \circ b$.
\end{prop}

\proof{} This is \cite[Proposition 7.4.6.1]{big}. \endproof

The conditions of Proposition~\ref{yo.prop} define the functor $y$
uniquely. To save space, we do not reproduce the full construction
here; however, Proposition~\ref{yo.prop} is not deep at all. On the
level of enhanced objects, $y$ sends some enhanced functor $\pi:\C
\to \E$ to the enhanced fibration $\sigma:\E \setminus^h_\pi \C \to
\E$ of \eqref{enh.comma.facto}. Just as in Subsection~\ref{yo.subs},
the adjunction maps are given by the enhanced functor $\eta:\C \to
\E \setminus^h_\pi \C$ of \eqref{enh.comma.facto}, its left-adjoint
$\tau:\E \setminus^h_\pi \C \to \C$, and its right-adjoint $\eta^\dg
= \sigma \circ (\pi \setminus^h \id)^\dg:\E \setminus^h_\pi \C \to
\C$ that exists as soon as $\pi$ is an enhanced fibration. The
required universal properties easily follow from
Corollary~\ref{semi.corr}, and then checking that $y$ and the
adjunction maps are functorial with respect to $\pi:\C \to I$ and
that $y$ is fully faithful becomes a purely formal
exercise. Combining the latter with the Grothendieck construction
equivalence \eqref{G.pr.eq}, we obtain a fully faithful enhanced
embedding
\begin{equation}\label{yo.h.cat.eq}
\Y:\cCat^h \bb^h \E \to \E^\iota\cCat^h.
\end{equation}
This is the enhanced version of the extended Yoneda embedding
\eqref{yo.cat.eq}. If we restrict out attention to the full enhanced
subcategory $\sSets^h \subset \cCat^h$ spanned by enhanced
groupoids, then \eqref{yo.h.cat.eq} restricts to a fully faithful
embedding $\sSets^h \bb^h \E \to \fFun^h(\E^\iota,\sSets^h)$, and if
we consider the enhanced fiber $(\sSets^h \bb^h \E)_{\ppt^h} \cong
\E$ over $\ppt^h \in \sSets^h \subset \cCat^h$, we obtain an
enhanced version
\begin{equation}\label{yo.h.eq}
\Y:\E \to \E^\iota\sSets^h
\end{equation}
of the usual Yoneda embedding \eqref{yo.eq}. It is still fully
faithful, and corresponds to an enhanced version
\begin{equation}\label{yo.h.pair}
\E^\iota \times^h \E \to \sSets^h
\end{equation}
of the usual $\Hom$-pairing. In particular, for any enhanced objects
$e$, $e'$ of the enhanced category $\E$, we obtain a small
``enhanced groupoid of maps'' $e \to e'$; the set of its connected
components is $\Hom_{\E_\ppt}(e,e')$.

As another application of Proposition~\ref{groth.prop}, one can
construct enhanced versions of the relative functor categories of
Subsection~\ref{fib.subs}. This reverses \eqref{cat.bb}. Namely,
assume given an enhanced cofibration $\gamma:\E' \to \E$ between
small enhanced categories. For any small enhanced category $\C$,
define the {\em relative functor category} $\fFun^h(\E'|\E,\C)$ by
the semicartesian square
\begin{equation}\label{enh.rel.dia}
\begin{CD}
\fFun^h(\E'|\E,\C) @>>> \cCat^h \bb^h \C\\
@VVV @VVV\\
\E @>{X}>> \cCat^h,
\end{CD}
\end{equation}
where $X:\E \to \cCat^h$ corresponds to $\E' \to \E$ by the
covariant Grothnedieck construction \eqref{G.iota.eq}. Then as in
the unenhanced setting, we have an enhanced {\em evaluation pairing}
$\ev:\E' \times^h_{\E} \fFun^h(\E'|\E,\C) \to \C$ such that for any
small enhanced category $\C'$ over $\C$, any enhanced functor
$\gamma:\E' \times^h_{\E} \C' \to \C$ factors as
\begin{equation}\label{rel.enh.fun.facto}
\begin{CD}
\E' \times^h_{\E} \C' @>{\id \times^h \wt{\gamma}}>> \E' \times_{\E}
\Fun(\E'|\E,\C)
@>{\ev}>> \C
\end{CD}
\end{equation}
for an enhanced functor $\wt{\gamma}:\C' \to \Fun(\E'|\E,\C)$ over
$\E$, uniquely up to an isomorphism. This follows from \cite[Lemma
  7.4.5.1]{big} that actually proves more. Namely, the functor
$\gamma^*:\Cat^h \bbi^h \E \to \Cat^h \bbi^h \E'$, $\C \mapsto \C
\times^h_{\E} \E'$ has an obvious tautological left-adjoint functor
$\gamma_{\trr}:\Cat^h \bbi^h \E' \to \Cat^h \bbi^h \E$ sending a
category $\C$ with an enhanced functor $\C \to \E'$ to the same
category $\C$ with the composition $\C \to \E' \to \E$. Then it is
easy to see that the adjoint pair $\langle \gamma^*,\gamma_\trr
\rangle$ has a natural enhancement --- see \cite[Subsection
  7.4.5]{big} --- and \cite[Lemma 7.4.5.1]{big} proves that
$\gamma^*:\cCat^h \bbi^h \E \to \cCat^h \bbi^h \E'$ also has a
right-adjoint enhanced functor $\gamma_\trl:\cCat^h \bbi^h \E' \to
\cCat^h \bbi^h \E$. The relative functor category
\eqref{enh.rel.dia} is recoved as $\fFun^h(\E'|\E,\C) \cong
\gamma_\trl(\C \times \E')$, and the universal property follows from
the adjunction between $\gamma^*$ and $\gamma_\trl$, with the
evaluation pairing given by the adjunction map.

If $\gamma:\E' \to \E$ is an enhanced fibration rather than an
enhanced cofibration, one can achieve the same results by defining
\begin{equation}\label{rel.fun.perp}
\fFun^h(\E'|\E,\C) = \fFun^h({\E'}^\iota|\E^\iota,\C^\iota)^\iota
\cong \fFun^h(\E'_{h\perp}|\E^\iota,\C)_{h\perp},
\end{equation}
where $\gamma^\iota:{\E'}^\iota \to \E^\iota$ is the
enhanced-opposite cofibration. Then the evaluation pairing for
$\gamma^\iota$ provides an enhanced functor $\ev^\iota:\gamma^{\iota
  *}\fFun^h({\E'}^\iota|\E^\iota,\C^\iota) \to \C^\iota$ that gives
evaluation pairing for $\gamma$ after going back to the
enhanced-opposite categories, with the same universal property as in
the cofibration case.

\subsection{Limits and Kan extensions.}

Let us now discuss limits and colimits in the enhanced context; the
story here is largely parallel to Subsection~\ref{lim.subs}.

An {\em enhanced cone} of an enhanced functor $E:\C \to \E$ between
enhanced categories $\C$, $\E$ is an enhanced functor $E_>:\C^{h>}
\to \E$ equipped with an isomorphism $\eps^*E_> \cong E$, where
$\C^{h>}$ is the enhanced category \eqref{h.gt.lt}, and $\eps:\C \to
\C^{h>}$ is the embedding. The {\em vertex} of a cone $E_>$ is the
enhanced object in $\E$ corresponding to $E_> \circ o:\ppt^h \to
\E$. If $\C$ and $\E$ are small, all enhanced cones of some $E$ form
an enhanced category $\cCone(E) = E \setminus^h_{\gamma^*} \E$,
where $\gamma:\C \to \ppt^h$ is the tautological projection, and a
cone is {\em universal} if it is the initial enhanced object in
$\cCone(E)$. Since semicartesian products preserve fully faithful
embeddings, a universal cone $E_>$ that factors through a full
enhanced subcategory $\E' \subset \E$ is also universal as a cone
for a functor to $\E'$. If the target enhanced category $\E$ is not
small, we say that a cone is {\em universal} if it is universal as a
cone for a functor to any small full enhanced subcategory $\E'
\subset \E$ through which it factors.  If a universal enhanced cone
exists, its vertex is called the {\em enhanced colimit} of the
functor $E$, and we denote it by $\colim^h E \in \E_\ppt$. An
enhanced category $\E$ is {\em enhanced-cocomplete} if any enhanced
functor $\C \to \E$ from a small enhanced category $\C$ has an
enhanced colimit.

More generally, for any functor $\gamma:\C \to \C'$ between small
enhanced categories, a {\em relative enhanced cone} of an enhanced
functor $E:\C \to \E$ is an enhanced functor $E_>:\Cyl_h(\gamma) \to
\E$ equipped with an isomorphism $s^*E_> \cong E$, where
$\Cyl_h(\gamma)$ is the enhanced cylinder of \eqref{enh.cyl.sq}, and
$s:\C \to \Cyl_h(\gamma)$ is the natural embedding. Again, if $\E$
is small, relative enhanced cones form a category $\cCone(E,\gamma)
= E \setminus^h_{\gamma^*} \fFun^h(\C',\E)$, and a relative enhanced
cone $E_>$ is {\em universal} if it is the initial enhaced object in
$\cCone(E,\gamma)$. If $\E$ is not small, then $E_>$ is universal if
it is universal as a functor to any small enhanced subcategory $\E_0
\subset \E$ through which it factors. If a universal relative cone
exists, then the {\em enhanced left Kan extension} $\gamma^h_!E$ is
$t^*E_>$.

An enhanced functor $F:\E \to \E'$ to some enhanced category $\C'$
{\em preserves the enhanced left Kan extension} $\gamma^h_!E$
corresponding to a universal enhanced relative cone $E_>$ if the
composition $F \circ E_>$ is a universal relative cone. An enhanced
left Kan extension $\gamma^h_!E$ is {\em universal} if it is
preserved by the Yoneda embedding $\Y:\E_0 \to \E_0^\iota\sSets^h$
for any small full enhanced subcategory $\E_0 \subset \E$ through
which it factors. All enhanced colimits are universal (\cite[Lemma
  7.5.3.8]{big}). An enhanced functor $\E \to \E'$ between
enhanced-cocomplete enhanced categories is {\em
  enhanced-right-exact} if it preserves all enhanced colimits.

\begin{lemma}\label{h.adj.kan.le}
For any enhanced category $\E$ and left-reflexive enhanced fun\-ctor
$\gamma:\C \to \C'$ between small enhanced categories $\C$, $\C'$,
with left-adjoint $\gamma_\dg:\C' \to \C$, the enhanced functor
$\gamma^*:\fFun^h(\C',\E) \to \fFun^h(\C,\E)$ is right-adjoint to
the enhanced functor $\gamma_\dg^*:\fFun^h(\C,\E) \to
\fFun^h(\C',\E)$, and for any enhanced functor $E:\C \to \E$, the
enhanced left Kan extension $\gamma^h_!E$ exists and is functorially
isomorphic to $\gamma_\dg^*E$.
\end{lemma}

\proof{} This is \cite[Lemma 7.5.3.6]{big} jointly with
\cite[Corollary 7.4.1.7]{big}. \endproof

\begin{prop}\label{h.cmpl.kan.prop}
Assume given an enhanced-cocomplete enhanced category $\E$. Then for
any enhanced functor $\gamma:\C \to \C'$ between enhanced small
categories, and any enhanced functor $E:\C \to \E$, an enhanced left
Kan extension $\gamma^h_!E$ exists. The enhanced functor
$\gamma^*:\fFun^h(\C',\E) \to \fFun^h(\C,\E)$ admits a left-adjoint
enhanced functor $\gamma^h_!:\fFun^h(\C,\E) \to \fFun^h(\C',\E)$
whose fiber $\gamma^h_!:\Fun^h(\C,\E) \to \Fun^h(\C',\E)$ over $\ppt
\in \Posf^+$ sends $E:\C \to \E$ to $\gamma^h_!E:\C' \to
\E$. Furthermore, for any diagram
\begin{equation}\label{h.bc.dia}
\begin{CD}
\C'_0 @>{\gamma'}>> \C'_1 @>{\pi'}>> \C'\\
@V{\phi_0}VV @VV{\phi_1}V @VV{\phi}V\\
\C_0 @>{\gamma}>> \C_1 @>{\pi}>> \C
\end{CD}
\end{equation}
of small enhanced categories and enhanced functors such that the
squares are commutative and semicartesian, and both $\pi$ and $\pi
\circ \gamma$ are enhanced cofibrations, the base change map
\begin{equation}\label{h.kan.bc.eq}
  {\gamma'}^h_! \circ \phi_0^* \to \phi_1^* \circ \gamma^h_!
\end{equation}
of enhanced functors $\fFun^h(\C_0,\E) \to \fFun^h(\C',\E)$ is an
isomorphism.
\end{prop}

\proof{} Combine \cite[Corollary 7.5.3.11]{big}, \cite[Lemma
  7.5.3.6]{big} and \cite[Proposition 7.5.1.9]{big}. \endproof

\begin{corr}\label{h.cmpl.corr}
  For any small enhanced category $\C$, and any enhanced-cocomplete
  enhanced category $\E$, the enhanced functor category
  $\fFun^h(\C,\E)$ is enhanced-cocomplete.
\end{corr}

\proof{} This is \cite[Corollary 7.5.1.10]{big}. \endproof

As in the unenhanced setting, Lemma~\ref{h.adj.kan.le},
Proposition~\ref{h.cmpl.kan.prop} and the decomposition
\eqref{enh.comma.facto} provide an isomorphism
\begin{equation}\label{h.kan.eq}
\gamma^h_!E(c) \cong \colim^h_{\C' /^h_\gamma c}E
\end{equation}
for any enhanced functor $\gamma:\C' \to \C$ between small enhanced
categories, enhanced functor $E:\C' \to \E$ to an enhanced
cocomplete target enhanced category $\E$, and enhanced object $c \in
\C_\ppt$.

For any enhanced functor $\gamma:\C \to \C'$ between small enhanced
categories, and any enhanced functor $E:\C \to \E$ that admits a
universal enhanced left Kan extension $\gamma^h_!E:\C' \to \E$, we
have an enhanced map $a:E \to \gamma^*\gamma^h_!E$, and the pair
$\langle \gamma_!E,a \rangle$ has the universal property:
\begin{enumerate}
\renewcommand{\labelenumi}{(\Roman{enumi})}
\item\label{I.it} for any enhanced
functor $E':\C' \to \E$ equipped with an enhanced morphism $a':E \to
\gamma^*E'$, there exists a unique enhanced morphism $b:\gamma^h_!E
\to E'$ such that $\gamma^*(b) \circ a = a'$.
\end{enumerate}
By itself, \thetag{I} does not insure that $\gamma^h_!E$ is an
enhanced Kan extension. Indeed, if $o \in \E_\ppt$ is a terminal
enhanced object in an enhanced category $\E$, then it is also a
terminal object in $\E_\ppt$ in the unenhanced sense, but the
converse is wrong --- there are plenty of enhanced categories such
that $\E_\ppt$ has a terminal object, but $\E$ does not have a
terminal enhanced object at all (e.g.\ take any non-trivial
enhancement for $\E_\ppt = \ppt$). However, once we know that a
universal left Kan extension exists, \thetag{I} characterizes it
completely. This allows to define enhanced left Kan extensions
$\gamma^h_!E$ along enhanced functors $\gamma:\C \to \C'$ whose
target is not necessairily small: one says (\cite[Subsection
  7.5.4]{big}) that $E:\C \to \E$ {\em admits a universal enhanced
  left Kan extension} $\gamma^h_!E$ if for any small enhanced full
subcategory $\C'_0 \subset \C'$ such that $\gamma$ factors through
an enhanced functor $\gamma_0:\C \to \C'_0$, $E$ admits a universal
enhanced left Kan extension $\gamma^h_{0!}E:\C'_0 \to \E$. In this
case, \thetag{I} insures that the enhanced functors $\gamma^h_{0!}E$
for various $\C'_0 \subset \C'$ patch together to a single enhanced
functor $\gamma^h_!E:\C' \to \E$; it comes equipped with an enhanced
map $a:E \to \gamma^*\gamma^h_!E$, and still has the universal
property \thetag{I}.

Dually, the {\em enhanced limit} $\lim^h E$ of an enhanced functor
$E:\C \to \E$, $\C$ small is given by $\lim^h E = \colim^h E^\iota$,
if it exists. The {\em enhanced right Kan extension} $\gamma^h_*E$
is $\gamma^h_*E = (\gamma^h_!E^\iota)^\iota$, if it exists, and an
enhanced category $\E$ is {\em enhanced-complete} if $\E^\iota$ is
enhanced-cocomplete. We also have the obvious dual counterparts of
Lemma~\ref{h.adj.kan.le}, Proposition~\ref{h.cmpl.kan.prop},
Corollary~\ref{h.cmpl.corr} and the isomorphisms \eqref{h.kan.eq}.

\begin{remark}
Technically, the definitions of an enhanced-complete and
enhanced-cocomplete enhanced categories given in \cite[Definition
  7.5.1.1]{big} are different from the ones we use here (but the two
are equivalent by \cite[Corollary 7.5.3.11]{big}).
\end{remark}

\begin{prop}\label{cmpl.exa.prop}
The enhanced categories $\sSets^h$ and $\cCat^h$ are both
en\-hanced-complete and enhanced-cocomplete, and the embedding
$\sSets^h \subset \cCat^h$ is both right and left-admissible. If a
model category $\C$ is complete resp.\ cocomplete, then its enhanced
localization $\Hh^W(\C)$ of Proposition~\ref{mod.prop} is
enhanced-complete resp.\ enhanced-cocomplete.
\end{prop}

\proof{} The first claim is \cite[Example 7.5.1.6]{big} and
\cite[Proposition 7.5.2.1]{big}. The second claim is \cite[Example
  7.5.1.8]{big}. \endproof

\begin{exa}
By Proposition~\ref{cmpl.exa.prop} and
Proposition~\ref{h.cmpl.kan.prop}, for any enhanced functor
$\gamma:\C' \to \C$ from a small enhanced category $\C'$, an
enhanced functor $X:\C' \to \sSets^h$ admits an enhanced left Kan
extension $\gamma^h_!X:\C \to \sSets^h$. In particular, if $\C' =
\ppt^h$, so that $\gamma=\eps^h(c)$ corresponds to an enhanced
object $c \in \C_\ppt$, we obtain a canonical enhanced functor
\begin{equation}\label{yo.c.eq}
\Y(c) = \eps^h(c)^h_!\ppt^h:\C \to \sSets^h,
\end{equation}
where $\ppt^h$ stands for the constant enhanced functor
$\eps^h(\ppt^h):\ppt^h \to \sSets^h$. If $\C$ is small, then
\eqref{yo.c.eq} is obtained by applying the enhanced Yoneda
embedding $\Y:\C^\iota \to \fFun^h(\C,\sSets^h)$ of \eqref{yo.h.eq}
to $c \in \C^\iota_\ppt$; however, even if $\C$ is not small,
\eqref{yo.c.eq} is well-defined.
\end{exa}

Explicitly, enhanced limits and colimits in $\Hh^W(\C)$ are given by
homotopy limits and colimits provided by the standard Quillen
Adjunction Theorem. For $\cCat^h$, the following description of
$\lim^h$ and $\colim^h$ is given in \cite[Lemma 7.5.2.9]{big}. Let
$\C$ be a small enhanced category equipped with an enhanced functor
$E:\C^\iota \to \cCat^h$ corresponding to an enhanced fibration $\E
= \C E \to \C$. Then $\lim^h_{\C^\iota}E \cong \sSec^{\Tot h}(\C,\E)
\subset \sSec^h(\C,\E)$ is the full enhanced subcategory spanned by
cartesian sections, and $\colim^h_{\C^\iota}E$ fits into an
enhanced-cocartesian square
\begin{equation}\label{enh.colim.sq}
\begin{CD}
\E_{h\flat} @>>> \E\\
@VVV @VVV\\
\Hh^\Tot(\E_{h\flat}) @>>> \colim^h_{\C^\iota}E,
\end{CD}
\end{equation}
where $\E_{h\flat} \to \C$ is the enhanced family of groupoids of
\eqref{h.flat.sq}, with its universal property, and
$\Hh^\Tot:\cCat^h \to \sSets^h$ is the {\em total localization
  functor} left-adjoint to the embedding $\sSets^h \subset \cCat^h$,

It is also useful to consder the situation when an enhanced category
$\E$ has some colimits but not all of them. Formally, for any full
subcategory $I \subset \Cat^h$, say that an enhanced category $\E$
is {\em $I$-cocomplete} if $\colim^h_\C E$ exists for any enhanced
functor $E:\C \to \E$ from a small enhanced category $\C \in I$, and
say that an enhanced functor $\E \to \E'$ between two $I$-cocomplete
enhanced categories is {\em $I$-right-exact} if it preserves all
these enhanced colimits. For example, for any regular cardinal
$\kappa$, we can consider the full subcategory $I = \Cat^h_\kappa
\subset \Cat^h$ of $\kappa$-bounded enhanced categories; then
$I$-cocomplete enhanced categories are called {\em
  $\kappa$-enhanced-cocomplete}, and {\em finitely
  enhanced-cocomplete} if $\kappa$ is the countable cardinal. In
general, we have the following extension result.

\begin{prop}\label{env.prop}
For any full subcategory $I \subset \Cat^h$ and enhanced category
$\E$, there exists an $I$-cocomplete enhanced category $\Env(\E,I)$
and an enhanced full embedding $\Y(\E,I):\E \to \Env(\E,I)$ such
that for any $I$-cocomplete enhanced category $\E'$, any enhanced
functor $\gamma:\E \to \E'$ factors as
\begin{equation}\label{env.facto}
\begin{CD}
\E @>{\Y(\E,I)}>> \Env(\E,I) @>{\gamma'}>> \E',
\end{CD}
\end{equation}
where $\gamma'$ is $I$-right-exact, and the factorization
\eqref{env.facto} is unique up to a unique isomorphism. The category
$\E$ is itself $I$-cocomplete if and only if $\Y(\E,I)$ is
right-reflexive, and in this case, if $\gamma=\id$, then $\gamma'$
in \eqref{env.facto} is right-adjoint to $\Y(\E,I)$.
\end{prop}

\proof{} This is \cite[Lemma 7.5.7.14]{big} and \cite[Lemma
  7.5.7.17]{big}. \endproof

\begin{exa}
Take $I = \Cat^h$. Then ``$I$-cocomplete'' is the same thing as
``enhanced-cocomplete'', and if $\E$ is small, then
$\Env(\E,I)=\E^\iota\sSets^h$, and $\Y(\E,I)$ is the Yoneda
embedding \eqref{yo.h.eq}. If $\E = \ppt^h$, then
Proposition~\ref{env.prop} says that for any enhanced-cocomplete
enhanced category $\E'$ and enhanced object $e \in \E'_\ppt$, the
corresponding enhanced functor $\eps^h(e):\ppt^h \to \E'$ extends to
an enhanced-right-exact enhanced functor $\sSets^h \to \E'$,
uniquely up to a unique isomorphism. Note, however, that
$\Env(\E,I)$ exists even for a large enhanced category $\E$.
\end{exa}

\begin{exa}\label{gt.exa}
Let $I=\{\emptyset\} \subset \Cat^h$; then an enhanced category $\E$
is $I$-cocomplete iff it has a terminal enhanced object, and
$\Env(\E,I) = \E^{h>}$ is the enhanced category \eqref{h.gt.lt}.
\end{exa}

\begin{exa}\label{kar.exa}
Let $P$ be the category with a single object $o$ that has a single
non-trivial endomorphism $p:o \to o$ with $p^2=p$. Consider the full
subcategory $\P = \{\Unf(P)\} \subset \Cat^h$. Then an enhanced
category is {\em enhanced-Karoubi-closed} if it is $\P$-cocomplete,
and $\Env(\E,\P)$ is the {\em enhanced Karoubi completion} of $\E$.
\end{exa}

\begin{remark}
The situation of Example~\ref{kar.exa} is quite special. In
particular, all enhanced functors are automatically $\P$-cocomplete
(\cite[Lemma 7.6.1.2]{big}, and because of this, $\Env(-,\P)$ is an
idempotent oparation. For other subcategories $I \subset \Cat^h$,
this is certainly not true. In Example~\ref{gt.exa}, $\Env(-,I)$
adds a new enhanced terminal object to the category $\E$; if it
already had one, it stops being terminal.
\end{remark}

\subsection{Colimits in $\cCat^h$ and localizations.}

One can summarize Proposition~\ref{h.cmpl.kan.prop} by saying that
everything works just as in the unenhanced setting, with the proviso
that the commutative squares in \eqref{h.bc.dia} are only
semicartesian, and all the functors have to be enhanced.

To illustrate the latter point, assume given a small category $I$
and an enhancement $\E$ for a category $\E_\ppt$. Then an enhanced
functor $\gamma:\Unf(I) \to \E$ induces a usual functor
$\gamma_\ppt:I \to \E_\ppt$, and we think of $\gamma$ as an
enhancement for $\gamma_\ppt$. If $I = [1]$, then since $\E$ is a
separated reflexive family, an enhancement for any $\gamma_\ppt$ is
unique up to an isomorphism, and since $\E$ is semiexact, the same
holds when $I = \V = [1] \copr_\ppt [1]$. If $I = [1]^2 \cong \V^>
\cong \V^{o<}$, the enhancement is no longer unique, but we can
still say that a commutative square $[1]^2 \to \E_\ppt$ is {\em
  homotopy cocartesian} resp.\ {\em homotopy cartesian} if it admits
an enhancement $\gamma:\Unf([1]^2) \to \E$ such that $\gamma$
resp.\ $\gamma^\iota$ is a universal enhanced cone. However, such
squares $[1]^2 \to \E_\ppt$ themselves have no universal properties,
and in particular, they are only unique up to an isomorphism
(determined by the choice of an enhancement for the original functor
$\V \to \E_\ppt$ resp.\ $\V^o \to \E_\ppt$). This behaviour is
completely parallel to the behaviour of cones in triangulated
categories, with the only exception that the collection of
distinguished triangles in a triangulated category has to be imposed
as an extra piece of data and axiomatized, while in the enhanced
setting, homotopy cartesian and cocartesian squares are determined
by enhancements.

In the particular case $\E = \cCat^h$, homotopy cartesian squares
are induced by semicartesian squares of Corollary~\ref{semi.corr}
(see \cite[Lemma 7.5.2.14]{big}), and those do have a universal
property. However, they are only universal with respect to
commutative squares of enhanced categories and enhanced functors,
and such a commutative square contains strictly more information
than an unenhanced commutative square $[1]^2 \to \Cat^h$. In fact,
for any partially ordered set $I$, one can define an {\em $I$-family
  of enhanced categories} as a fibration $\C \to \Pos^+ \times I$
whose restriction $\C_i$ to $\Posf^+ \times \{i\}$ is an enhanced
category for any $i \in \I$, and such that for any $i \leq i'$, the
transition functor $\C_{i'} \to \C_i$ is an enhanced functor. Then
such families form a category $\Cat^h(I|I)$, with morphisms given by
cartesian functors over $\Posf^+ \times I$ up to isomorphisms over
$\Posf^+ \times I$, and we have natural comparison functors
\begin{equation}\label{cat.I.eq}
\begin{CD}
\Cat^h(I) @>{\alpha}>> \Cat^h(I|I) @>{\beta}>> I^o\Cat^h,
\end{CD}
\end{equation}
where $\alpha$ is the pullback with respect to the functor $\Posf^+
\times I \to \Unf(I)$ sending $J \times i$ to the constant map $J
\to I$ with value $i$, and $\beta$ sends an $I$-family $\C$ to the
collection of $\C_i$ and transition functor $\C_{i'} \to \C_i$. In
general, both functors in \eqref{cat.I.eq} are very far from being
an equivalence. However, if $\dim I \leq 1$ --- in particular, for
$I=\V$ --- then $\alpha$ in \eqref{cat.I.eq} is an equivalence
(\cite[Corollary 7.3.5.3]{big}), while $\beta$ is an epivalence
(\cite[Lemma 4.1.1.9]{big}) but not necessarily an equivalence. Then
if $\dim I \leq 1$, $\lim^h_{I^o}$ factors through $\alpha$, and in
the case $I=\V$, $\lim_{V^o}$ is obtained by taking the
semicartesian product. However, unless $I$ is discrete, $\lim_{I^o}$
does not factor through $\beta$.

\begin{remark}
For a general enhanced limit $\lim^h_{I^o}E$ of an enhanced functor
$E:\Unf(I) \to \E$ to some enhanced category $\E$, one can pass to a
small full enhanced subcategory $\E_0 \subset \E$ that contains the
essential image of the corresponding universal enhanced dual cone,
and apply the Yoneda embedding $\Y:\E_0 \to \E_0^\iota\sSets^h$. Then
enhanced objects in $\E_0$ become enhanced families of groupoids
over $\E_0$, and as in the case $\E=\cCat^h$, when $I = \V$, the
limit is obtained by taking the semicartesian product of such
families.
\end{remark}

As a formal corollary of Proposition~\ref{cmpl.exa.prop}, $\Cat^h$
has arbitrary enhanced localizations: for any small enhanced
category $\C$ and a set $W$ of enhanced morphisms $\Unf([1]) \to
\C$, we have an enhanced-cocartesian square
\begin{equation}\label{enh.loc.sq}
\begin{CD}
W \times \Unf([1]) @>{w}>> \C\\
@VVV @VVV\\
W \times \ppt^h @>>> \Hh^W(\C)
\end{CD}
\end{equation}
of enhanced categories and enhanced functors that comes from a
universal enhanced cone in $\Cat^h$, and has the same universal
properties as \eqref{loc.sq} and \eqref{pre.sq}. When $W=\Tot$ is
the set of (isomorphism classes of) all enhanced morphisms, then
$\Hh^\Tot$ is the total localization functor left-adjoint to the
embedding $\sSets^h \subset \cCat^h$. Altogether, we have a
commutative diagram
\begin{equation}\label{cat.set.sq}
\begin{CD}
\Unf(\Sets) @>>> \sSets^h\\
@VVV @VVV\\
\cCat @>>> \cCat^h
\end{CD}
\end{equation}
of enhanced categories and enhanced full embeddings. All the full
embeddings in \eqref{cat.set.sq} are left-admissible; the adjoint
functor $\sSets^h \to \Unf(\Sets)$ is the truncation functor that
sends an enhanced groupoid to the set of its connected components,
and the adjoint functor $\cCat^h \to \cCat$ is the truncation
functor \eqref{trunc.eq}. In principle, for any essentially small
category $\C$, we also have the naive full localization $h^\Tot(\C)$
with respect to the set of all maps, but this is {\em not the same}
as the enhanced localization $\Hh^\Tot(\C)$ -- even if $\C$ is an
unenhanced essentially small category, the groupoid $h^\Tot(\C)$ is
only the truncation \eqref{trunc.eq} of the enhanced groupoid
$\Hh^{\Tot}(\C)$.

As another application of enhanced colimits in $\cCat^h$, one can
construct an enhanced version of the nerve functor
\eqref{nerve.eq}. Namely, for any small enhanced category $\C$,
define an enhanced category $\Delta \bbi^h \C$ by the cartesian
square
\begin{equation}\label{del.enh.sq}
\begin{CD}
\Delta \bbi^h \C @>>> \cCat^h \bbi^h \C\\
@VVV @VVV\\
\Unf(\Delta) @>>> \cCat^h,
\end{CD}
\end{equation}
where the bottom arrow is the standard full embedding $[n] \mapsto
\Unf([n])$, and the arrow on the right is the enhanced family of
groupoids of Proposition~\ref{ccat.h.prop}. Then the arrow on the
left in \eqref{del.enh.sq} is also an enhanced family of groupoids,
thus corresponds to an enhanced functor $\Unf(\Delta)^\iota \to
\sSets^h$ by Proposition~\ref{groth.prop}, and sending $\C$ to
$\Delta \bbi^h \C \to \Unf(\Delta)$ defines an enhanced functor
\begin{equation}\label{enh.nerve.eq}
\cN:\cCat^h \to \Delta^o_h\sSets^h, \qquad \C \mapsto \Delta \bbi^h
\C.
\end{equation}
The enhanced functor \eqref{enh.nerve.eq} is left-reflexive and
fully faithful, and one can also describe its essential image (it
consists of ``complete Segal enhanced families of groupoids'', see
\cite[Definition 7.5.6.1]{big} and \cite[Proposition
  7.5.6.2]{big}). The enhanced {\em contraction functor}
$\Con:\Delta^o_h\sSets^h \to \cCat^h$ left-adjoint to
\eqref{enh.nerve.eq} is constructed in \cite[Lemma 7.5.6.5]{big};
for any small enhanced family of groupoids $\E \to \Unf(\Delta)$,
$\Con(\E)$ fits into an enhanced-cocartesian square
\begin{equation}\label{enh.con.sq}
\begin{CD}
\Unf(\overline{\nu}_\idot)^*\E @>>>
\Hh^{\Tot}(\Unf(\overline{\nu}_\idot)^*\E)\\
@VVV @VVV\\
\Unf(\nu_\idot)^*\E @>>> \Con(\E),
\end{CD}
\end{equation}
where $\Hh^\Tot(-)$ is the total localization functor,
$\nu_\idot:\Delta_\idot \to \Delta$ is the cofibration with fibers
$[n]$ induced from the tautological cofibration $\Cat_\idot \to
\Cat$ via the standard embedding $\Delta \to \Cat$,
$\Delta_{\idot\Tot} \cong [0] \setminus \Delta \subset \Delta_\idot$
is the subcategory with the same object and those maps that are
cocartesian over $\Delta$, and
$\overline{\nu}_\idot:\Delta_{\idot\Tot} \to \Delta$ is the discrete
cofibration induced by $\nu_\idot$. The square \eqref{enh.con.sq} is
a version of the localization square \eqref{enh.loc.sq}, and we in
fact have
\begin{equation}\label{del.loc.eq}
\Con(\E) \cong \Hh^W(\Unf(\nu_\idot)^*\E),
\end{equation}
where $W$ is the class of enhanced maps in $\Unf(\nu_\idot)^*\E$
cocartesian over $\E$. If $\E = \Delta \bbi^h \C$ for a small
enhanced category $\C$, then $\C \cong \Con(\E)$. One can let
$\Delta \bb^h \C = \fFun^h(\Delta_\idot|\Delta,\C)$ be the relative
functor category \eqref{enh.rel.dia}; we then have the embedding
$\Delta \bbi^h \C \to \Delta \bb^h \C$, and the bottom arrow in
\eqref{enh.con.sq} is induced by the evaluation pairing
$\ev:\Unf(\nu_\idot)^*(\Delta \bb^h \C) \to \C$.

As an alternative, one can observe that explicitly, $\Delta_\idot$
is the category of pairs $\langle [n],l \rangle$, $[n] \in \Delta$,
$l \in [n]$, with maps $\langle [n],l \rangle \to \langle [n'],l'
\rangle$ given by maps $f:[n] \to [n']$ such that $f(l) \leq
l'$. Then $\nu_\idot:\Delta_\idot \to \Delta$ is the forgetful
functor $\langle [n],l \rangle \mapsto [n]$, and it has a
right-adjoint $\nu^\dg:\Delta \to \Delta_\idot$, $[n] \mapsto
\langle [n],n \rangle$. The functor $\nu^\dg$ lifts to a
left-admissible enhanced embedding $\nu^\dg:\E \to
\Unf(\nu_\idot)^*\E$, and \eqref{del.loc.eq} induces an equivalence
\begin{equation}\label{loc.sp.del.eq}
\Con(\E) \cong \Hh^+(\E),
\end{equation}
where $+$ is the class of enhanced morphisms in $\E$ whose image in
$\Delta$ is special in the sense of Subsection~\ref{del.subs}. If
$\E = \Delta \bbi^h \C$ for a small enhanced category $\C$, then $\ev$
induces an enhanced functor
\begin{equation}\label{xi.eq}
\xi = \ev \circ \nu^\dg:\Delta \bbi \C \to \C,
\end{equation}
and \eqref{loc.sp.del.eq} becomes an equivalence $\C \cong
\Hh^+(\Delta \bbi^h \C)$, an enhanced version of the identification
$h^+(\Delta \bbi I) \cong I$ of Subsection~\ref{del.subs}. We refer
to \cite[Subsection 7.5.6]{big} for further details.

As a further application, one can give a characterization of the
transpose enhanced cofibrations that does not depend on the
Grothendieck construction of Proposition~\ref{groth.prop}. Namely,
let $\gamma:\C \to \E$ be an enhanced fibration of small enhanced
categories, with the enhanced-transpose cofibration $\C_{h\perp} \to
\E$. Then by \cite[Corollary 7.5.6.7]{big}, we have a functorial
enhanced-cocartesian square
\begin{equation}\label{perp.sq}
\begin{CD}
\Unf(\iota \circ \overline{\nu}_\idot)^*\xi^*\C @>>>
\C_{hv}\\
@VVV @VVV\\
\Unf(\iota \circ \nu_\idot)^*\xi^*\C @>>> \C_{h\perp},
\end{CD}
\end{equation}
where $\iota:\Delta \to \Delta$ is the involution $[n] \mapsto
[n]^o$, $\xi:\Delta \bbi^h \E \to \E$ is the enhanced functor
\eqref{xi.eq}, and $\C_{hv}$ is the enhanced category
\eqref{h.v.eq}, for the class $v$ of enhanced morpfisms $f$ in $\C$
such that $\gamma(f)$ is invertible.

\subsection{Large categories.}

The theory of large enhanced categories is even more deficient than
its unenhanced prototype --- while in the unenhanced setting, it is
the functor categories that present a problem, in the enhanced
situation, even semicartesian products may not exist. A possible
solution to this difficulty is an enhanced version of the machinery
of accessible categories described in
Subsection~\ref{acc.subs}. This is the content of \cite[Section
  7.6]{big}, and the story is again largely parallel to what happens
in the unenhanced setting.

\begin{defn}
For any regular cardinal $\kappa$, an enhanced category $\E$ is {\em
  $\kappa$-filtered} if any enhanced functor $\C \to \E$ from a
$\kappa$-bounded small enhanced category $\C$ admits an enhanced
cone $\C^{h>} \to \E$. An enhanced category $\E$ is {\em
  $\kappa$-filtered-cocomplete} if any enhanced functor $\C \to \E$
from a $\kappa$-filtered small enhanced category $\C$ has an
enhanced colimit. An enhanced object $c \in \C_\ppt$ in a
$\kappa$-filtered-cocomplete enhanced category $\C$ is {\em
  $\kappa$-compact} if the enhanced functor \eqref{yo.c.eq}
preserved $\kappa$-filtered enhaned colimits, and $\cComp_\kappa(\C)
\subset \C$ is the enhanced full subcategory spanned by
$\kappa$-compact enhanced objects.
\end{defn}

For any regular cardinal $\kappa$ and enhanced category $\E$, the
{\em enhanced $\kappa$-inductive completion} is
the envelope $\iInd_\kappa(\E) = \Env(\E,\F_\kappa)$ of
Proposition~\ref{env.prop} with respect to the full subcategory
$\F_\kappa \subset \cCat^h$ spanned by $\kappa$-filtered small
enhanced categories. By definition, it comes equipped with a full
embedding $\E \to \iInd_\kappa(\E)$, and for any
$\kappa$-filtered-cocomplete enhanced category $\C$, the full
embedding $\cComp_\kappa(\C) \to \C$ canonically extends to an
enhanced functor
\begin{equation}\label{iind.comp.eq}
\iInd_\kappa(\cComp_\kappa(\C)) \to \C
\end{equation}
that preserves $\kappa$-filtered enhanced colimits. This is an
enhanced version of \eqref{ind.comp.eq}, and while we do not have
\eqref{ind.eq} in the enhanced case, it is still true that the
enhanced functor \eqref{iind.comp.eq} is fully faithful (\cite[Lemma
  7.5.7.21]{big}).

\begin{defn}\label{enh.acc.def}
For any regular cardinal $\kappa$, an enhanced category $\C$ is {\em
  $\kappa$-accessible} if it is $\kappa$-filtered-cocomplete,
$\cComp_\kappa(\C)$ is small, and \eqref{iind.comp.eq} is essentally
surjective. A $\kappa$-accessible category is {\em
  $\kappa$-presentable} if it is enhanced-cocomplete. An enhanced
functor $\C \to \C'$ between $\kappa$-accessible enhanced categories
is {\em $\kappa$-accessible} if it preserves $\kappa$-fitered
enhenaced colimits. An enhanced category $\C$ is {\em accessible}
resp.\ {\em presentable} if it is $\kappa$-accessible
resp.\ $\kappa$-presentable for some regular cardinal $\kappa$, and
an enhanced functor is {\em accessible} if it $\kappa$-accessible
for some $\kappa$.
\end{defn}

As in the unehnaced case, a $\kappa$-presentable enhanced category
$\C$ is $\kappa'$-presentable for any regular cardinal $\kappa' >
\kappa$, and there is an order relation $\trr$ on regular cardinals,
\cite[Definition 2.1.6.4]{big}, such that a $\kappa$-accessible
enhanced category $\C$ is $\mu$-accessible for any regular cardinal
$\mu \trr \kappa$. Moreover, for any set of accessible enhanced
categories and enhanced functors, one can choose a single cardinal
$\kappa$ such that they are all $\kappa$-accessible. There is a
bunch of structural theorems for accessible and presentable enhanced
categories that repeat the standard results in the unenhanced
setting, see \cite[Subsection 7.6.5]{big}. In particular, by
\cite[Corollary 7.6.5.8]{big}, $\cCat^h$ and $\sSets^h$ are
presentable, by \cite[Lemma 7.6.5.7]{big}, a $\kappa$-accessible
enhanced category is $\kappa$-presentable as soon as it is
$\kappa$-cocomplete, and it is automatically enhanced-complete,
while by \cite[Corollary 7.6.3.7]{big}, a small enhanced category is
accessible iff it is enhanced-Karoubi-closed. Then there is the
following enhanced version of Subsection~\ref{acc.subs}~\thetag{i}
that also serves as an accessible version of
Corollary~\ref{semi.corr}.

\begin{prop}\label{acc.semi.prop}
Assume given accessible enhanced categories $\C$, $\C_0$, $\C_1$ and
accessible enhanced functor $\gamma_l:\C_l \to \C$, $l=0,1$. Then
there exists a semicartesian square \eqref{cat.sq} of accessible
enhanced categories and accessible enhanced functors. Moreover,
assume given another commutative square \eqref{cat.1.sq} of
accessible enhanced categories and accessible enhanced functor.
Then there exists an accessible enhanced functor $\phi:\C'_{01} \to
\C_{01}$ and enhanced isomorphisms $a_l:\gamma^l_{01} \circ \phi
\cong {\gamma'}^l_{01}$, $l = 0,1$, and the triple $\langle
\phi,a_0,a_1 \rangle$ is unique up to an enhanced isomorphism.
\end{prop}

\proof{} This is \cite[Proposition 7.6.5.9]{big}. \endproof

For Subsection~\ref{acc.subs}~\thetag{ii}, we note that for any two
$\kappa$-accessible categories $\C$, $\E$, $\kappa$-accessible
enhanced functors $\C \to \E$ form an enhanced category
$\fFun^h_\kappa(\C,\E) = \fFun^h(\cComp_\kappa(\C),\E)$, and all
acceesible enhanced functors form an enhanced category
\begin{equation}\label{enh.acc.fun}
\fFun^h(\C,\E) = \bigcup_\kappa\fFun^h_\kappa(\C,\E),
\end{equation}
where the union is over all regular cardinals $\kappa$ such that
both $\C$ and $\E$ are $\kappa$-accessible. Each indivual enhanced
category in the right-hand side of \eqref{enh.acc.fun} is accessible
(\cite[Corollary 7.6.5.9]{big}), but unfortunately, the whole
category $\fFun^h(\C,\E)$ need not be accessible.

Generalizing Definition~\ref{enh.acc.def} even further, one can say
that an enhanced category $\C$ is {\em tame} if its enhanced Karoubi
completion is accessible, and an enhanced functor between tame
enhanced categories is tame iff its canonical extension to enhanced
Karoubi completions is accessible. Then any small enhanced category
is tame, but Proposition~\ref{acc.semi.prop} breaks down: one can
always construct a semicartesian square \eqref{cat.sq}, but
$\C_{01}$ and $\gamma^l_{01}$, $l=0,1$ need not be tame. The problem
does not occur if either of the tame enhanced functors $\gamma_0$,
$\gamma_1$ is an enhanced fibration or cofibration.

\begin{prop}\label{tame.prop}
For any tame enhanced functor $\pi:\C \to \E$ between tame enhanced
categories that is an enhanced fibration or cofibration, so is the
enhanced Karoubi closure $\Env(\pi,\P)$. For any semicartesian
square \eqref{cat.sq} such that $\C$, $\C_0$, $\C_1$, $\gamma_l$,
$l=0,1$ are tame, and $\gamma_0$ is an enhanced fibration
resp.\ cofibration, $\C_{01}$ and $\gamma^l_{01}$, $l=0,1$ are tame,
$\gamma^0_{01}$ is an enhanced fibration resp.\ cofibration, and the
square has the same universal property with respect to tame
commutative square as in Proposition~\ref{acc.semi.prop}.
\end{prop}

\proof{} This is \cite[Lemma 7.6.6.2]{big}. \endproof

As in the unenhanced setting of Subsection~\ref{acc.subs}, the main
problem with the theory of acceesible enhanced categories is that
the enhanced-opposite $\C^\iota$ to an accessible enhanced category
$\C$ is usually not accessible. In particular, statements about
accessible and tame enhanced fibrations and cofibrations such as
Proposition~\ref{tame.prop} do not imply each other by passing to
the enhanced opposite categories, and need separate
proofs. Moreover, there is no meaningful Grothendieck
construction. One of its corollaries that survives is the
construction of the transpose-opposite enhanced cofibration
$\C_{h\perp}$ for an enhanced fibration $\C \to \E$, and dually for
an enhanced cofibration. Namely, for any tame enhanced fibration $\C
\to \E$, with small $\E$, it is proved in \cite[Proposition
  7.6.6.6]{big} that there exists a tame enhanced cofibration
$\C_{h\perp} \to \E^\iota$ that fits into a commutative square
\eqref{perp.sq} that is cocartesian with respect to tame commutative
squares, and for any tame enhanced cofibration $\C \to \E$, again
with small $\E$, there is a dual construction and characterization
of the transpose tame enhanced fibration $\C^{h\perp} \to
\E^\iota$. The cocartesian property insures that $\C_{h\perp}$
resp.\ $\C^{h\perp}$ is unique; moreover, for any tame enhanced
fibrations $\C,\C' \to \E$, an enhanced functor $\gamma:\C \to \C'$
cartesian over $\E$ induces an enhanced functor
$\gamma_{h\perp}:\C_{h\perp} \to \C'_{h\perp}$ cocartesian over
$\E^\iota$, and similarly for tame enhanced cofibrations.

Finally, one can consider the enhanced fibration $\cCat^h \bb^h \E
\to \cCat^h$ of \eqref{ccat.h.bb}. Here the result --- namely,
\cite[Proposition 7.6.6.8]{big} --- is as follows. If $\E$ is
$\kappa$-presentable for some regular cardinal $\kappa$, then
$\cCat^h \bb \E$ is also $\kappa$-presentable, and \eqref{ccat.h.bb}
is preserves enhanced colimits. If $\E$ is only accessible, then in
general, we cannot say that $\cCat^h \bb \E$ is accessible, nor even
tame. However, for any small enhanced full subcategory $\I \subset
\cCat^h$, the enhanced category $\I \bb^h \E = \I \times_{\cCat^h}
(\cCat^h \bb^h \E)$ is accessible, and so is the enhanced fibration
$\I \bb^h \E \to \I$. While this result is rather limited, it still
allows for some applications. In particular, for any enhanced
cofibration $\E' \to \E$ between small enhanced categories, and any
accessible enhanced category $\C$, we can define the accessible
relative functor category $\fFun^h(\E'|\E,\C)$ by
\eqref{enh.rel.dia}, and it has the same universal property with
respect to tame enhanced functors as in the case when $\C$ is
small. Analogously, for any enhanced fibration $\E' \to \E$ between
small enhanced categories, one can achieve the same result by
setting $\fFun^h(\E'|\E,\C) =
\fFun^h(\E'_{h\perp}|\E^\iota,\C)_{h\perp}$, as in
\eqref{rel.fun.perp}.

\section{How it works.}\label{app.sec}

While our stated purpose in this paper is just to give a toolkit for
working with enhanced categories, with all the constructions and
statements, but precise references to \cite{big} instead of proofs,
it might be worthwhile to at least mention the main ideas behind
those proofs. This section is but an appendix; an uninterested
reader is adviced to just skip it.

\subsection{Brown representability.}

Our starting point is the classic Brown Representability Theorem
that characterizes representable functors from the homotopy category
$h^W(\Delta^o\Sets_\idot)$ of pointed sets to $\Sets$ (\cite{brown},
\cite[Chapter 9]{switz}). A functor $Y:(\Delta^o\Sets_\idot)^o \to
\Sets$ defines a representable functor $h^W(\Delta^o\Sets_\idot)^o
\to \Sets$ iff it is homotopy-invariant (that is, inverts weak
equivalences, thus descends to $h^W(\Delta^o\Sets_\idot)^o$),
additive (that is, sends arbitrary coproducts to products) and
semiexact: for any $X \in \Delta^o\Sets_\idot$ equipped with subsets
$X_0,X_1 \subset X$, with intersection $X_{01} = X_0 \cap X_1$, the
map $Y(X) \to Y(X_0) \times_{Y(X_{01})} Y(X_1)$ is surjective. The
version for unpointed sets is wrong, essentially because
$h^W(\Delta^o\Sets)$ is not Karoubi-closed (there is a mistake to
this effect in \cite{switz} which is corrected in the Russian
translation). Our first observation is that the theorem holds for
unpointed sets, provided one replaces functors $(\Delta^o\Sets)^o
\to \Sets$ --- that is, discrete small fibrations over
$\Delta^o\Sets$ --- with small families of groupoids over
$\Delta^o\Sets$.

The fact that passing to groupoids helps one to get rid of
basepoints is not surprising at all -- actually, it is hard to
resist quoting \cite{prog}:
\begin{itemize}
\item {\em Ceci est li\'e notamment au fait que les gens s’obstinent
  encore, en calculant avec des groupes fondamentaux, \`a fixer un
  seul point base, plut\v{o}t que d’en choisir astucieusement
  tout un paquet qui soit invariant par les sym\'etries de la
  situation, lesquelles sont donc perdues en route.}
\end{itemize}
However, the more important advantage of passing to groupoids is
that while a map of sets can only be surjective, a functor between
groupoids can be either just essentially surjective, or also full,
thus an epivalence. It is the latter stronger condition that goes
into our notion of semiexactness.

For reasons of backward compatibility, \cite[Section 5.3.1]{big}
provides a proof of Brown Representability for families groupoids
along the lines of \cite[Chapter 9]{switz}, where the representing
object is constructed by gluing cells. Then we give an alternative
proof that is somewhat dual: the representing object is obtained as
the limit of a Postnikov tower. More generally, in \cite[Section
  3.3]{big}, we introduce a class of Reedy categories that we call
{\em cellular}; these include $\Delta$ (but not $\Delta^o$),
left-bounded partially ordered sets, products $I \times I'$ of two
cellular Reedy categories $I$, $I'$, and the category of elements
$IX$ of a functor $X:I^o \to \Sets$ for a cellular Reedy category
$I$. Then for any cellular Reedy category $I$, we introduce the
notion of a {\em liftable object} $c \in \C_X \subset \C$ in a
family of groupoids $\C \to I^o\Sets$ described in terms of
extensions with respect to injective maps $X \to X'$, and we prove
that if the family is small, additive and semiexact, a liftable
object exists. If $I = \Delta$ and $\C$ is constant along weak
equivalences, then a standard argument shows that a liftable object
represents $\C$, so this proves the representability
theorem. Moreover, for any cellular $I$, the product $I \times
\Delta$ is cellular, so we also obtain representability for families
over $I^o\Delta^o\Sets$, with the Reedy model structure.

A useful additional feature of cellular Reedy categories is that any
functor $X:I^o \to \Sets$ admits a functorial increasing {\em
  filtration by skeleta} $\sk_nX \subset X$, $n \geq 0$ --- for $I =
\Delta$, $\sk_n X \subset X$ is the union of non-degenerate
simplices of dimension $\leq n$. This is used in the construction of
a liftable object, but it gives more: it turns out that in order to
find this liftable object, it suffices to restrict the small
additive semiexact family $\C$ to the full subcategory $I^o_f\Sets
\subset I^o\Sets$ of functors $X$ such that $\sk_nX = X$ for some
$n$ (if $I=\Delta$, these correspond to finite-dimensional
simplicial sets). Any such family over $I^o_f\Sets$ then extends to
a family over the whole $I^o\Sets$, uniquely up to an equivalence
unique up to a non-unique isomorphism. This is the origin of
Proposition~\ref{pl.prop}.

\subsection{Enhanced groupoids.}

In our enhancement formalism, objects in $h^W(\Delta^o\Sets)$
correspond to small enhanced groupoids, and these are given by small
families of groupoids $\C \to \Posf$ over the category of
finite-dimensional partially ordered sets. The category $\Posf$ does
not have a model structure, but it has two of the three ingredients
of one, namely, classes $C$ and $W$ (but not $F$). The class $C$ is
formed by left-closed embeddings, and the class $W$ is the class of
{\em reflexive maps}, defined as finite compositions of
left-reflexive and right-reflexive maps. We have standard pushout
squares \eqref{st.sq}, and if the top arrow in \eqref{st.sq} is in
$W$, then so is the bottom arrow. The resulting structure is
axiomatized in \cite[Section 5.1]{big} under the name of a {\em
  CW-category} (there is also an additional axiom that is a cut-down
version of the model category factorizations). While the notion of a
CW-category is not particularly deep --- it is even more of a purely
technical gadget than the notion of a model category --- it allows
one to prove something. In particular, there is the following
result. A CW-category structure on some category $I$ already allows
one to define additive semiexact families of groupoids $\C \to I$
constant along $W$. Assume given such a family $\C \to \Posf$, and
say that it is {\em stably constant} along a map $f$ in $\Posf$ if
it is constant along $\id_J \times f$ for any finite $J \in
\Posf$. Then \cite[Lemma 5.2.2.7]{big} implies that
\begin{enumerate}
\renewcommand{\labelenumi}{(\Roman{enumi})}
\stepcounter{enumi}
\item for any additive semiexact family of groupoids $\C \to \Posf$
  constant along $W$, and a square \eqref{st.sq} in $\Posf$ such
  that $\C$ is stably constant along the top arrow, $\C$ is also
  stably constant along the bottom arrow.
\end{enumerate}
This is useful because it shows that while the class $W$ of
reflexive maps in $\Posf$ is really small, an additive semiexact
family $\C$ constant along $W$ is automatically constant along a
much bigger class of maps. Namely, say that a class $v$ of maps in a
category $\C$ is {\em saturated} if it is closed under retracts and
has the two-out-of-three property: for any composable pair $f$, $g$
of maps in $\C$, if two out of the three maps $f$, $g$, $f \circ g$
are in $W$, then so is the third.

\begin{defn}\label{ano.def}
The class of {\em anodyne maps} is the smallest saturated class of
maps in $\Posf$ that contains all the projections $J \times [1] \to
J$, $J \in \Posf$, and has the following property:
\begin{enumerate}
\renewcommand{\labelenumi}{(\Roman{enumi})}
\stepcounter{enumi}
\stepcounter{enumi}
\item for any cofibrations $J',J' \to J$ in $\Posf$, a map $f:J' \to
  J''$ cocartesian over $J$ is anodyne if so are the fibers
  $f_j:J'_j \to J''_j$ for all $j \in J$.
\end{enumerate}
\end{defn}

The class of anodyne maps contains all reflexive maps --- in fact,
\thetag{III} is not even needed for this, see \cite[Lemma
  2.1.3.11]{big} -- but it is strictly bigger. For example, for any
excision square \eqrefi{J.S.sq}{i} in $\Posf$, the bottom arrow is
reflexive, while the top one is usually only anodyne. It turns out
that this is the universal example: by \cite[Proposition
  2.1.9.7]{big}, the smallest saturated class of maps in $\Posf$ that
contains the top arrows for all excision squares  \eqrefi{J.S.sq}{i}
and is closed under standard pushouts contains all anodyne
maps. Therefore any additive semiexact family of groupoids $\C \to
\Posf$ that is constant along the projections $J \times [1] \to J$
and satisfies excision is constant along all anodyne maps.

To go from families of groupoids over $\Delta^o_f\Sets$ to families
over $\Posf$, one can use the nerve functor $N:\Posf \to
\Delta^o_f\Sets$; it is easy to check that it sends all anodyne maps
to weak equivalences, and it also preserves coproducts and pushout
squares \eqref{st.sq}, so that for any representable small family of
groupoids $\C \to \Delta^o_f\Sets$, the family $N^*\C$ is a small
enhanced groupoid. To go in the other directions, one uses another
property of cellular Reedy categories: by \cite[Proposition
  3.3.3.10]{big}, for any finite-dimensional cellular Reedy category
$I$, there exists a finite-dimensional partially ordered set $Q(I)$
and a functor $Q(I) \to I$ that has anodyne right comma-fibers $i
\setminus Q(I)$. Moreover, this construction is sufficiently
functorial to be made into a functor $Q:\Delta_f^o\Sets \to \Posf$,
$X \mapsto Q(\Delta X)$, and we have a functorial pointwise-anodyne
map $Q \circ N \to \id$ and a functorial map $N \circ Q \to \id$
that is a pointwise weak equivalence. This identifies appropriate
families of groupoids over $\Posf$ and $\Delta^o_f\Sets$, and proves
Theorem~\ref{repr.thm} for small enhanced groupoids (the precise
statement is \cite[Proposition 6.1.2.3]{big}).

Moreover, one can extend the correspondence between $\Posf$ and
$\Delta^o_f\Sets$ to a correspondence between $\Posf^+$ and the
whole $\Delta^o\Sets$. To do this, one defines the class of {\em
  $+$-anodyne maps} as the smallest saturated class of maps in
$\Posf^+$ that again contains all the projections $J \times [1] \to
J$ and has the property \thetag{III}. An example of a map that
is $+$-anodyne but not anodyne is the top arrow in a square
\eqref{J.N.sq}, and as before, one shows that this new example is
universal: by \cite[Proposition 2.1.10.6]{big}, the class of
$+$-anodyne maps is the smallest saturated class of maps in
$\Posf^+$ that is closed under coproducts and standard pushouts, and
contains the projections $J \times [1] \to J$ together with top
arrows of squares \eqrefi{J.S.sq}{i} and \eqref{J.N.sq} in
$\Posf^+$. With a little additional effort, this establishes a
bijection between small enhanced groupoids $\C \to \Posf^+$ and
representable families over $\Delta^o\Sets$.

\subsection{Enhanced categories.}

To go from small enhanced groupoids to small enhanced categories,
one has to replace simplicial sets $X \in \Delta^o\Sets$ with
complete Segal spaces $X \in \Delta^o\Delta^o\Sets$, and we need an
appropriate $\Pos$-based model for
those. Proposition~\ref{enh.loc.prop} suggests looking at partially
ordered sets $J$ equipped with a set $W$ of maps $[1] \to J$, in the
spirit of relative categories of \cite{BK}, but it turns out that
there is a better choice.

\begin{defn}
A {\em biordered set} is a partially ordered set $J$ equipped with
two additional orders $\leq^l$, $\leq^r$ such that $j \leq^l j'$, $j
\leq^r j'$ imply $j \leq j'$, and for any $j \leq j'$, there exists
a unique $j'' \in J$ such that $j \leq^l j'' \leq^r j'$. A {\em
  biordered map} is a map between biordered sets that preserves the
orders $\leq$ and $\leq^l$, and a biordered map is {\em strict} if
it also preserves $\leq^r$.
\end{defn}

Categorically, a biorder on a partially ordered set $J$ is a
factorization system in the sense of \cite{bou2}, and it behaves
accordingly. In particular, any of the orders $\leq^l$, $\leq^r$
uniquely determines the other one, if it exists, and for any
cofibration $\pi:J' \to J$, we have a natural biorder on $J'$ such
that $j \leq^r j'$ iff $j \leq j'$ and $\pi(j)=\pi(j')$ (and then $j
\leq^l j'$ iff $j \leq j'$, and the corresponding map in $J'$ is
cocartesian over $J$). Moreover, by a stroke of luck, we have a
reasonable notion of a cofibration: a biordered map $J' \to J$ is a
{\em bicofibration} if it is strict, and the underlying map is a
cofibration. We denote the category of biordered sets by $\Rel$, and
we let $R,L:\Pos \to \Rel$ be the functors sending $J \in \Pos$ to
itself with discrete biorder $\leq^r$ resp.\ $\leq^l$. We have the
embedding $\Delta \times \Delta \to \Rel$, $[n] \times [m] \mapsto
L([n]) \times R([m])$, and this gives rise to the biordered nerve
functor
\begin{equation}\label{rel.nerve.eq}
\begin{aligned}
&\Nn_\dm:\Rel \to \Delta^o\Delta^o\Sets,\\
&\Nn_\dm(J)([n] \times [m]) = \Hom_{\Rel}(L([n]) \times R([m]),J).
\end{aligned}
\end{equation}
We then let $\Relf \subset \Relf^+ \subset \Rel$ be the full
subcategories spanned by finite-dimensional resp.\ left-bounded
biordered sets, and we define the classes of {\em bianodyne}
resp.\ {\em $+$-bianodyne} maps by repeating
Definition~\ref{ano.def}, with the projections $J \times [1] \to J$
replaced by projections $J \times R([1]) \to J$, and cofibrations in
\thetag{III} replaced by bicofibrations. The standard anodyne
resp.\ $+$-anodyne maps in \eqrefi{J.S.sq}{i} resp.\ \eqref{J.N.sq}
carry natural biorders that makes them bianodyne
resp.\ $+$-bianodyne, and there is an additional class of bianodyne
maps that appear as top arrows in squares \eqref{J.Z.sq}. Then one
proves (\cite[Proposition 2.2.7.3]{big}) that these generate all
bianodyne resp.\ $+$-bianodyne maps by the same procedure of taking
saturation and standard pushouts, and that for any small
family of groupoids $\C \to \Delta^o\Delta^o\Sets$ represented by a
Segal space $X$, the pullback $\Nn_\dm^*\C$ is an additive semiexact
family constant along $+$-bianodyne maps.

To go back, one again takes the anodyne resolution $Q(X) \to
\Delta\Delta X$ of the cellular Reedy category $\Delta\Delta X$
associated to a bisimplicial set $X$, and notes that it comes
equipped with the projection $Q(X) \to \Delta$ to the first factor
$\Delta$. One then considers the fibered product $Q(X)_\idot =
\Delta_\idot \times_\Delta Q(X)$, where $\Delta_\idot \to \Delta$ is
the cofibration with fibers $[n]$ induced by the cofibration
$\Cat_\idot \to \Cat$, and defines $\QQ_\dm(X)$ as $Q(X)_\idot$ with
the biorder induced by the cofibration $Q(X)_\idot \to Q(X)$. This
gives a functor $\QQ_\dm:\Delta^o\Delta^o\Sets \to \Relf^+$, and one
proves that $\Nn_\dm^*$ and $\QQ_\dm^*$ establish a bijection
between small families of groupoids $\C \to \Delta^o\Delta^o\Sets$
represented by Segal spaces, and small families of groupoids $\C \to
\Relf^+$ that are additive, semiexact, and satisfy the appropriate
biordered versions of semicontinuity, excision and the cylinder
axiom. The precise statement is \cite[Proposition 6.2.4.1]{big}.

To deduce Theorem~\ref{repr.thm} in full generality, it remains to
do two things: firstly, restrict our families of groupoids over
$\Relf^+$ to $\Posf^+$, and secondly, extend them to families of
categories. The first is done by pullback with respect to the
functor $L:\Posf^+ \to \Relf^+$, and there is a general
reconstruction result, \cite[Proposition 6.3.2.7]{big} that shows
that under some additional assumptions (called ``coherence'' and
``reflexivity''), a family of groupoids $\C \to \Relf^+$ obtained
from a Segal space is uniquely determined by its restriction $L^*\C
\to \Posf^+$. On the Segal space side, the additional assumptions
roughly correspond to requiring that $X$ is complete. Finally, to go
from families of groupoids to families of categories, one uses a
completely general reconstruction result, \cite[Proposition
  7.1.5.2]{big} that shows that a separated non-degenerate reflexive
family of categories $\C \to \Posf^+$ is uniquely and functorialy
determined by the family of groupoids $\C_\flat \to \Posf^+$; on the
level of functors, this is Lemma~\ref{fl.le}

\subsection{Enhanced category theory.}

In actual reality, the proof of Theorem~\ref{repr.thm} sketched
above requires quite a lot of space; at the end of the day, it takes
up the whole of \cite[Chapter 5]{big} and \cite[Chapter 6]{big},
with the additional preliminary theory of \cite[Chapter
  3]{big}. However, once it has been done, the rest of the story
becomes quite straightforward. In particular,
Proposition~\ref{enh.loc.prop} is almost immediate. The actual
statement is \cite[Proposition 7.3.3.4]{big}, and it is slightly
more precise. Effectively, one first extends an enhanced category
$\C \to \Posf^+$ to a family $\C^\dm \to \Relf^+$ called the {\em
  unfolding} of $\C$, and given by the cartesian square
\begin{equation}\label{unf.sq}
\begin{CD}
\C^\dm @>>> U^*\C\\
@VVV @VVV\\
\Unf^\dm(\C_\ppt) @>>> U^*\Unf(\C_\ppt),
\end{CD}
\end{equation}
where $U:\Relf^+ \to \Posf^+$ sends a biordered set to its
underlying partially ordered set, the arrow on the right is the
truncation functor \eqref{trunc.eq}, and $\Unf^\dm(\C_\ppt) \subset
U^*\Unf(\C_\ppt)$ is the full subcategory of functors $U(J)^o \to
\C_\ppt$ that invert all maps corresponding to order relations $j
\leq^l j'$. Then one defines a {\em universal object} for $\C$ as an
object $c \in \C^\dm_J$ such that for any enhanced category $\E$ and
object $e \in \E^\dm_J$, there exists a unique enhanced functor
$\gamma:\C \to \E$ such that $\gamma^\dm:\C^\dm \to \E^\dm$ sends
$c$ to $e$, and for any map $f:e \to e'$, there exists a unique
enhanced map $\gamma \to \gamma'$ between the corresponding enhanced
functors $\gamma,\gamma':\C \to \E$ that evaluates to $f$ at
$c$. One then proves that for any small enhanced category $\C$, a
universal object exists, and this gives
Proposition~\ref{enh.loc.prop}, with $W$ given by all biordered maps
$R([1]) \to J$. Corollary~\ref{fun.h.corr} is rather easy, and
Corollary~\ref{semi.corr}~\thetag{ii} -- namely, the functoriality
of semicartesian products -- is immediate. For
Corollary~\ref{semi.corr}~\thetag{i}, one uses
Theorem~\ref{repr.thm} -- the product $\C_0 \times^h_\C \C_1$ is
represented by the homotopy fibered product of the representing
objects for $\C$, $\C_0$ and $\C_1$.

As far as examples of Subsection~\ref{exa.subs} are concerned, then
in principle, Proposition~\ref{ccat.prop},
Proposition~\ref{mod.prop} and Proposition~\ref{ho.prop} can be
proved directly, by checking all the axioms. However, there is a
general trick that allows one not to do that. Namely, one defines a
{\em complete family of categories} $\C \to \Posf^+$ as a family of
categories whose transition functors $f^*$ for all maps $f$ in
$\Posf^+$ admit right-adjoint functors $f_*$ such that the base
change maps for diagrams \eqref{bc.dia} are isomorphisms
(\cite[Definition 7.1.4.1]{big}). Then one shows (\cite[Section
  7.1.4]{big}) that for a complete family, two rather mild
conditions called {\em tightness} and {\em weak semicontinuity}
insure that the family is reflexive and defines an enhanced category
(and moreover, a complete one). This makes the proofs of
Proposition~\ref{mod.prop} and Proposition~\ref{ho.prop} easy, and
also greatly simplifies the proof of the $\cCat^h$ part of
Proposition~\ref{ccat.h.prop} (in fact, while we have moved
Proposition~\ref{ccat.prop} to a separate statement, its actual
proof in \cite{big} is simply a part of
Proposition~\ref{ccat.h.prop}). For the other enhanced categories of
Proposition~\ref{ccat.h.prop}, the proofs are not hard either, but
they require a stronger form of Proposition~\ref{enh.loc.prop}: one
needs to check that if a small enhanced category $\C$ is equipped
with an $I$-augmentation $\C \to \Unf(I)$, then one can choose a
universal object $c \in \C^{\iota\dm}_J$ for the opposite enhanced
category $\C^\iota$ in such a way that the map $U(J) \to I$ is a
fibration, and $f^*c$ is universal for $f^*\C^\iota$ for any map
$f:I' \to I$. This is what \cite[Proposition 7.3.3.4]{big} actually
says, and the proof is rather delicate: it involves proving a
version of the representabilty theorem for $I$-augmented enhanced
categories (this is \cite[Section 6.4]{big}).

The rest of the story is straightforward. The last place where the
actual construction of the universal objects is used is the
characterization of adjunctions of
Lemma~\ref{adj.cyl.le}. Everything else in Section~\ref{enh.cat.sec}
is done by standard categorical arguments, with the necessary
functoriality provided by Corollary~\ref{semi.corr}, and one does
not need to go into the gory details of the anatomy of representing
objects at all.

\subsection{Towards axiomatization.}

To finish the paper, let us discuss the {\em uniqueness} of the
enhancement technology we have described.

Typically, uniqueness statements for various notions of an
enhancement have the following form: assuming that the answer has
the form X, and satisfies some natural conditions, it must be X. In
particular, for complete Segal spaces, this was very convincingly
done in \cite{toen-inf}. It is probably possible to prove a
statement like this in our approach, too, but in fact, we can do
better. Namely, while it seems hopeless to try to axiomatize the
category $\Sets^h$ of homotopy types, as simply a category, it is
possible to do the same for $\Cat^h$. An ideal statement here would
be something like Lawvere's description of Grothendieck's toposes in
terms of the subobject classifier, see e.g.\ \cite{J}, where a
couple of simple axioms generates an amazingly rich theory. At the
moment, such an ``enhanced topos theory'' describing $\Cat^h$ with
similar elegance is just a dream; but as a proof of concept, let us
show that if one forgets elegance and simplicity, it can be done. At
the end of the day, the only data needed to identify a category $\C$
with $\Cat^h$ is an object $[1] \in \C$, two maps $s,t:\ppt \to
[1]$, and an autoequivalence $\iota:\C \to \C$; these should satisfy
a lot of conditions. Let us give a list.

For simplicity, let us index all enhanced categories by the whole
category $\Pos$, by virtue of the canonical extension of
Proposition~\ref{lf.prop}. Assume given a category $\Cat^?$ that has
all products and coproducts, and is cartesian-closed. Let $\ppt \in
\Cat^?$ be the terminal object. Assume fixed an object $[1] \in
\Cat^?$ equipped with two maps $s,t:\ppt \to [1]$, and let $e:[1]
\to \ppt$ be the tautological map. Assume given an autoequivalence
$\iota:\Cat^? \to \Cat^?$ and an isomorphism $\iota([1]) \cong [1]$
that intertwines $s$ and $t$ (for $? = h$, $\iota(\C) =
\C^\iota$). Let $\Ar^?:\Cat^? \to \Cat^?$ be the internal $\Hom$
functor $\Hhom([1],-)$ (if $?=h$, this is the arrow category functor
$\Ar^h$). Let $\Sets \to \Cat^?$ be the functor sending a set $S$ to
the coproduct of copies of $\ppt$ numbered by the elements $s \in S$,
and assume that the functor is fully faithful, so we have a full
embedding $\Sets \subset \Cat^?$. Let $\Sets^?  \subset \Cat^?$ be
the full subcategory of objects $X$ such that $e^*:X \to \Ar^?(X)$
is an isomorphism, and assume that $\Sets \subset \Sets^?$. Assume
that $\Sets^? \subset \Cat^?$ is right-admissible, and let $R:\Cat^?
\to \Sets^?$ be the right-adjoint to the embedding functor. Say that
an object $X \in \Cat^?$ is {\em rigid} if $R(X) \in \Sets$ (if
$?=h$, then $\Sets^h$ is the full subcategory of enhanced groupoids,
$R$ sends an enhanced category to its isomorphism groupoid, and an
enhnaced category is rigid iff it is of the form $\Unf(I)$ for a
rigid small category $I$). Let $\Pos^?  \subset \Cat^?$ be the full
subcategory spanned by rigid objects $X$ such that $s^* \times
t^*:\Ar^?(X) \to X \times X$ is a monomorphism, and assume that for
any $X \in \Pos^?$, the subset $\Hom([1],X) \subset \Hom(\ppt,X)
\times \Hom(\ppt,X)$ defines a partial order on the set
$\Hom(\ppt,X)$, and that the resulting functor $\Pos^? \to \Pos$ is
an equivalence, so that we have a full embedding $\Pos \subset
\Cat^?$.

Next, assume that for any map $f:X \to Y$ between rigid objects in
$\Cat^?$, and any map $Z \to Y$, there exists a fibered product
$f^*Z = X \times_Y Z \in \Cat^?$, so that we have a pullback functor
$f^*:\Cat^?/Y \to \Cat^?/X$ (if $?=h$, this holds by
\eqref{U.prod}). Note that since $\Cat^?$ is assumed to be
cartesian-closed, $f^*$ automatically has a right-adjoint
$f_*:\Cat^?/X \to \Cat^?/Y$. Let $f^{\iota*} = \iota \circ f^* \circ
\iota$ and $f^\iota_* = \iota \circ f_* \circ \iota$. For any $I \in
\Pos \subset \Cat^?$, consider the map $\tau=t^*:\Ar^?(I) \to I$,
and let $\Sets^?(I) \subset \Cat^?/I$ be the full subcategory
spanned by maps $Z \to I$ such that $\Ar^?(Z) \to \tau^*Z$ is an
isomorphism. Let $\sSets^?  \subset \Pos \setminus \Cat^?$ be the
full subcategory spanned by all $\Sets^?(I)$, and note that the
induced functor $\sigma:\sSets^? \to \Pos$ is a fibration. Assume
that it is an enhanced category (if $?=h$, this is
$\sSets^h$). Moreover, define a category $\sSets^?_\idot$ by the
cartesian square
\begin{equation}\label{sets.q.sq}
\begin{CD}
\sSets^?_\idot @>{c}>> \Ar^h(\sSets^?)\\
@VVV @VV{\sigma}V\\
\Pos @>{\eta}>> \sSets^?,
\end{CD}
\end{equation}
where the bottom arrow sends $I \in \Pos$ to $\id:I \to I$, and
assume that $\sSets^?_\idot$ with the functor $\tau \circ
c:\Sets^?_\idot \to \Pos$ is an enhanced category, and
$c:\Sets^?_\idot \to \Sets^?$ is an enhanced cofibration (if $?=h$,
then $\Sets^?_\idot = \Sets^h_\idot$, and $c$ is induced by the
enhanced cofibration \eqref{ccat.h.dot}).

Now consider the product $\Unf(\Pos) \times^h \sSets^?$, denote the
projections by $\pi_0:\Unf(\Pos) \times^h \sSets^? \to \Unf(\Pos)$
and $\pi_1:\Unf(\Pos) \times^h \sSets^? \to \sSets^?$, identify
$\Unf(\Pos) \cong \Arf(\Pos)$, and note that we have a functor
\begin{equation}\label{la.q.eq}
\lambda:\Unf(\Pos) \times^h \sSets^? \to \Sets^?
\end{equation}
sending a pair $\langle f,Z \rangle$ of a fibration $f:I' \to I$ in
$\Arf(\Pos) \cong \Unf(\Pos)$ and an object $Z \in \Sets^?(I)$ to
$f^{\iota*}\iota(Z) \in \Sets^?$. Then the product functor $\lambda
\times \pi_0$ has a right-adjoint
\begin{equation}\label{rho.q.eq}
\rho:\Unf(\Pos) \times \Sets^? \to \Unf(\Pos) \times^h
\sSets^?.
\end{equation}
Explicitly, \eqref{rho.q.eq} sends a pair $\langle f,Z \rangle$ of a
fibration $f:I' \to I$ in $\Arf(\Pos)$ and an object $Z \in \Sets^?$
to $\langle f,f^\iota_*(I' \times Z)\rangle$. In particular, for any
fixed $Z$, we have a functor $\rho(Z) = \pi_1 \circ
\rho|_Z:\Unf(\Pos) \to \sSets^?$ over $\Pos$. Assume that it is an
enhanced functor, and define an enhanced cofibration $\nu^c(Z) \to
\Unf(\Pos)$ by the semicartesian product
\begin{equation}\label{nu.q.eq}
\nu^c(Z) = \Unf(\Pos) \times^h_{\sSets^?} \sSets^?_\idot,
\end{equation}
where $c:\sSets^?_\idot \to \sSets^?$ is the enhanced cofibration of
\eqref{sets.q.sq}. Then in particular, $\nu^c(Z)_\ppt \to
\Unf(\Pos)_\ppt \cong \Pos$ is a cofibration; let $\nu(Z) \to \Pos$
be the transpose fibration, and assume that it is an enhanced
groupoid. By Corollary~\ref{semi.corr}, $\nu(Z)$ is functorial with
respect to $Z$ as an object in $\Sets^h$, so that we have a
comparison functor $\nu:\Sets^? \to \Sets^h$. Assuming that this
functor is an equivalence, we obtain an identification $\Sets^?
\cong \Sets^h$ and a full embedding $\Sets^h \subset \Cat^?$.

To extend this to a comparison functor $\Cat^? \to \Cat^h$, we
follow the same plan, but with one modification. Namely, the functor
\eqref{la.q.eq} can be composed with the full embedding $\Sets^?
\subset \Cat^?$, and then if we further assume that the embedding
$\Sets^?/I \to \Cat^?/I$ has a right-adjoint $R$ for any $I \in
\Pos$, $\lambda \times \pi_0$ still has a right-adjoint
\eqref{rho.q.eq} sending $\langle f,Z \rangle$ to $\langle
f,R(f^\iota_*(I' \times Z))\rangle$. However, morphisms in
$\Cat^?/I$ are all morphisms over $I$, and this is not want we want:
if $?=h$, then the map $f^\iota_*(I' \times Z) \to I$ is an
$I$-augmentation, and to get the right answer, we need to consider
only $I$-augmented maps. To correct for this, we consider the
decomposition \eqref{comma.facto} of the fibration $f:I' \to I$,
with the corresponding projections $\sigma:I \setminus I' \to I$,
$\tau:I \setminus I' \to I'$, embedding $\eta:I' \to I \setminus
I'$, and the map $\eta_\dg:I \setminus I' \to I'$ adjoint to
$\eta$. We then let $\wt{\eta}:(I \setminus I') \times [1] \to I'$
be the map defined by the adjunction map $\eta_\dg \to \tau$, and we
define $\rho'(Z) \in \Cat^?(I)$ by the cartesian square
$$
\begin{CD}
\rho'(Z) @>>> \sigma^\iota_*((I \setminus I') \times Z))\\
@VVV @VV{e^*}V\\
f^\iota_*(I' \times Z) @>{(\id \times \wt{\eta})^{\iota*}}>>
\Ar^?(\sigma^\iota_*((I \setminus I' \times Z)) \cong (\sigma
\times e)^\iota_*((I \setminus I') \times [1] \times Z).
\end{CD}
$$
To ensure that such a cartesian square exists in $\Cat^?$, we need
to further assume that for any $X \in \Cat^?$, with the map $\eta =
e^*:X \to \Ar^?(X)$, $\eta^*:\Cat^?/\Ar^?(X) \to \Cat^?/X$ exists
(if $? = h$, this holds by Corollary~\ref{semi.corr}, since $\eta$
is a fully faithful embedding). Then as for $\Sets^?$, we obtain a
functor $\rho'(Z):\Unf(\Pos) \to \sSets^?$ over $\Pos$ for any $Z \in
\Cat^?$, we assume that it is an enhanced functor, restrict the
induced enhanced cofibration \eqref{nu.q.eq} to $\Pos \cong
\Unf(\Pos)_\ppt$, and denote by $\nu(Z) \to \Pos$ the transpose
fibration. This is a family of groupoids; we then require that it is
of the form $\C_\flat$ for an enhanced category $\C$, unique and
functorial by Lemma~\ref{fl.le}, and this gives a functor $\Cat^?
\to \Cat^h$. It remains to assume that it is an equivalence.

\bigskip

{\noindent
Affiliations:
\begin{enumerate}
\renewcommand{\labelenumi}{\arabic{enumi}.}
\item Steklov Mathematics Institute (main affiliation).
\item National Research University Higher School of Economics.
\end{enumerate}}

{\noindent
{\em E-mail address\/}: {\tt kaledin@mi-ras.ru}
}

\end{document}